\documentclass[final,3p,times,fleqn]{elsarticle}
\usepackage{amssymb} 
\usepackage{amsmath} 
\usepackage{dcolumn}
\usepackage{bm}
\usepackage{booktabs}
\usepackage{amssymb}
\usepackage{amsmath}
\usepackage{array}
\usepackage{colortbl}
\usepackage{color} 
\usepackage{tabularx}
\usepackage{multirow} 
\usepackage{graphicx}
\usepackage{epstopdf}
\usepackage{subfig} 
\usepackage{bbding}
\usepackage{pifont}
\usepackage{wasysym}
\usepackage{utfsym}

\usepackage{xcolor}
\usepackage{algorithm}
\usepackage{algpseudocode} 

\begin{document}
	
	\begin{frontmatter}
		
		\title{A bound-preserving hybrid lattice Boltzmann method for scalar hyperbolic conservation laws}
		
		\author[a]{Ying Chen} 
		\author[a,b,c,d]{Zhenhua Chai \corref{cor1}}
		\ead{hustczh@hust.edu.cn}
		\address[a]{School of Mathematics and Statistics, Huazhong University of Science and Technology, Wuhan 430074, China}
		\address[b]{Institute of Interdisciplinary Research for Mathematics and Applied Science, Huazhong University of Science and Technology, Wuhan 430074, China}
		\address[c]{Hubei Key Laboratory of Engineering Modeling and
			Scientific Computing, Huazhong University of Science and Technology,
			Wuhan 430074, China}
		\address[d]{The State Key Laboratory of Intelligent Manufacturing Equipment and Technology, Huazhong University of Science and Technology, Wuhan 430074, China }
		\cortext[cor1]{Corresponding author.}
		\begin{abstract}  
	
		Solutions to scalar hyperbolic conservation laws may develop discontinuities within finite time even from smooth initial data, while higher-order lattice Boltzmann (LB) methods tend to generate non-physical oscillations and out-of-bounds behavior near steep gradients. To address this issue, we first construct a first-order LB method with the upwind-type equilibrium distribution functions under the \(DdQ(2d+1)\) lattice structure, and establish its bound-preserving (BP) conditions. We then develop a second-order multiple-relaxation-time LB (MRT-LB) method by introducing the auxiliary source distribution functions, and analyze its \(L^2\) stability. Building on these, we propose a hybrid LB method that preserves the collision-streaming structure, by conservatively mixing first- and second-order updates through a convex constraint mechanism. We further establish the global BP (GBP), local BP (LBP), and relaxed LBP (RLBP) limiters in the hybrid LB method, together with the associated macroscopic BP constraints. Numerical experiments examine the convergence order, interface capturing performance, and discontinuity resolution capability of the proposed BP hybrid LB method. The results indicate that it effectively suppresses out-of-bounds behavior in the second-order LB method and reduces numerical dissipation in the first-order LB method, while maintaining the prescribed bounds and discrete mass conservation.
			\end{abstract}
		
	\end{frontmatter}
	
	\section{Introduction}\label{introduction}
	 
	 Hyperbolic conservation laws arise in the modeling of nonlinear wave propagation \cite{LeVeque2002FiniteVolume}, compressible flows \cite{LandauLifshitz1987}, and interface transport \cite{OsherSethian1988}. Their solutions may develop discontinuities in finite time, even from smooth initial data, which poses challenges for the accuracy and stability of numerical methods. First-order monotone methods are generally robust but may smear discontinuities through excessive numerical dissipation. Higher-order methods provide better resolution in smooth regions but may generate spurious oscillations and violate solution bounds near discontinuities. Developing conservative methods that combine low numerical dissipation with bound preservation is therefore an important objective in the numerical approximation of hyperbolic conservation laws \cite{ZhangShu2011}.
	 
	 In this work, we consider the following $d$-dimensional scalar hyperbolic conservation law:
	 \begin{equation}\label{chapintro-eq}
	 	\partial_t u+\sum_{\alpha=1}^{d}\partial_{x_\alpha}F_\alpha(u)=0,
	 	\qquad
	 	u(\boldsymbol{x},0)=u_0(\boldsymbol{x}),
	 \end{equation}
	 where $\boldsymbol{x}\in\Omega\subseteq\mathbb{R}^d$, $u$ is the conserved variable, and $F_\alpha$ denotes the flux in the $x_\alpha$ direction. Under periodic boundary conditions, the entropy solution satisfies
	 \begin{align}\label{moM0}
	 m_0\leq u(\boldsymbol{x},t)\leq M_0,
	 \end{align}
	 almost everywhere whenever $m_0\leq u_0(\boldsymbol{x})\leq M_0$. A bound-preserving (BP) numerical method can maintain these macroscopic bounds at the discrete level and prevents nonphysical overshoots and undershoots. This property is particularly important for problems with shocks and transported interfaces.
	 
	 The lattice Boltzmann (LB) method provides a kinetic framework in which distribution functions evolve through collision and streaming, and macroscopic quantities are recovered from their moments \cite{kruger2017lattice,guo2013lattice,aidun2010lattice}. Its application to hyperbolic conservation laws is closely related to relaxation and discrete kinetic methods. Jin et al. \cite{JinXin1995} introduced relaxation approximations with linear transport and local relaxation, and Aregba-Driollet et al. \cite{AregbaDriolletNatalini2000} established convergence of multidimensional discrete kinetic approximations to entropy solutions for scalar problems. Together with the kinetic entropy conditions developed in \cite{Bouchut1999}, these results provide a theoretical basis for vectorial LB (VLB) methods.
	 
	The choice of equilibrium distribution functions in the LB method affects both bound preservation and numerical dissipation. A common choice is to express these functions as linear combinations of the conserved variable and its fluxes \cite{Chai2016,Chai2020} (hereafter referred to as linear equilibrium distribution functions). Their monotonicity imposes constraints on the weight coefficients and lattice velocity in terms of the flux derivatives. When the lattice velocity is determined by a global bound on the characteristic speeds, excessive numerical dissipation may arise in regions of slow local propagation. Upwind-type equilibrium distribution functions provide an alternative in which directional propagation is represented through the positive and negative parts of the flux derivatives. Anandan et al. \cite{AnandanRao2024} analyzed the \(H\)-inequality, total variation boundedness, and positivity of VLB methods with upwind-type equilibrium distribution functions. However, for scalar conservation laws, bound preservation requires the numerical solution to remain within the interval determined by the initial data. Positivity of the distribution functions alone does not ensure preservation of this interval, so the available positivity results do not directly establish the BP property. Moreover, positivity arguments based on a convex collision update do not directly extend to over-relaxation. Therefore, for the more general multiple-relaxation-time LB (MRT-LB) method \cite{Lallemand2000,Pan2006,Cui2016,Luo2011}, a BP analysis must establish bounds on the distribution functions that imply the prescribed macroscopic bounds, and find constraints on the relaxation parameters under which the collision step preserves these bounds.

	 Higher-order LB methods offer improved accuracy. Recently, Bellotti \cite{Bellotti2024Parasitic} developed a fourth-order LB method with the relaxation parameter fixed at \(s=2\) for the one-dimensional linear transport equation, and Chen et al. \cite{chen2026fourth-order} further extended this study to multidimensional cases. Although the characteristic polynomial of the amplification matrix need not be a simple von Neumann one, these methods can remain stable for linear problems under suitable conditions.  However, such linear stability, does not guarantee nonlinear stability at \(s=2\). For nonlinear hyperbolic problems, Bellotti et al. \cite{BellottiHelluyNavoret2025} constructed a
	 fourth-order LB method through the composition of time-symmetric second-order operators initially defined with $s=2$, and introduced
	 a local adjustment of the relaxation parameter to conserve microscopic entropy. However, this entropy control alone does not
	 guarantee pointwise bounds. Zhao et al. \cite{Zhao2026MRT} developed an MRT-LB method for the special nonlinear advection equations with reduced oscillations, but they did not establish the admissible ranges of the relaxation parameters and the lattice velocity that guarantee bound preservation.

	 Hybrid methods combine a BP low-order update with high-order corrections that are restricted by a  limiter. Recently, Wissocq et al. \cite{WissocqLiuAbgrall2025} developed a conservative hybrid VLB method in which convex limiters enforce positivity and local bounds for the two-dimensional Euler equations. Their work provides a basis for hybrid LB methods. However, the BP constraints associated with the MRT collision operator considered in this work require a further analysis. The performance of a hybrid method also depends on its low-order and high-order components. The dissipation in the low-order method affects the resolution in the regions where the limiter is applied, whereas the stability of the high-order method influences the corrections that the limiter must control. These considerations motivate a joint study of bound preservation, numerical dissipation, and linear stability.
	
	In this work, we will develop a conservative hybrid LB method for the \(d\)-dimensional scalar hyperbolic conservation law \eqref{chapintro-eq}. It combines a BP first-order LB method based on upwind-type equilibrium distribution functions with a second-order MRT-LB method that incorporates an auxiliary source term. Global or local BP limiters control the high-order corrections to preserve the prescribed bounds with reduced numerical dissipation, while the method can retain the collision-streaming structure.
	
	The remainder of this paper is organized as follows. Section~\ref{chap1} first presents the first-order LB method under the $DdQ(2d+1)$ lattice structure and establishes its BP conditions, and then introduces the second-order MRT-LB method and analyzes its $L^2$ stability for linear equations. Section~\ref{chap3} develops the conservative hybrid LB method with the global BP (GBP), local BP (LBP), and relaxed LBP (RLBP) limiters. Section~\ref{chap4} presents numerical tests involving Burgers' equation, multidimensional interface transport, and rotating waves with nonconvex fluxes. Finally, Section~\ref{chap5} summarizes the results and discusses future work.

	\section{Bound-preserving first-order lattice Boltzamnn method}\label{chap1}
	In this section, for the \(d\)-dimensional scalar hyperbolic conservation law \eqref{chapintro-eq}, we adopt the \(DdQ(2d+1)\) lattice structure and upwind-type equilibrium distribution functions. First, we present the first-order LB method, analyze its sufficient conditions for the BP property, and compare the BP constraints based on the upwind-type and the classical linear equilibrium distribution functions. We then develop a second-order MRT-LB method with the auxiliary source distribution functions, followed by an \(L^2\) stability analysis.

	\subsection{First-order lattice Boltzamnn method}\label{chapfirst}
	\subsubsection{Evolution equation of first-order lattice Boltzamnn method}\label{chapfirst-lb}
	For \eqref{chapintro-eq}, the classical first-order LB method with the general MRT collision model can be written as \cite{Chai2020}
	\begin{align}\label{chapfirst-mrt-lbm}
		f_i(\mathbf{x}+\mathbf{c}_i\Delta t,t+\Delta t)=f_i(\mathbf{x},t)-\bm{\Lambda}_{ij}\big[f_j-f_j^{\text{eq}}\big](\mathbf{x},t),\qquad i=0,1,\ldots,2d,
	\end{align}
	where $f_i$ and $f_i^{\text{eq}}$ denote the distribution functions and equilibrium distribution functions, respectively. The MRT collision matrix \(\bm{\Lambda}\) is given by \(\bm{\Lambda} = \mathbf{M}^{-1}\mathbf{S}\mathbf{M}\), where \(\mathbf{M}\) is the transformation matrix and \(\mathbf{S}\) is the relaxation matrix. \(\Delta x\) is the uniform lattice spacing on the square grid employed in this work, and \(\Delta t\) is the time step.
	
	In this study, we adopt the \(DdQ(2d+1)\) lattice structure, and the discrete velocity set is expressed as
	\begin{align*}
		\mathbf{c}_0 = (0, 0, \ldots, 0)^T\in\mathbb{R}^{d\times 1}, \qquad \mathbf{c}_{\alpha} = (0, 0, \ldots, \underbrace{c}_{i\text{-th entry}}, \dots, 0)^T\in\mathbb{R}^{d\times 1}, \qquad \alpha = 1, 2, \ldots, d.
	\end{align*}
	where \(c=\Delta x/\Delta t\) is the lattice velocity.
	Based on this lattice structure, we adopt the collision model in a natural moment space \cite{Chai2020}, in which the transformation matrix \(\mathbf{M}\) reads
	\begin{align}\label{chapfirst-m}
		\mathbf{M} = \left( \mathbf{1}, \hat{\mathbf{c}}_{1}, \hat{\mathbf{c}}_{2}, \dots, \hat{\mathbf{c}}_{d}, \hat{\mathbf{c}}_{1}^{.2}, \hat{\mathbf{c}}_{2}^{.2}, \dots, \hat{\mathbf{c}}_{d}^{.2} \right)^T,
	\end{align}
	where  \(\mathbf{1}\in\mathbb{R}^{(2d+1)\times 1}\) is the all-ones column vector,    \(\hat{\mathbf{c}}_{\alpha}\) is given by
	\begin{align*}
		\hat{\mathbf{c}}_{\alpha}=(c_{0\alpha},c_{1\alpha},\ldots,c_{d\alpha},-c_{1\alpha},\ldots,-c_{d\alpha})^T\in\mathbb{R}^{(2d+1)\times 1},\qquad \alpha=1,2,\ldots,d.
	\end{align*} 
	and  \(\hat{\mathbf{c}}_{\alpha}^{.2}\in\mathbb{R}^{(2d+1)\times 1}\) is the column vector whose \(\beta\)-th entry is \(\hat{c}_{\alpha\beta}^2\) for \(\beta = 1, 2, \dots, (2d+1)\).

	For the relaxation matrix \(\mathbf{S}\), we consider the following diagonal form for brevity:
	\begin{align}\label{chapfirst-s}
		\mathbf{S} = \operatorname{diag}\left( s_0, s_{x_1}, s_{x_2}, \dots, s_{x_d}, p_{x_1}, p_{x_2}, \dots, p_{x_d} \right),
	\end{align}
	where \(s_0\) corresponds to the conserved moment, i.e., \(u\) in \eqref{chapintro-eq}, \(s_{x_\alpha}\) and \(p_{x_\alpha}\) correspond to the first- and second-order non-conserved moments of the distribution functions, respectively.

	For \(f_i^{\text{eq}}\) in \eqref{chapfirst-mrt-lbm}, the most widely adopted forms in the LB method are the classical linear equilibrium distribution functions \cite{Chai2016,Chai2020}, i.e.,
		\begin{align}\label{chapfirst-linear-feq}
		f_i^{\text{eq}}=\omega_i\left(u+\frac{\mathbf{c}_i\cdot\mathbf{F}(u)}{c_s^2}\right),\qquad i=0,1,\ldots,2d,
	\end{align}
	where \(\omega_i\) are the weight coefficients satisfying the symmetry condition \(\omega_\alpha = \omega_{\alpha+d}\) (\(\alpha = 1, 2, \dots, d\)), and \(c_s\) is the lattice sound speed determined by  \(c\) and \(\omega_i\).
	
	The choice of equilibrium distribution functions affects both the admissible parameter range and the numerical dissipation of LB methods. This motivates the use of upwind-type equilibrium distribution functions based on the characteristic wave propagation directions   \cite{AnandanRao2024} in this work, which are given by
	\begin{align}\label{chapfirst-upwind-feq}
		\begin{cases}
			f_0^{\text{eq}}(u)=u-\sum_{\alpha=1}^{d}\frac{F_\alpha^+(u)+F_{\alpha}^-(u)}{c},\\
			f_{\alpha}^{\text{eq}}(u)=\frac{F_{\alpha}^{+}(u)}{c},\qquad f_{\alpha+d}^{\text{eq}}(u)=\frac{F_{\alpha}^{-}(u)}{c},\qquad \alpha=1,2,\ldots,d,
		\end{cases}
	\end{align}
	where  \(F_i^+(u)\) and \(F_i^-(u)\) are the forward and backward upwind flux components along the \(i\)-th direction, satisfying \(F_i(u) = F^+_i(u) - F^-_i(u)\).
	
	Standard flux splitting schemes include the Lax-Friedrichs \cite{Lax1954}, Steger-Warming \cite{StegerWarming1981}, and van Leer \cite{vanLeer1974} schemes, as well as kinematic flux vector splitting (KFVS) \cite{MandalDeshpande1994}. Engquist-Osher flux splitting \cite{EngquistOsher1981} is also applicable to scalar conservation laws. Different schemes feature different wave-speed treatments and numerical dissipation mechanisms, and their stability and structure-preserving properties must be analyzed together with sub-characteristic conditions and relaxation parameters.
	
	In this work, we adopt the Engquist-Osher flux splitting, with \(F_i^\pm(u)\) given by
	\begin{align*}
		F_i^\pm(u)=\int_0^u[\pm F_i'(s)]_+\mathrm ds,
	\end{align*}
with \([z]_+ = \max\{z, 0\}\).	Under this construction, it is obvious that the above upwind-type equilibrium distribution functions \eqref{chapfirst-upwind-feq} satisfy the follwing moment conditions:
	\begin{align}\label{chapfirst-moment}
		\sum_{i=0}^{2d}f_i^{\text{eq}}=u,\qquad \sum_{i=0}^{2d} \mathbf{c}_i f_i^{\text{eq}}=\mathbf{F}(u).
	\end{align}
	The macroscopic quantity \(u\) at the lattice node \(K\) and time \(t_n = n \Delta t\), denoted \(u_K^n\) for brevity, is then given by
	\[
	u_K^n = \sum_{i=0}^{2d} f_{i,K}^n.
	\]
	
	By direct Taylor expansion, at the acoustic scaling (\(\Delta x\propto \Delta t\)), it is easy to verify that the above LB method \eqref{chapfirst-mrt-lbm} with \(f_i^{\text{eq}}\) in \eqref{chapfirst-upwind-feq} is first-order consistent with the scalar hyperbolic conservation law \eqref{chapintro-eq} (see Section \ref{chapsecond} for details). The corresponding modified equation is given by
	\begin{align}\label{chapfirst-me}
		\partial_tu+\sum_{\alpha=1}^d\partial_{x_\alpha}F_{\alpha}(u)=\Delta t \sum_{\alpha,\beta=1}^d\partial_{x_\alpha}\left[\Big(\frac{1}{s_{x_\alpha}}-\frac{1}{2}\Big)\big(c\delta_{\alpha\beta}G_{\alpha}'(u)-F_{\alpha}(u)'F_{\beta}'(u)\big)\partial_{\beta}u\right]+O(\Delta t^2).
	\end{align}
	with \(G_\alpha=F_\alpha^++F_\alpha^-\).
	\subsubsection{Bound preservation for first-order lattice Boltzmann method}
	
	The LB method consists of collision and streaming steps. The collision step preserves the local conserved quantity $u$,
	whereas the streaming step reconstructs $u$ from distribution functions at different lattice nodes. Thus, bounds on $u$ before streaming do not automatically imply the same bounds after streaming. To establish the BP property, we first rewrite the first-order MRT-LB method \eqref{chapfirst-mrt-lbm} as
	\begin{align*}
		\begin{cases}
			\text{Collision: }
			f_i^\star(\mathbf x,t)
			=\big[\mathcal C(\mathbf f\big]_i(\mathbf x,t)),\\
			\text{Streaming: }
			f_i(\mathbf x+\mathbf c_i\Delta t,t+\Delta t)
			=f_i^\star(\mathbf x,t),
		\end{cases}
	\end{align*}
	where \(i=0,1,\ldots,2d\) and the collision operator \(\mathcal C\) is defined by
	\[\mathcal C(\mathbf f)
	=\mathbf f-\bm\Lambda\bigl[\mathbf f-\mathbf f^{\mathrm{eq}}(u)\bigr].
	\]
	Here, \(\mathbf f=(f_0,f_1,\ldots,f_{2d})^T\) and \(\mathbf f^{\mathrm{eq}}(u)	=(f_0^{\mathrm{eq}}(u),f_1^{\mathrm{eq}}(u),\ldots,f_{2d}^{\mathrm{eq}}(u))^T\).
	
	Our analysis proceeds by constructing componentwise bounds for the distribution functions and deriving conditions under
	which these bounds are preserved in collision step. Since streaming step only shifts each distribution component between lattice nodes, which means that constant componentwise bounds are also preserved under periodic boundary conditions. From \eqref{moM0}, we define the following candidate bounds:
	\begin{equation}\label{chapfirst-up-low}
		\underline{\mathbf f}=\mathbf f^{\mathrm{eq}}(m_0),
		\qquad
		\overline{\mathbf f}=\mathbf f^{\mathrm{eq}}(M_0).
	\end{equation}
	Their componentwise ordering will follow from the monotonicity
	of \(\mathbf f^{\mathrm{eq}}\) established below.
	
	We next consider the order-preserving property of  collision operator $\mathcal C$. Since $u=\sum_{i=0}^{2d}f_i$, we have $\partial u/\partial f_i=1$. For a constant relaxation matrix $\bm\Lambda$, the Jacobian of  $\mathcal C$ is expressed as
	\[	\text{D}\mathcal C(\mathbf f)	=\mathbf I-\bm\Lambda	+\bm\Lambda	\frac{\text{d}\mathbf f^{\mathrm{eq}}(u)}{\text{d}u}\mathbf 1^T.
	\]
	Equivalently,
	\begin{equation*}
		\frac{\partial\mathcal C_i}{\partial f_j}
		=\delta_{ij}-\Lambda_{ij}
		+\sum_{k=0}^{2d}\Lambda_{ik}
		\frac{\text{d} f_k^{\mathrm{eq}}(u)}{\text{d}u},
		\qquad i,j=0,1,\ldots,2d.
	\end{equation*}
	Here, $\delta_{ij}$ is the Kronecker delta. A sufficient condition for $\mathcal C$ to be order-preserving
	on a convex set is that its Jacobian be componentwise non-negative throughout that set. Indeed, for any two vectors
	$\mathbf f\le\mathbf g$ in the set, one can obtain
	\begin{equation}\label{chapfirst-order-preserve}
		\mathcal C(\mathbf g)-\mathcal C(\mathbf f)
		=
		\int_0^1
		\text{D}\mathcal C\bigl(\mathbf f+\theta(\mathbf g-\mathbf f)\bigr)
		(\mathbf g-\mathbf f)\,d\theta
		\ge\mathbf 0.
	\end{equation}
	 In the following theorem, we derive conditions ensuring this
	property for states satisfying $u\in[m_0,M_0]$, and further utilize it to establish bound preservation of the first-order LB method 	\eqref{chapfirst-mrt-lbm}. \\
\noindent\textbf{Theorem 1 (Bound preservation).}
Assume that the lattice velocity $c>0$, the relaxation parameters satisfy $s_{x_\alpha}>0$ and $p_{x_\alpha}>0$, and $F_\alpha,G_\alpha\in C^1([m_0,M_0])$ for $\alpha=1,2,\ldots,d$. Suppose that periodic boundary conditions are imposed and the initial data satisfy \(m_0\le u_0(\mathbf x)\le M_0\) at every lattice node \(\mathbf{x}\). Initialize the distribution functions by
\(f_i^0(\mathbf x)=f_i^{\mathrm{eq}}\bigl(u_0(\mathbf x)\bigr)\), \(i=0,1,\ldots,2d\). If, for every $u\in[m_0,M_0]$,
\begin{equation}\label{chapfirst-monoton}
	\begin{cases} 
		\sum_{\alpha=1}^{d}p_{x_\alpha}G_\alpha'(u)
		\le
		c\min\{1,p_{x_1},\ldots,p_{x_d}\},
		\\
		p_{x_\alpha}G_\alpha'(u)
		-s_{x_\alpha}|F_\alpha'(u)|
		\ge
		c\max\{0,\,
		p_{x_\alpha}+s_{x_\alpha}-2,\,
		p_{x_\alpha}-s_{x_\alpha}\},
		\quad \alpha=1,2,\ldots,d,
	\end{cases}
\end{equation}
 then the numerical solution generated by \eqref{chapfirst-mrt-lbm} satisfies \(m_0\le u^n(\mathbf x)\le M_0\) at every lattice node and for all $n\ge0$.

\medskip
\noindent\textbf{Proof.}
We first express the collision step in terms of the moments
\begin{equation*}
	u=\sum_{i=0}^{2d}f_i,\qquad
	j_\alpha=c(f_\alpha-f_{\alpha+d}),\qquad
	k_\alpha=c^2(f_\alpha+f_{\alpha+d}),
\end{equation*}
where $\alpha=1,2,\ldots,d$. The corresponding inverse relations are
\[
f_0=u-\frac{1}{c^2}\sum_{\alpha=1}^{d}k_\alpha,
\qquad
f_\alpha=\frac{k_\alpha}{2c^2}+\frac{j_\alpha}{2c},
\qquad
f_{\alpha+d}=\frac{k_\alpha}{2c^2}-\frac{j_\alpha}{2c}.
\]
Since $j_\alpha^{\mathrm{eq}}=F_\alpha(u)$ and $k_\alpha^{\mathrm{eq}}=cG_\alpha(u)$, the post-collision moments satisfy
\[
\begin{aligned}
	u^\star=u,\qquad 	j_\alpha^\star=(1-s_{x_\alpha})j_\alpha+s_{x_\alpha}F_\alpha(u),\qquad
	k_\alpha^\star=(1-p_{x_\alpha})k_\alpha+p_{x_\alpha}cG_\alpha(u).
\end{aligned}
\]
Substituting these expressions into the inverse moment relations gives
\begin{align}\label{chapfirst-collision-components}
	\begin{cases}
		\mathcal C_0(\mathbf f)
		=f_0+\sum_{\alpha=1}^{d}p_{x_\alpha}(f_\alpha+f_{\alpha+d})
		-\frac{\sum_{\alpha=1}^{d}p_{x_\alpha}G_\alpha(u)}{c},\\[1mm]
		\mathcal C_\alpha(\mathbf f)
		=\left(1-\frac{s_{x_\alpha}+p_{x_\alpha}}{2}\right)f_\alpha
		+\frac{s_{x_\alpha}-p_{x_\alpha}}{2}f_{\alpha+d}+\frac{p_{x_\alpha}G_\alpha(u)+s_{x_\alpha}F_\alpha(u)}{2c},\\[1mm]
		\mathcal C_{\alpha+d}(\mathbf f)=\frac{s_{x_\alpha}-p_{x_\alpha}}{2}f_\alpha
		+\left(1-\frac{s_{x_\alpha}+p_{x_\alpha}}{2}\right)f_{\alpha+d}+\frac{p_{x_\alpha}G_\alpha(u)-s_{x_\alpha}F_\alpha(u)}{2c},
	\end{cases}
\end{align}
where \(\alpha=1,2,\ldots,d\). In particular, $\sum_{i=0}^{2d}\mathcal C_i(\mathbf f)=u$.

To examine the signs of the Jacobian entries, we define
\[
\begin{aligned}
	H(u)=\frac{1}{c}\sum_{\alpha=1}^{d}p_{x_\alpha}G_\alpha'(u),\qquad 	B_\alpha(u)=\frac{p_{x_\alpha}G_\alpha'(u)+s_{x_\alpha}F_\alpha'(u)}{2c},
	\qquad
	B_{\alpha+d}(u)=\frac{p_{x_\alpha}G_\alpha'(u)-s_{x_\alpha}F_\alpha'(u)}{2c}.
\end{aligned}
\]
Differentiating \eqref{chapfirst-collision-components} yields
\[
\frac{\partial\mathcal C_0}{\partial f_0}=1-H(u),
\qquad
\frac{\partial\mathcal C_0}{\partial f_\beta}
=\frac{\partial\mathcal C_0}{\partial f_{\beta+d}}
=p_{x_\beta}-H(u),
\quad \beta=1,2,\ldots,d,
\]
and
\[
\begin{aligned}
	&\frac{\partial\mathcal C_\alpha}{\partial f_j}
	=\left(1-\frac{s_{x_\alpha}+p_{x_\alpha}}{2}\right)\delta_{j\alpha}
	+\frac{s_{x_\alpha}-p_{x_\alpha}}{2}\delta_{j\alpha+d}+B_\alpha(u),\\
	&\frac{\partial\mathcal C_{\alpha+d}}{\partial f_j}
	=\frac{s_{x_\alpha}-p_{x_\alpha}}{2}\delta_{j\alpha}
	+\left(1-\frac{s_{x_\alpha}+p_{x_\alpha}}{2}\right)\delta_{j\alpha+d}
	+B_{\alpha+d}(u),
\end{aligned}
\]
where \(j=1,2,\ldots,2d\) and \(\alpha=1,2,\ldots,d\). Conditions \eqref{chapfirst-monoton} imply
\begin{align*}
	\begin{cases}
		H(u)\le\min\{1,p_{x_1},\ldots,p_{x_d}\},\\
		\min\{B_\alpha(u),B_{\alpha+d}(u)\}
		\ge\max\{0,p_{x_\alpha}+s_{x_\alpha}-2,
		p_{x_\alpha}-s_{x_\alpha}\}/2.
	\end{cases}
\end{align*}
Thus, $\partial\mathcal C_i/\partial f_j\ge0$ for all $i,j=0,1,\ldots,2d$ whenever $u=\sum_i f_i\in[m_0,M_0]$.

We next verify that the equilibrium distribution functions \(f_i^{\text{eq}}\) are nondecreasing on $[m_0,M_0]$. The second inequality in \eqref{chapfirst-monoton} gives
\[
p_{x_\alpha}G_\alpha'(u)\ge s_{x_\alpha}|F_\alpha'(u)|\ge0,
\]
and hence $G_\alpha'(u)\ge0$. Let $p_{\min}=\min_{\alpha=1,2,\ldots,d}p_{x_\alpha}>0$. The first inequality in \eqref{chapfirst-monoton} then yields
\[
p_{\min}\sum_{\alpha=1}^{d}G_\alpha'(u)
\le\sum_{\alpha=1}^{d}p_{x_\alpha}G_\alpha'(u)
\le c\min\{1,p_{\min}\}\le cp_{\min}.
\]
Consequently, \(\sum_{\alpha=1}^{d}G_\alpha'(u)\le c\) and \(0\le G_\alpha'(u)\le c\). In particular, if $p_{x_\alpha}\le s_{x_\alpha}$, we have
\[
G_\alpha'(u)\ge\frac{s_{x_\alpha}}{p_{x_\alpha}}|F_\alpha'(u)|
\ge|F_\alpha'(u)|.
\]
If $p_{x_\alpha}>s_{x_\alpha}$, the second inequality in \eqref{chapfirst-monoton} gives
\[
s_{x_\alpha}\bigl(G_\alpha'(u)-|F_\alpha'(u)|\bigr)
\ge(p_{x_\alpha}-s_{x_\alpha})\bigl(c-G_\alpha'(u)\bigr)
\ge0.
\]
Therefore, $G_\alpha'(u)\ge|F_\alpha'(u)|$ in both cases, and
	\begin{align}\label{chapfirst-feq-mono}
		\begin{cases}
				\frac{\text{d} f_0^{\mathrm{eq}}}{\text{d}u}
			=1-\frac1c\sum_{\alpha=1}^{d}G_\alpha'(u)\ge0,\\
			\frac{\text{d} f_\alpha^{\mathrm{eq}}}{\text{d}u}
			=\frac{G_\alpha'(u)+F_\alpha'(u)}{2c}\ge0,\qquad \frac{\text{d} f_{\alpha+d}^{\mathrm{eq}}}{\text{d}u}
			=\frac{G_\alpha'(u)-F_\alpha'(u)}{2c}\ge0,\qquad \alpha=1,2,\ldots,d.
		\end{cases} 
	\end{align} 

Using the bounds in \eqref{chapfirst-up-low}, we define \(\mathcal K=\{\mathbf f\in\mathbb R^{2d+1}:
\underline{\mathbf f}\le\mathbf f\le\overline{\mathbf f}\}\). By \eqref{chapfirst-feq-mono} and the equilibrium initialization, $\mathbf f^0(\mathbf x)\in\mathcal K$. Moreover, every $\mathbf f\in\mathcal K$ satisfies
\[
m_0=\sum_{i=0}^{2d}\underline f_i
\le\sum_{i=0}^{2d}f_i
\le\sum_{i=0}^{2d}\overline f_i=M_0.
\]
The Jacobian of $\mathcal C$ is therefore componentwise non-negative throughout the convex set $\mathcal K$. It follows from \eqref{chapfirst-order-preserve} that $\mathcal C$ is order-preserving on $\mathcal K$.

Since equilibrium distribution functions are fixed nodes of the collision operator, i.e.,
\[
\mathcal C\bigl(\mathbf f^{\mathrm{eq}}(\kappa)\bigr)
=\mathbf f^{\mathrm{eq}}(\kappa),
\qquad \kappa\in\{m_0,M_0\},
\]
we can obtain, for every $\mathbf f\in\mathcal K$,
\[
\underline{\mathbf f}
=\mathcal C(\underline{\mathbf f})
\le\mathcal C(\mathbf f)
\le\mathcal C(\overline{\mathbf f})
=\overline{\mathbf f}.
\]
Thus, the collision step preserves $\mathcal K$. Under periodic boundary conditions, the streaming step \(f_i^{n+1}(\mathbf x)=f_i^{\star,n}(\mathbf x-\mathbf c_i\Delta t)\) also preserves the constant bounds $\underline f_i\le f_i\le\overline f_i$. By induction, $\mathbf f^n(\mathbf x)\in\mathcal K$ at every lattice node and for all $n\ge0$. Summing over the velocity index gives
\[
m_0=\sum_{i=0}^{2d}\underline f_i
\le u^n(\mathbf x)=\sum_{i=0}^{2d}f_i^n(\mathbf x)
\le\sum_{i=0}^{2d}\overline f_i=M_0,
\]
which completes the proof.\hfill$\square$

	Now we give some discussions on the BP conditions in \eqref{chapfirst-monoton}.
	\begin{itemize}
		\item We would like to point out that, for our proposed first-order LB method in Section \ref{chapfirst-lb}, the BP conditions in \eqref{chapfirst-monoton} require that all relaxation parameters must be no greater than 1. Specifically, under the upwind-type equilibrium distribution functions, we have \(G_\alpha=F_\alpha^++F_\alpha^-\) which satisfies \(G_\alpha'(u)=|F_\alpha(u)|\), then \eqref{chapfirst-monoton} becomes	
	\begin{subequations}
	\begin{align}
			& \sum_{\alpha=1}^dp_{x_{\alpha}}|F_{\alpha}'(u)|\leq c\min\{1,p_{x_1},p_{x_2},\ldots,p_{x_d}\},\label{chapfirst-con1}\\
			&(p_{x_\alpha}-s_{x_{\alpha}})|F_{\alpha}'(u)|\geq c\max\{0,p_{x_{\alpha}}+s_{x_{\alpha}}-2,p_{x_{\alpha}}-s_{x_{\alpha}}\},\qquad\alpha=1,2,\ldots,d.\label{chapfirst-con2}
		\end{align}
		\end{subequations}	
		which is equivalent to
		\begin{align*}
			\begin{cases}
				\sum_{\alpha=1}^ds_{x_\alpha}|F_{\alpha}'(u)|\leq c\min\{1,s_{x_1},s_{x_1},\ldots,,s_{x_d}\}&\text{from \eqref{chapfirst-con1}},\\
				0<p_{x_\alpha}=s_{x_{\alpha}}\leq 1, \qquad \alpha=1,2,\ldots,d,&\text{from \eqref{chapfirst-con2}}.
			\end{cases}
		\end{align*}	
		In particular, the monotonic and nondecreasing conditions for \(f_i^{\text{eq}}\) in \eqref{chapfirst-feq-mono} is essentially		
		\begin{align}\label{chapfirst-feq-mono-reduce}
			\sum_{\alpha=1}^d|F_{\alpha}'(u)|\leq c.
		\end{align}	
		
		What is more, the present LB method is equivalent to a \((2d+2)\)-level macroscopic finite-difference scheme in the sense of the conservative moment \cite{bellotti2022finite,chen2024gplb,bellotti2023truncation}. And when all the relaxation parameters are set to 1, the corresponding macroscopic finite-difference scheme is a two-level one, i.e.,
		\begin{align*}
			u_K^{n+1}=u_K^n-\frac{\Delta t}{\Delta x}\sum_{\alpha=1}^d\Big[F_{\alpha}^+(u_K^n)-F_{\alpha}^+(u_{K-\mathbf{e}_\alpha}^n)
			+F_{\alpha}^-(u_{K}^n)-F_{\alpha}^-(u_{K+\mathbf{e}_\alpha}^n)\Big],
		\end{align*}
		where \(K\pm\mathbf{e}_i\) (\(\mathbf{e}_i = \mathbf{c}_i / c\)) denotes the positive and negative adjacent nodes of \(K\). And it is obvious that the above equation is the classic first-order Engquist-Osher upwind scheme \cite{engquist1981one,crandall1980monotone,vanleer1984relation}. 
		\item When the MRT collision model reduces to a TRT model, the BP conditions \eqref{chapfirst-monoton} simplify to 
		\begin{align*}
			\begin{cases}
				\sum_{\alpha=1}^ds_{e}G'_{\alpha}(u)\leq \min c\{1,s_{e}\},\\
				s_{e}G_{\alpha}'(u)-s_{\alpha}|F_{\alpha}'(u)|\geq c\max\{0,s_{e}+s_{o}-2,s_{e}-s_{o}\},
			\end{cases}
		\end{align*}	
		where  \(s_o\) and \(s_e\) are the relaxation parameters corresponding the odd- and even-order moments of distribution functions, respectively. In particular, when considering the linear equilibrium distribution functions \eqref{chapfirst-linear-feq}, the BP conditions are consistent with that in Ref. \cite{AregbaDriolletBellotti2026}. 
		\item If the equilibrium distribution functions are the linear forms shown in \eqref{chapfirst-linear-feq},   the BP conditions \eqref{chapfirst-monoton} become
		\begin{align*}
			\begin{cases}
				2\sum_{\alpha=1}^d\omega_{\alpha}p_{x_{\alpha}}\leq \min\{1,p_{x_1},p_{x_2},\ldots,p_{x_d}\},\\
				2cp_{x_\alpha}\omega_{\alpha}-s_{x_{\alpha}}|F_{\alpha}'(u)|\geq c\max\{0,p_{x_{\alpha}}+s_{x_{\alpha}}-2,p_{x_{\alpha}}-s_{x_{\alpha}}\},\qquad\alpha=1,2,\ldots,d.
			\end{cases}
		\end{align*}	
		And the monotonic and nondecreasing conditions for \(f_i^{\text{eq}}\) in \eqref{chapfirst-feq-mono} become
		\begin{align}\label{chapfirst-feq-upwind-monot}
			\begin{cases}
				0\leq\omega_0<1,\\
				|F_{\alpha}'(u)|\leq 2\omega_{\alpha}c,\qquad\alpha=1,2,\ldots,d,
			\end{cases}
		\end{align}	
		from which one can obtain \(\sum_{\alpha=1}^d|F_{\alpha}'(u)|\leq c\).  This indicates that the conditions for monotonicity of the linear \(f_i^{\text{eq}}\) are stricter than those for the upwind-type \(f_i^{\text{eq}}\) [see \eqref{chapfirst-feq-mono-reduce}], which further suggests that the upwind-type \(f_i^{\text{eq}}\) offers superior stability [see the numerical results in Sections \ref{chapnum-ex1} and \ref{chapnum-ex3} for details].
	\end{itemize}
	\subsection{Second-order lattice Boltzamnn method}\label{chap2}
	\subsubsection{Second-order lattice Boltzmann method with a auxiliary source term}\label{chapsecond}
	According to the modified equation \eqref{chapfirst-me}, the scalar hyperbolic conservation laws \eqref{chapintro-eq} can be rewritten as
	\begin{align}\label{chapsecond-eq}
		\partial_tu+\nabla\cdot\mathbf{F}(u)= \nabla\cdot\left[\mathbf{K}\mathbf{A}(u;c)\nabla u\right]+(-\nabla\cdot\left[\mathbf{K}\mathbf{A}(u;c)\nabla u\right]),
	\end{align} 
	with
	\begin{align*}
		\mathbf{K}=\text{diag}(\kappa_1,\kappa_2,\ldots,\kappa_d),\qquad 
		\mathbf{A}(u;c)=\text{diag}\big(c\mathbf{G}'(u)\big)-\mathbf{F}'(u)\mathbf{F}'(u)^T.
	\end{align*}
	It is readily seen that \eqref{chapsecond-eq} is essentially a convection-diffusion equation with diffusion matrix \(\mathbf{K}\in\mathbb{R}^{d\times d}\) and a source term of the form \(-\nabla\cdot\left[\mathbf{K}\mathbf{A}(u;c)\nabla u\right]\). Then the corresponding MRT-LB method can be constructed as \cite{Chai2020}
	\begin{align}\label{chapsecond-mrt-lbm}
		f_i(\mathbf{x}+\mathbf{c}_i\Delta t,t+\Delta t)=f_i(\mathbf{x},t)-\bm{\Lambda}_{ij}\big[f_j-f_j^{\text{eq}}\big](\mathbf{x},t)+\Delta t\left(\mathbf{I}-\frac{\bm{\Lambda}}{2}\right)_{ij}R_j(\mathbf{x},t),\qquad i=0,1,\cdots,2d,
	\end{align}
	where the lattice structure, collision matrix, and equilibrium distribution functions are identical to those in the first-order case presented in Section \ref{chapfirst-lb}. The source distribution functions \(R_i\) are designed as
			\begin{align}\label{chapsecond-source}
				\begin{cases}
			R_0=-\left[1-\frac{1}{c}\sum_{\alpha=1}^d  \big((F_{\alpha}^+)'+(F_{\alpha}^-)'\big) \right]\mathcal{A}(u),\\
			R_{\alpha}=\partial_{x_\alpha}F_{\alpha}^+-\frac{(F_{\alpha}^+)'}{c}\mathcal{A}(u),\qquad 			R_{\alpha+d}=-\partial_{x_\alpha}F_{\alpha}^--\frac{(F_{\alpha}^-)'}{c}\mathcal{A}(u), \qquad \alpha=1,2,\ldots,d,
		\end{cases}			
		\end{align}
	where 
	\begin{align}\label{chapsecond-au}
		\mathcal A(u)=\sum_{\beta=1}^d\partial_{x_\beta}F_\beta(u)=\sum_{\beta=1}^dF_\beta'(u)\partial_{x_\beta}u.
	\end{align}
	Under this construction, it is obvious that the source distribution functions satisfy the following moment conditions:
	\begin{align}\label{chapsecond-source-moment}
		\sum_{i=0}^{2d}R_i=0,\qquad \sum_{i=0}^{2d}c_{i\alpha} R_i=cG_{\alpha}'(u)\partial_{\alpha}u-F_{\alpha}'(u)\mathcal{A}(u).
	\end{align}

	Under the acoustic scaling, we perform an accuracy analysis of \eqref{chapsecond-mrt-lbm} via direct Taylor expansion \cite{Chai2020}. For notational simplicity, we write \(\psi(\mathbf{x}, t)\) simply as \(\psi\). Expanding \eqref{chapsecond-mrt-lbm} in a Taylor series up to \(O(\Delta t^3)\) yields
	\begin{align*}
		\Delta t\left(\partial_tf_i+\mathbf{c}_i\cdot \nabla f_i\right)+\frac{\Delta t^2}{2}\left(\partial_t^2f_i+
		\mathbf{c}_i^T\mathbf{c}_i:\nabla^2f_i+2\mathbf{c}_i\cdot\nabla \partial_tf_i\right)=-\bm{\Lambda}_{ij}f_j^{\text{ne}}+\Delta t\left(\mathbf{I}-\frac{\bm{\Lambda}}{2}\right)_{ij}R_j+O(\Delta t^3),
	\end{align*}
	where the non-equilibrium distribution functions are defined as \(f_i^{\text{ne}} = f_i - f_i^{\text{eq}}\), which satisfy \(\sum_{i=0}^{2d} f_i^{\text{ne}} = 0\). From the above equation, we have
	\begin{subequations}
		\begin{align}
			&f_i^{\text{ne}}=O(\Delta x),\\
			&\left(\partial_tf_i^{\text{eq}}+\mathbf{c}_i\cdot \nabla f_i^{\text{eq}}\right)=-\frac{1}{\Delta t}\bm{\Lambda}_{ij}f_j^{\text{ne}}+\left(\mathbf{I}-\frac{\bm{\Lambda}}{2}\right)_{ij}R_j+O(\Delta t),\label{chapsecond-f-1st}\\
			&\left(\partial_tf_i+\mathbf{c}_i\cdot\nabla  f_i\right)+\frac{\Delta t}{2}\left(\partial_t^2f_i^{\text{eq}}+
			\mathbf{c}_i^T\mathbf{c}_i:\nabla^2f_i^{\text{eq}}+2\mathbf{c}_i\cdot\nabla \partial_tf_i^{\text{eq}}\right)=-\frac{1}{\Delta x}\bm{\Lambda}_{ij}f_j^{\text{ne}}+\left(\mathbf{I}-\frac{\bm{\Lambda}}{2}\right)_{ij}R_j+O(\Delta t^2).
		\end{align}
	\end{subequations}
	Taking the zeroth-order moment of the above equation yields
	\begin{subequations}
		\begin{align}
			&\partial_tu+\nabla\cdot \mathbf{F}(u)=O(\Delta x),\label{chapsecond-1me}\\
			&\partial_tu+\nabla\cdot \mathbf{F}(u)+\sum_{i=1}^{2d}\mathbf{c}_i\cdot \nabla f_i^{\text{ne}}+\frac{\Delta t}{2}\left(\partial_t^2u+
			\nabla^2cG_i(u)+2\nabla \cdot\partial_t\mathbf{F}\right)=O(\Delta t^2).\label{chapsecond-2me}
	\end{align}
	\end{subequations}
	Substituting \eqref{chapsecond-f-1st} and \eqref{chapsecond-1me} into \eqref{chapsecond-2me} gives
	\begin{align*}
		\partial_tu+\nabla\cdot \mathbf{F}(u)-\Delta x\sum_{i=0}^{2d}\mathbf{c}_i\cdot \nabla (\bm{\Lambda}^{-1})_{ij}\left(\partial_tf_j^{\text{eq}}+\mathbf{c}_j\cdot\nabla f_j^{\text{eq}}-\Big(\mathbf{I}-\frac{\bm{\Lambda}}{2}\Big)_{jk}R_k\right)+\frac{\Delta x}{2}\left(\partial_t\nabla\cdot\mathbf{F}(u)+
		\nabla^2cG_i(u)\right)=O(\Delta x^2),
	\end{align*}
	For the MRT model in \eqref{chapfirst-m} and \eqref{chapfirst-s}, we have \(\mathbf{c}_{i\alpha}(\bm{\Lambda}^{-1})_{ij}=\mathbf{c}_{j\alpha}/s_{x_{\alpha}}\) \cite{Chai2020}. Thus, above equation  becomes
	\begin{align*}
		\partial_tu+\nabla\cdot \mathbf{F}(u)+\Delta t\sum_{\alpha=1}^{d} \partial_{x_\alpha}\Big(\frac{1}{s_{x_\alpha}}-\frac{1}{2}\Big)\left[ \Big(cG_{\beta}'(u)\partial_{x_{\beta}}u-F_{\alpha}'(u)\mathcal A(u)\Big)-\Big(\partial_tF_{\alpha}(u)+
		c\partial_{x_{\beta}}G_{\beta}(u)\Big)\right]=O(\Delta t^2),
	\end{align*}
	which is second-order consistent with \eqref{chapsecond-eq} under the relation of \(\kappa_{i}=\Delta t(1/s_{x_i}-1/2)\). And particularly, from the expression of \(\mathcal A(u)\) in \eqref{chapsecond-au}, we obtain
	\begin{align*}
		\partial_tu+\sum_{\alpha=1}^d\partial_{x_\alpha}F_{\alpha}(u)=O(\Delta t^2).
	\end{align*}

	Several remarks on the second-order MRT-LB scheme \eqref{chapsecond-mrt-lbm} are given below.
	\begin{itemize}
		\item We would like to point out that, to guarantee stability of the above second-order MRT-LB method, the relaxation parameters \(s_{x_\alpha}\) and \(p_{x_\alpha}\) must lie within the interval \((0,2)\). Moreover, the method is equivalent to a \((2d+2)\)-level macroscopic finite-difference scheme for the conserved quantity \(u\) \cite{bellotti2022finite,chen2024gplb,bellotti2023truncation}. And particularly, when all relaxation parameters are set to 1, this finite-difference scheme reduces to a two-level one, given by
		\begin{equation*}
			\begin{aligned}
				u_K^{n+1}=u_K^n&-\frac{\Delta t}{\Delta x}\sum_{\alpha=1}^d\left[F_{\alpha}^+(u_K^n)-
				F_{\alpha}^+(u_{K-\mathbf{e}_\alpha}^n)+
				F_{\alpha}^-(u_K^n)-F_{\alpha}^-(u_{K+\mathbf{e}_{\alpha}}^n)\right]\\
				&+\frac{\Delta t  }{2        }\left[\partial_{\alpha}F_{\alpha}^+(u_K^n)-\partial_{\alpha}F_{\alpha}^+(u_{K-\mathbf{e}_\alpha}^n)
				+\partial_{\alpha}F_{\alpha}^-(u_K^n)-\partial_{\alpha}F_{\alpha}^-
				(u_{K+\mathbf{e}_\alpha}^n)\right]\\
				&-\frac{\Delta t^2}{2\Delta x}\Big[		(	F_{\alpha}^+)'(u_K^n)\mathcal A(u_K^n)
				-(F_{\alpha}^+)'(u_{K-\mathbf{e}_\alpha}^n)\mathcal A(u_{K-\mathbf{e}_\alpha}^n)
				+(F_{\alpha}^-)'(u_{K}^n)\mathcal A(u_{K}^n)
				-(F_{\alpha}^-)'(u_{K+\mathbf{e}_{\alpha}}^n)\mathcal A(u_{K+\mathbf{e}_{\alpha}}^n)\big)
				\Big].
			\end{aligned}
		\end{equation*}
		\item For numerical implementation, the spatial derivatives in the auxiliary source distribution functions \eqref{chapsecond-source} are discretized by a second-order central-difference based on the chain rule, and computed explicitly from the current-time macroscopic quantities, i.e.,
		\begin{align}\label{chapsecond-centeral-fd}
			\partial_{x_{\alpha}}F_{\beta}\big(u(\mathbf{x},t_n)\big)
			=F_{\beta}'\big(u(\mathbf{x},t_n)\big)
			\partial_{x_{\alpha}}\big(u(\mathbf{x},t_n)\big)\approx F_{\beta}		'\big(u(\mathbf{x},t_n)\big)\frac{u(\mathbf{x}+\mathbf{c}_{\alpha}\Delta t,t_n)-u(\mathbf{x}-\mathbf{c}_{\alpha}\Delta t,t_n)}{2\Delta x}.
		\end{align}
		In addition, to improve the regularity of the upwind-split at the point where the characteristic velocity changes sign, we replace \(|F_\alpha'(u)|\) with \(\sqrt{F_\alpha'(u)^2+\delta^2}\) (\(\delta\) is the smoothing parameter) in numerical simulations, then the equilibrium distribution functions are constructed consistently in the following forms:
		\begin{align}\label{chapfirst-upwind-feq-delta}
			&\begin{cases}
				f_0^{\text{eq},\delta}=u-\sum_{\alpha=1}^{d}\frac{F_{\alpha,\delta}^+(u)+F_{\alpha,\delta}^-(u)}{c},\\
				f_{\alpha}^{\text{eq},\delta}=\frac{F_{\alpha,\delta}^{+}(u)}{c},\qquad			f_{\alpha+d}^{\text{eq},\delta}=\frac{F_{\alpha,\delta}^{-}(u)}{c},\qquad {\alpha}=1,2,\ldots,d,
			\end{cases}
		\end{align}
		where 
		\begin{align*}
		F^{\pm}_{\alpha,\delta}=\frac{1}{2}\int_0^u\left[\sqrt{\big(F_{\alpha}'(v)\big)^2+\delta^2}\pm F_{\alpha}'(v)\right]dv,\qquad F^+_{\alpha,\delta}-F^-_{\alpha,\delta}=F_{\alpha}.
		\end{align*} 
		And the corresponding  auxiliary source distribution functions are given by
		\begin{align}\label{chapsecond-source-delta}
		\begin{cases}
			S_0^{\delta}=-\left[1-\frac{1}{c}\sum_{\alpha=1}^d  \big((F_{\alpha,\delta}^+)'+(F_{\alpha,\delta}^-)'\big) \right]\mathcal{A}(u),\\
			S_{\alpha}^{\delta}=\partial_{x_\alpha}F_{\alpha,\delta}^+-\frac{(F_{\alpha,\delta}^+)'}{c}\mathcal{A}(u),\qquad S_{\alpha+d}^{\delta}=-\partial_{x_\alpha}F_{\alpha,\delta}^+-\frac{(F_{\alpha,\delta}^-)'}{c}\mathcal{A}(u),\qquad{\alpha}=1,2,\ldots,d,
		\end{cases}
		\end{align}
	It is clearly that this treatment also preserves the moment conditions \eqref{chapfirst-moment} and \eqref{chapsecond-source-moment}.
	\item When the equilibrium distribution functions are the linear forms \eqref{chapfirst-linear-feq}, to ensure that the MRT-LB method \eqref{chapsecond-mrt-lbm} is second-order, the corresponding source distribution functions can be designed as \cite{Fu2024LevelSet}
	\begin{align}\label{chapsecond-source-linear}
			R_i=\omega_i\left[\frac{\mathbf{c}_i\cdot\big(c_s^2\nabla u+\partial_t\mathbf{F}\big)}{c_s^2}\right],\qquad i=0,1,\ldots,2d.
	\end{align}
	\end{itemize}
	 
	\subsection{$L^2$ stability of second-order lattice Boltzamnn method}
	In this part, we focus on the following linear scalar hyperbolic conservation laws:
	\begin{align*}
		\partial_tu+\mathbf{a}\cdot\nabla u= \nabla\cdot\left(\mathbf{K}\tilde{\mathbf{A}}\nabla u\right)+\left[-\nabla\cdot\big(\mathbf{K}\tilde{\mathbf{A}}\nabla u\big)\right],
	\end{align*} 
	with \(\mathbf{a}\in\mathbb{R}^{d\times1 }\)  and \(\tilde{\mathbf{A}}=\text{diag}\big(c\mathbf{a}\big)-\mathbf{a}\mathbf{a}^T\). The corresponding upwind-taype \(f_i^{\text{eq}}\) and \(R_i\) in \eqref{chapsecond-mrt-lbm} are given by
	\begin{align*}
		\begin{cases}
			f_0^{\text{eq}}=u-\sum_{\alpha=1}^{d}\frac{|a_{\alpha}|}{c},\\
			f_{\alpha}^{\text{eq}}=\frac{|a_{\alpha}|+a_{\alpha}}{2c},\qquad f_{\alpha+d}^{\text{eq}}=\frac{|a_{\alpha}|-a_{\alpha}}{2c}, 
			\qquad \alpha=1,2,\ldots,d,		\end{cases}
	\end{align*}
	and
	\begin{align*}
		\begin{cases}
			R_0=-\left(1-\frac{1}{c}\sum_{\alpha=1}^d |a_{\alpha}| \right)\mathbf{a}\cdot\nabla u,\\
			R_{\alpha}=\frac{|a_{\alpha}|+a_{\alpha}}{2c}\left(c\partial_{x_\alpha}u-\mathbf{a}\cdot\nabla u\right),\qquad
			R_{\alpha+d}=-\frac{|a_{\alpha}|-a_{\alpha}}{2c}\left(c\partial_{x_\alpha}u+\mathbf{a}\cdot\nabla u\right),\qquad \alpha=1,2,\ldots,d.
		\end{cases}			
	\end{align*}
	
	To analyze the \(L^2\) satbility of the second-order MRT-LB method, we now rewritten  \eqref{chapsecond-mrt-lbm} in the moment space:
	\begin{align}\label{chapsecond-mrt-lbm-moment}
		\mathbf{m}(\mathbf{x},t+\Delta t)=\mathbf{M}\mathbf{T}\mathbf{M}^{-1}(\mathbf{I}-\mathbf{S})\mathbf{m}+\mathbf{M}\mathbf{T}\mathbf{M}^{-1}\mathbf{S}\mathbf{m}^{\text{eq}}+\Delta t\mathbf{M}\mathbf{T}\mathbf{M}^{-1}\left(\mathbf{I}-\frac{\mathbf{S}}{2}\right)\mathbf{R}(\mathbf{x},t),
	\end{align}
	where the moment vectors are given by \(\mathbf{m}=\mathbf{Mf}\), \(\mathbf{m}^{\text{eq}}=\mathbf{Mf}^{\text{eq}}\), and  \(\mathbf{R}=\mathbf{M}(R_0,R_1,\ldots,R_{2d})^T\). The diagonal space-shift matrix \(\mathbf{T}\) is defined as
	\begin{align*}
		\mathbf{T} = \text{diag}\left(T_0^{\Delta x}, T_1^{\Delta x},\dots, T_{2d}^{\Delta x}\right) {\quad \text{with} \quad} T^{\Delta x}_{\mathbf{c}_{\alpha}}h(\mathbf{x},t)=h(\mathbf{x}-\mathbf{c}_{\alpha}\Delta t,t),
	\end{align*}
	here \(h\) denotes an arbitrary spatially dependent function.
	
	Takeing the discrete Fourier transfrom to \eqref{chapsecond-mrt-lbm-moment} yields the following amplification matrix \(\mathbf{G}\) of the second-order MRT-LB method:
	\begin{align}\label{chapsecond-g}
		\mathbf{G}(\bm{\xi})=\mathcal{T}(\bm{\xi})\mathcal{P}(\bm{\xi}),\qquad\bm{\xi}\in[-\pi,\pi]^2,
	\end{align}
	with
	\begin{subequations}
		\begin{align}
			&\mathcal{T}(\bm{\xi})=\mathbf{M}\hat{\mathbf{T}}\mathbf{M}^{-1},\label{chapsecond-T}\\
			&\mathcal{P}(\bm{\xi})=(\mathbf{I}-\mathbf{S})+\left[\mathbf{S}\bm{\varepsilon}^{\text{eq}}+\bigg(\mathbf{I}-\frac{\mathbf{S}}{2}\bigg)\bm{\varepsilon}^{\text{r}}\right]\mathbf{e}_0^T,\label{chapsecond-P}
		\end{align}
	\end{subequations}	
	where  
		\begin{align*}			
			&\hat{\mathbf{T}}=\text{diag}\left(1,e^{-i\xi_1},e^{-i\xi_2},\ldots,e^{-i\xi_d},e^{i\xi_1},\dots, e^{i\xi_d}\right),\\
			&\bm{\varepsilon}^{\text{eq}}=\left(1,\mathbf{a}^T,c|\mathbf{a}|^T\right)^T,\qquad \bm{\varepsilon}^{\text{r}}=i\left(0,|\mathbf{a}|.\mathbf{z}-\mathbf{a}\mathbf{a}^T\mathbf{z}/c,c\mathbf{a}.\mathbf{z}-|\mathbf{a}|\mathbf{a}^T\mathbf{z}\right)^T, 
		\end{align*}
	here the vector \(\mathbf{z}=(\sin\xi_{1},\sin\xi_{2},\ldots,\sin\xi_{d})^T\) originates from the central-difference approximation \eqref{chapsecond-centeral-fd} of the spatial derivative \(\partial_{x_\alpha}\), whose Fourier symbol in Fourier space reads \(\widehat{\partial_{x_\alpha}^h} = i(\sin\xi_\alpha) / \Delta x\).
	
	In particular, the explicit expressions for the matrices in \eqref{chapsecond-T} and \eqref{chapsecond-P} are given by
	\begin{align*}
		\mathcal{T}(\bm{\xi})=\begin{pmatrix}
			1&\frac{-i\mathbf{z}^T}{c}&\frac{(\mathbf{b}-\mathbf{1}_d)^T}{c^2}\\
			\mathbf{0}&\text{diag}(\mathbf{z})&-\frac{i}{c}\text{diag}(\mathbf b)\\
			\mathbf{0}&-ic\text{diag}(\mathbf b)&\text{diag}(\mathbf{z})
		\end{pmatrix},
	\end{align*}
	and
	\begin{align}\label{chapsecond-2-P} 
		\mathcal{P}(\bm{\xi})=\mathbf{I}+\begin{pmatrix}
			0&0&0\\
			i\mathbf{S}_1\mathbf{a}^T+i(\mathbf{I}_d-\mathbf{S}_1/2)(|\mathbf{a}|.\mathbf{z}-\mathbf{a}\mathbf{a}^T\mathbf{z}/c)&-\mathbf{S}_1&\mathbf{0}\\
			ic\mathbf{S}_1|\mathbf{a}|^T+i(\mathbf{I}_d-\mathbf{S}_1/2)(c\mathbf{a}.\mathbf{z}-|\mathbf{a}|\mathbf{a}^T\mathbf{z})&\mathbf{0}&-\mathbf{S}_2
		\end{pmatrix},
	\end{align}
	where \(\mathbf{b}=\left(\cos\xi_1,\cos\xi_2,\ldots,\cos\xi_d\right)^T\), \(\mathbf{S}_1=\text{diag}\left(s_{x_1},s_{x_2},\ldots,s_{x_d}\right)\), and \(\mathbf{S}_2=\text{diag}\left(p_{x_1},p_{x_2},\ldots,p_{x_d}\right)\).

\begin{figure}[h]
		\centering
		\subfloat[{$s_x=p_x$}]
		{\begin{minipage}{0.27\linewidth}
				\centering
				\includegraphics[width=\linewidth]{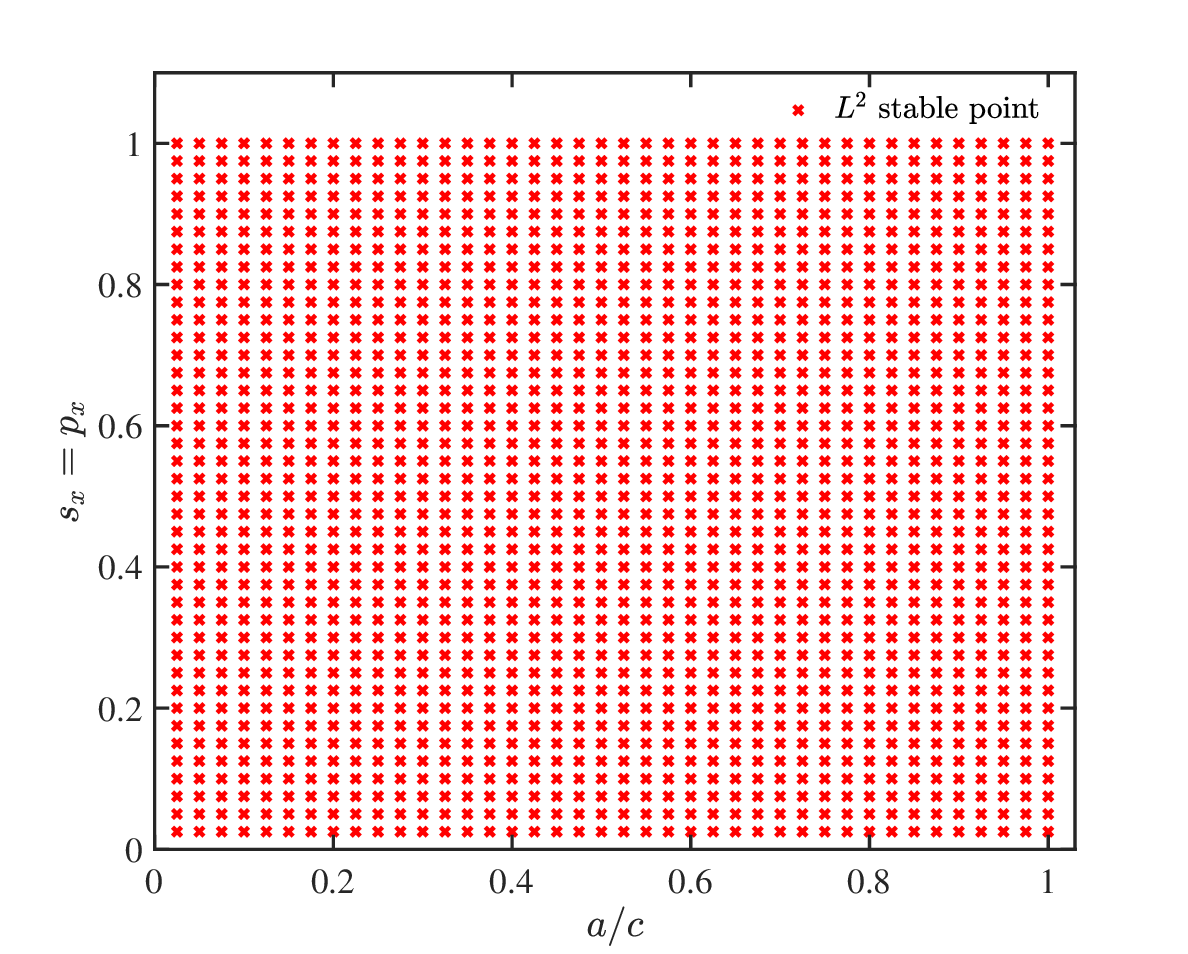}   
		\end{minipage}}   
		\subfloat[{$s_x=1$} ]
		{\begin{minipage}{0.27\linewidth}
				\centering
				\includegraphics[width=\linewidth]{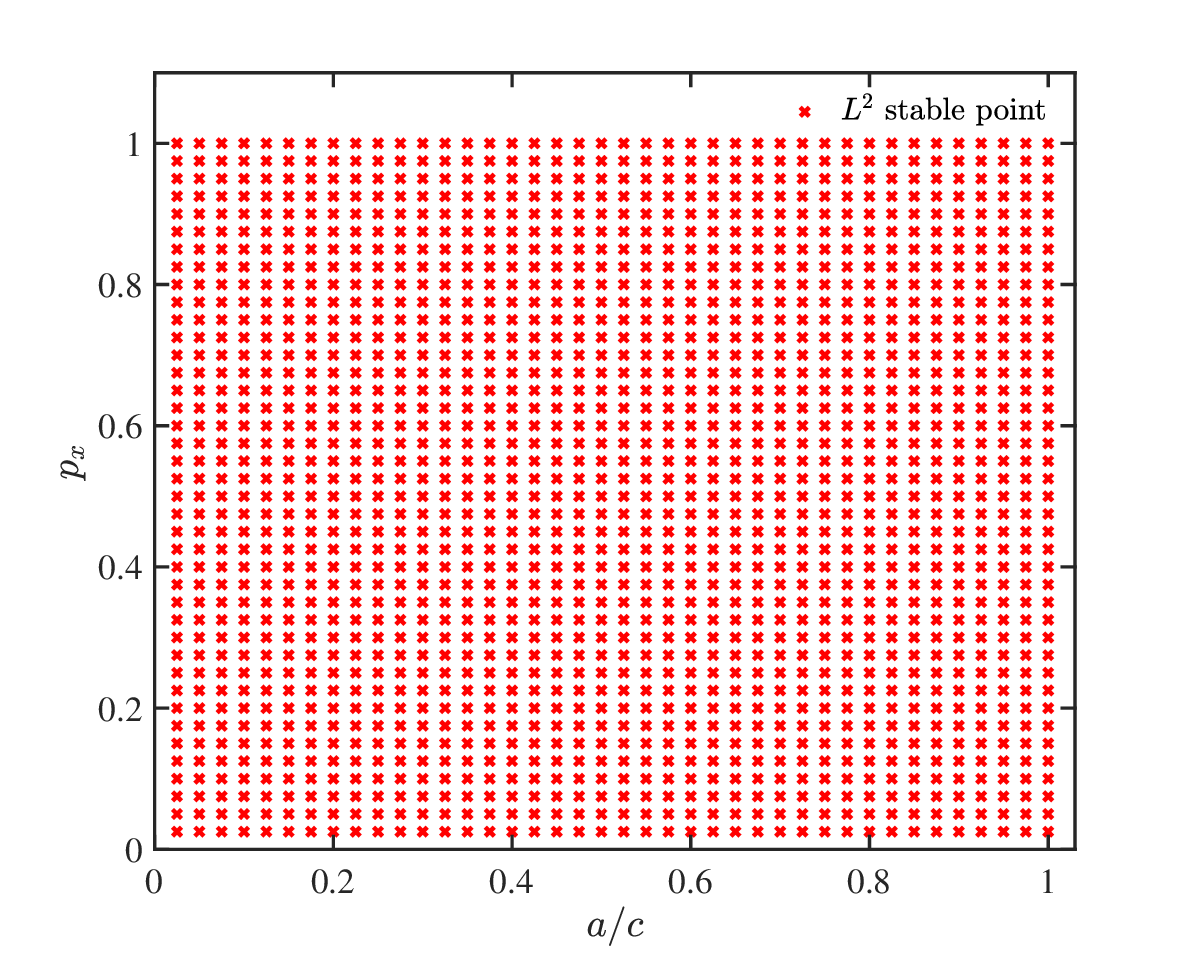} 	
		\end{minipage}}    
		\subfloat[$p_x=1$]
		{\begin{minipage}{0.27\linewidth}
				\centering
				\includegraphics[width=\linewidth]{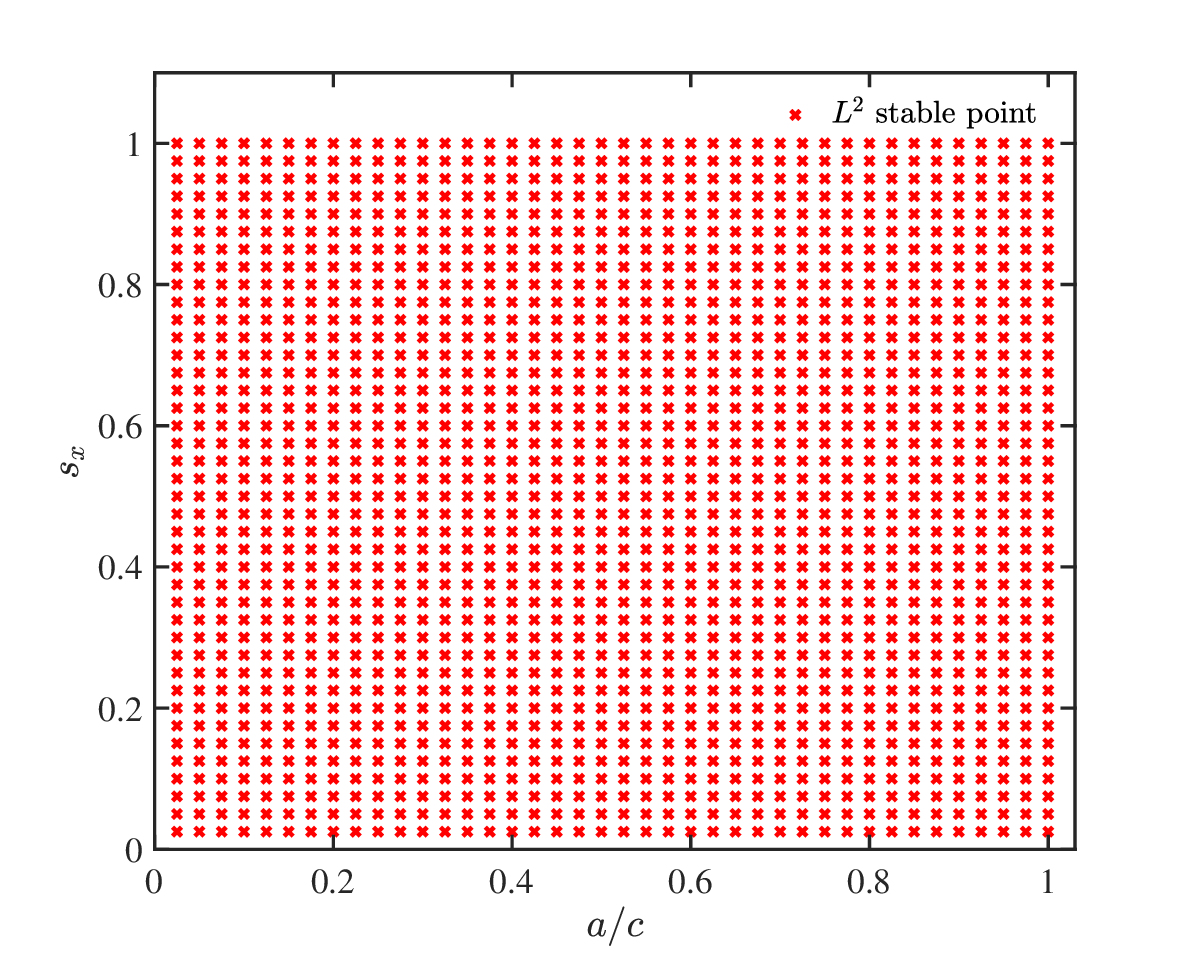} 	
		\end{minipage}}     		
		\caption{\(L^2\)-stable regions under different choices of the relaxation parameters in one-dimensional case.}
		\label{chapsecond-stable-1d}
	\end{figure}

\begin{figure}[h]
 		\centering
 		\subfloat[{$s_x=p_x=0.4,\quad s_y=p_y=0.8$}]
 		{\begin{minipage}{0.27\linewidth}
 				\centering
 				\includegraphics[width=\linewidth]{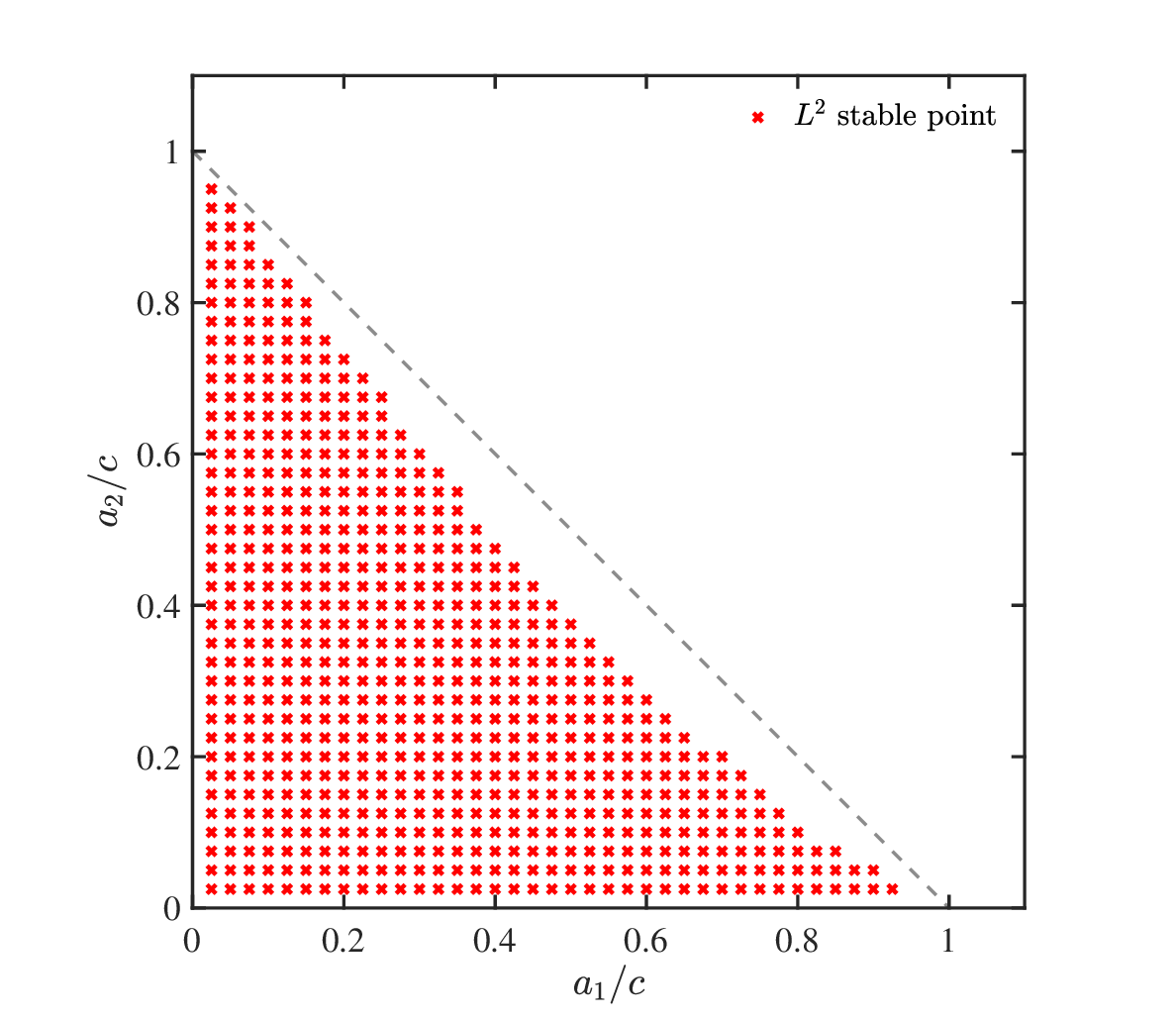}   
 		\end{minipage}}   
 		\subfloat[{$s_x=p_x=1.2,\quad s_y=p_y=1.6$} ]
 		{\begin{minipage}{0.27\linewidth}
 				\centering
 				\includegraphics[width=\linewidth]{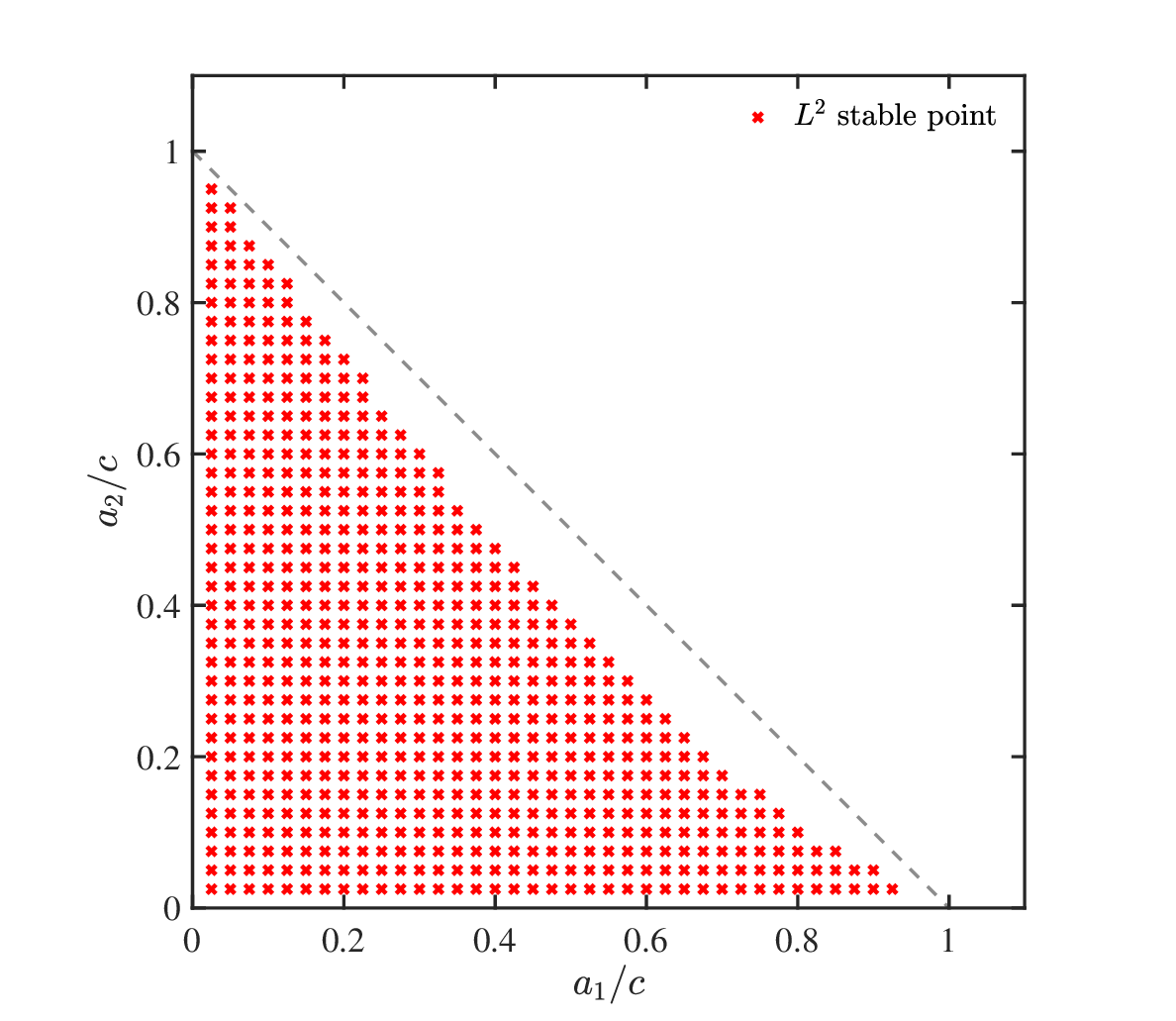} 	
 		\end{minipage}}    
 		\subfloat[$s_x=p_x=0.6,\quad s_y=p_y=1.4$]
 		{\begin{minipage}{0.27\linewidth}
 				\centering
 				\includegraphics[width=\linewidth]{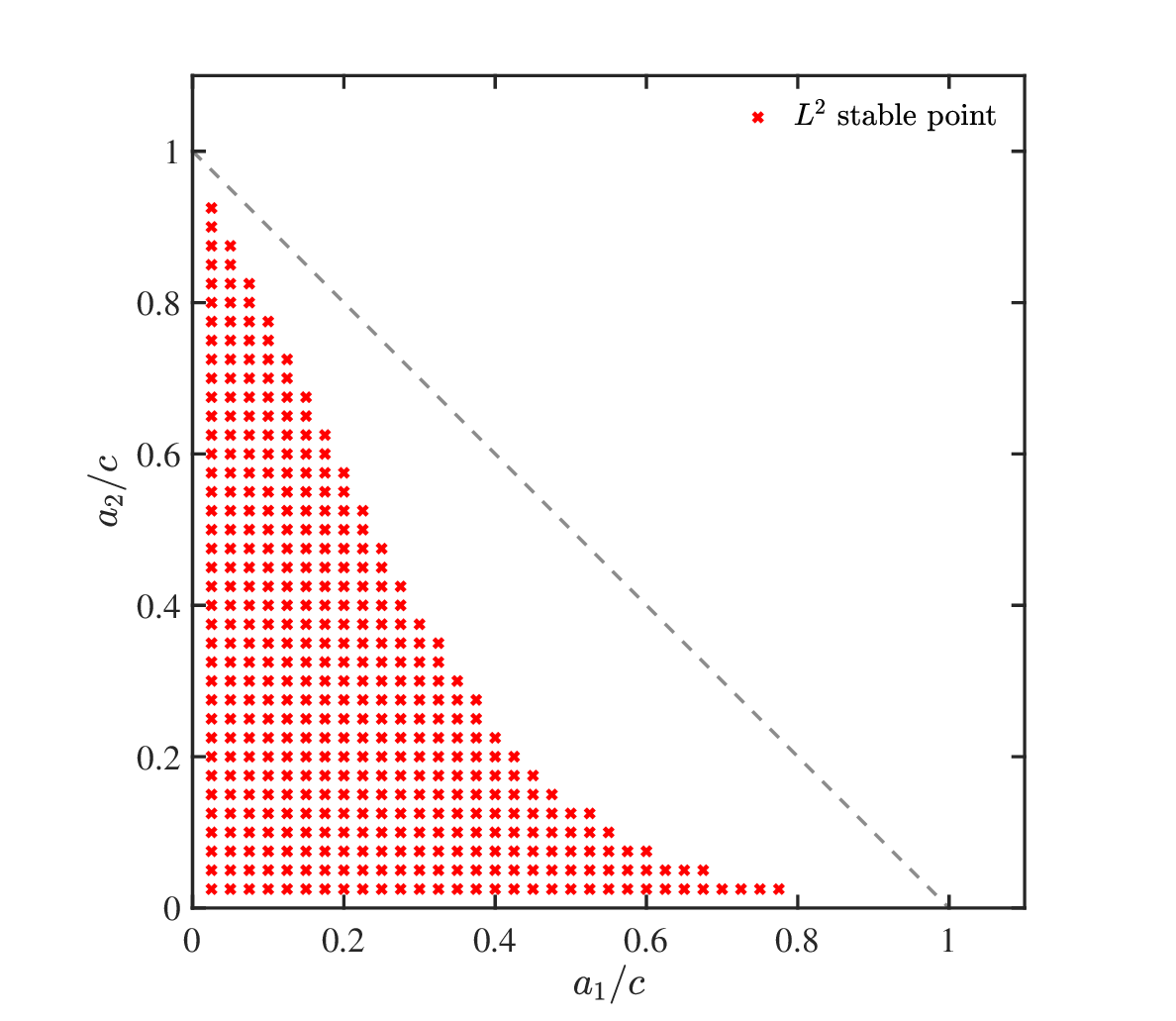} 	
 		\end{minipage}}     		
 		\caption{\(L^2\)-stable regions under different choices of the relaxation parameters in two-dimensional case.}
 		\label{chapsecond-stable-2d}
 	\end{figure}

	For the MRT model, a rigorous analytical derivation of \(L^2\) stability conditions is intractable. We therefore numerically investigate the \(L^2\) stability regions of the second-order MRT-LB method in one- and two-dimensional case, as shown in Figures  \ref{chapsecond-stable-1d} and \ref{chapsecond-stable-2d}. The results are obtained under different choices of the relaxation parameters, and without loss of generality, we consider only the case \(a_{\alpha} > 0\) here. It is observed that the stability region is bounded by \(\sum_{\alpha=1}^d |a_\alpha| \leq c\). In the following, we focus on the second-order LB method with all relaxation parameters set to 1, and derive the corresponding necessary and sufficient \(L^2\) stability conditions.

	\noindent
	\textbf{Theorem 2 (Necessary and sufficient condition for \(L^2\) stability).}
	The second-order LB method \eqref{chapsecond-mrt-lbm} with \(\mathbf{S}=\mathbf{I}\) is \(L^2\)-stable if and only if \(\sum_{\alpha=1}^{d} \lvert a_\alpha \rvert \le c\).\\	
	\textbf{Proof.} For the case \(s_{x_\alpha} = p_{x_\alpha} = 1\), \(\alpha = 1, 2, \dots, d\), the matrix \(\mathcal{P}(\bm{\xi})\) in \eqref{chapsecond-2-P} reduces to a rank-one matrix, i.e.,
	\[
	\mathcal{P}(\bm{\xi}) = \left( \bm{\varepsilon}^{\text{eq}} + \frac{\bm{\varepsilon}^{\text{r}}}{2} \right) \mathbf{e}_0^T,
	\]
	so that the amplification matrix \(\mathbf{G}\) in \eqref{chapsecond-g} has exactly one nonzero eigenvalue, given by
	\begin{equation}\label{chapsecond-root}
		\lambda(\bm{\xi}) = \mathbf{e}_0^T \mathcal{T}(\bm{\xi}) \left( \bm{\varepsilon}^{\text{eq}} + \frac{\bm{\varepsilon}^{\text{r}}}{2} \right) = C + \frac{Y}{2} - i \left( V - \frac{Q}{2} \right),
	\end{equation}
	where
	\begin{align*}
		C &= 1 + \frac{1}{c} (\mathbf{b} - \mathbf{1})^T \lvert \mathbf{a} \rvert, &
		V &= \frac{\mathbf{a}^T \boldsymbol{z}}{c}, &
		Y &= \frac{1}{c} \lvert \mathbf{a} \rvert^T (\mathbf{z}. \mathbf{z}) - V^2, &
		Q &= \frac{1}{c} \mathbf{a}^T (\mathbf{b}. \mathbf{z}) - C V.
	\end{align*}
	Under the assumption \(\sum_{\alpha=1}^{d} \lvert a_\alpha \rvert \le c\), we define the weight coefficients
	\begin{align*}
		\gamma_0 = 1 - \frac{1}{c} \sum_{\alpha=1}^d \lvert a_\alpha \rvert \ge 0, \qquad
		\gamma_\alpha = \frac{\lvert a_\alpha \rvert + a_\alpha}{2c} \ge 0, \qquad
		\gamma_{\alpha+d}= \frac{\lvert a_{\alpha+d} \rvert + a_{\alpha+d}}{2c} \ge 0,  
	\end{align*}
	where \(\alpha=1,2,\ldots,d\). Let \(\bm{\phi} = (0, \mathbf{z}, -\mathbf{z})\), we then obtain
	\begin{align*}
		C = \mathbb{E}(\cos \bm{\phi}), \qquad
		V = \mathbb{E}(\sin \bm{\phi}), \qquad
		Y = \operatorname{Var}(\sin \bm{\phi}), \qquad
		Q = \operatorname{Cov}(\cos \bm{\phi}, \sin \bm{\phi}).
	\end{align*}
	Due to the fact of \(\cos^2 \phi_i + \sin^2 \phi_i = 1\) for \(i = 0, 1, \dots, 2d\), we have \(X + Y = 1 - C^2 - V^2\) with \(X = \operatorname{Var}(\cos \bm{\phi})\). Then \eqref{chapsecond-root} gives
	\begin{align*}
		1 - \lvert \lambda(\bm{\xi}) \rvert^2
		&= 1 - \left( C + \frac{Y}{2} \right)^2 - \left( V - \frac{Q}{2} \right)^2
		= X + (1-C)Y + VQ - \frac{Y^2 + Q^2}{4} \\
		&= \operatorname{Var}\left( \cos \bm{\phi} + \frac{V}{2} \sin \bm{\phi} \right)
		- \frac{Y}{4}(1-C)(3-C) + \frac{1}{4}(XY - Q^2).
	\end{align*}
	Since the variance is nonnegative, \(\lvert C \rvert \le 1\), \(Y \ge 0\), and \(Q^2 \le XY\) (the Cauchy-Schwarz inequality), we can obtain \(1 - \lvert \lambda \rvert^2 \ge 0\).
	
	To establish necessity, we take \(\bm{\xi} =\pi\mathbf{1}^T\). Then \(\mathbf{z} = \mathbf{0}\) and \(\bm{\varepsilon}^{\text{r}} = \mathbf{0}\), and hence
	\[
	\lambda(\bm{\xi}) = 1 - \frac{2}{c} \sum_{\alpha=1}^d \lvert a_\alpha \rvert.
	\]
	If \(\sum_{\alpha=1}^d \lvert a_\alpha \rvert > c\), then \(\lambda(\bm{\xi}) < -1\), i.e., \(\lvert \lambda(\bm{\xi}) \rvert > 1\), so the method is unstable. In summary, for the second-order LB method \eqref{chapsecond-mrt-lbm} with all relaxation parameters set to 1, the necessary and sufficient condition for \(L^2\) stability is \(\sum_{\alpha=1}^{d} \lvert a_\alpha \rvert \le c\). This completes the proof.\hfill$\square$

	 It should be noted here that for the second-order LB method \eqref{chapsecond-mrt-lbm}, in which all relaxation parameters are set to 2 (meaning there is no auxiliary source distribution functions), it can be proven that the roots of its amplification matrix on the unit circle are multiple roots. Therefore, the \(L^2\) stability is not strictly guaranteed \cite{Bellotti2024Parasitic,chen2026fourth-order}. Consequently, while this second-order LB method may be sufficient for linear problems, it exhibits severe instability in nonlinear cases \cite{BellottiHelluyNavoret2025}.

 	\section{Bound-preserving hybrid lattice Boltzmann method}\label{chap3}
 	In this section, based on the first-order and second-order LB methods developed above, we first construct a conservative hybrid LB method using a convex limiter. For the interface coefficients in the convex limiter, we further adopt the GBP, LBP, and RLBP limiters \cite{WissocqLiuAbgrall2025}, together with the proofs of macroscopic bound preservation. We finally present the algorithmic flowchart of the proposed hybrid LB method for \eqref{chapintro-eq}.
 	
 	\subsection{Construction of the hybrid LB method}\label{chap3-part-1}
 	
	 	In the  framework of the LB method, the first-order and second-order schemes are adopted from \eqref{chapfirst-mrt-lbm} (with \(\mathbf{S} = \mathbf{I}\)) and \eqref{chapsecond-mrt-lbm}, respectively. For brevity, we define their updated distribution vectors at node \(K\) as
 	\begin{subequations}
 		\begin{align}
 			&\mathbf{L}_K = \mathbf{f}^{\text{eq}}, \label{chap3-first-scheme} \\
 			&\mathbf{H}_K = \mathbf{f}_K - \bm{\Lambda}\left(\mathbf{f} - \mathbf{f}^{\text{eq}}\right)_K + \Delta t \left( \mathbf{I} - \frac{\bm{\Lambda}}{2} \right) \mathbf{R}_K, \label{chap3-second-scheme}
 		\end{align}
 	\end{subequations}

 	Now, let \(\mathbf{D}_K = \mathbf{H}_K - \mathbf{L}_K\). For the interfaces \(e^\pm_{\alpha} = K \pm \mathbf{e}_\alpha/2\) (\(\alpha = 1,2, \dots, d\)), we define the following dimensionless interface corrections:
 	\begin{align*}
 		A_{e_\alpha^+} = D_{\alpha,K} - D_{\alpha+d, K+\mathbf{e}_{\alpha}},\qquad A_{e_\alpha^-} = D_{\alpha+d, K-\mathbf{e}_{\alpha}}-D_{\alpha,K},
 	\end{align*}
 	and design the hybrid LB method as
 	\begin{subequations}\label{chap3-collision}
 			\begin{align}
 			&f^{\star}_{0,K}=u_K-\sum_{\alpha=1}^d\left({f}_{\alpha,K}+{f}_{\alpha+d,K}\right),\label{chap3-collision-0}\\
 			&f^{\star}_{\alpha,K}=L_{\alpha,K}+\eta_{\alpha}^+D_{\alpha,K},\qquad f^{\star}_{\alpha+d,K}=L_{\alpha+d,K}+\eta_{\alpha}^-D_{\alpha+d,K},\label{chap3-collision-fi}
 		\end{align}
	\end{subequations}  
 where the interface coefficient \(\eta_{\alpha}^{+}\) (or \(\eta_{\alpha}^{-}\)) is assigned to the interface shared by the two adjacent nodes \(K\) and \(K + \mathbf{e}_{\alpha}\) (or \(K - \mathbf{e}_{\alpha}\)). The streaming step then reads
	\begin{align}\label{chap3-stream}
 		f_{0,K}^{n+1}=f_{0,K}^{\star},\qquad f_{\alpha,K}^{n+1}=f_{\alpha,K-\mathbf{e}_{\alpha}}^{\star},\qquad f_{\alpha+d,K}^{n+1}=f_{\alpha+d,K+\mathbf{e}_{\alpha}}^{\star}.
 	\end{align}
 	Setting all \(\eta_{e_\alpha^{\pm}} = 0\) recovers the first-order scheme \eqref{chap3-first-scheme}, whereas setting all \(\eta_{e_\alpha^{\pm}} = 1\) reduces the scheme to the second-order update \(f_{i,K} = H_{i,K}\) by \eqref{chap3-second-scheme}.
 	
 	We would like to point out that, the hybrid LB method in \eqref{chap3-collision}-\eqref{chap3-stream}  can also be expressed in the conservation form of a numerical flux for the macroscopic quantity \(u\). Combining the streaming step \eqref{chap3-stream} with the expression for \(f^{\star}_0\) in \eqref{chap3-collision-0} gives
 	\begin{align}\label{chap3-flux-eq-0}
 	u_K^{n+1}=u_K+\sum_{\alpha=1}^d
 	(f^{\star}_{\alpha,K-\mathbf e_\alpha}-f^{\star}_{\alpha,K}
 	+f^{\star}_{\alpha+d,K+\mathbf e_\alpha}-f^{\star}_{\alpha+d,K}),
 	\end{align}
 	and at each interface \(e^\pm_{\alpha}\), we can define the numerical flux shared by adjacent nodes \(K\) and \(K\pm\mathbf{e}_{\alpha}\) as
 	\begin{align*}
 		\begin{cases}
 			\mathcal F_{e_{\alpha}^+}=c(f^{\star}_{\alpha,K}-f^{\star}_{\alpha+d,K+\mathbf{e}_{\alpha}})
 			=c(L_{\alpha,K}-L_{\alpha+d,K+\mathbf{e}_{\alpha}})+c\eta_\alpha^+(D_{\alpha,K}-D_{\alpha+d,K+\mathbf{e}_{\alpha}})
 			,\\
 			\mathcal F_{e_{\alpha}^-}=c(f^{\star}_{\alpha+d,K-\mathbf{e}_{\alpha}}-f^{\star}_{\alpha,K})
 			=c(L_{\alpha+d,K-\mathbf{e}_{\alpha}}-L_{\alpha,K})+c\eta_\alpha^-(D_{\alpha+d,K-\mathbf{e}_{\alpha}}-D_{\alpha,K})
 			.
 		\end{cases}\end{align*}
 	Substituting the above equations into \eqref{chap3-flux-eq-0} yields
 	\begin{equation*}
 		u_K^{n+1}=u_K-\frac{\Delta t}{\Delta x}\sum_{\alpha=1}^d
 		(\mathcal F_{e^+_\alpha}-\mathcal F_{e^-_\alpha}).
 	\end{equation*}
 	Equivalently,
 	\begin{equation}\label{eq:uk-update}
 		u_K^{n+1}=u_K^L+\sum_{\alpha=1}^d
 		\left(\eta_{\alpha}^-A_{e_{\alpha}^-}
 		-\eta_{\alpha}^+A_{e_{\alpha}^+}\right).
 	\end{equation}
  
 	\subsection{Determination of interface coefficients}
 	The conservation property of the hybrid LB method presented in Section \ref{chap3-part-1} holds for any shared interface coefficients. In this part, by determining the permissible intervals and interface coefficients, we ensure that this hybrid LB method simultaneously satisfies the corresponding macroscopic BP constraints.
 	
 	First, we define the low-order intermediate states on both sides of the interfaces \(e_{\alpha}^\pm=K\pm\mathbf{e}_{\alpha}/2\) (\(\alpha=1,2,\ldots,d\)) as
 	\begin{align*}
 		B_{\alpha,K}^+=u_K+2d(f^{\text{eq}}_{\alpha+d}(u_{K+\mathbf{e}_{\alpha}})-f^{\text{eq}}_{\alpha+d}(u_K)),\qquad B_{\alpha,K}^-=u_K+2d(f^{\text{eq}}_{\alpha}(u_{K-\mathbf{e}_{\alpha}})-f^{\text{eq}}_{\alpha}(u_K)).
 	\end{align*}
  	Then the macroscopic quantity \(u_K\) from the first-order LB method \eqref{chap3-first-scheme} is given by
 	\begin{align}\label{eq:1st-update}
 		u_K^L=f^{\text{eq}}_0(u_K)+\sum_{\alpha=1}^d\left(f^{\text{eq}}_\alpha(u_K-\mathbf{e}_\alpha)+f^{\text{eq}}_{\alpha+d}(u_K+\mathbf{e}_\alpha)\right)
 		=\frac{1}{2d}\sum_{\alpha=1}^d\left(B_{\alpha,K}^++B_{\alpha,K}^-\right).
 	\end{align}
 	Now we present a lemma characterizing the BP property of the interface states \(B_{\alpha,K}^{\pm}\).\\
 	\textbf{Lemma 1} \textbf{(Two-point bound preservation for interface states).} If, for all $u\in I$, the equilibrium distribution functions satisfy
 	\begin{equation}\label{eq:facecondition}
 		0 \le 2d \frac{\mathrm{d} f_i^{\text{eq}}(u)}{\mathrm{d} u} \le 1, \qquad i = 1, 2, \dots, 2d,
 	\end{equation}
 	then each interface state \(B_{K,\alpha}^{\pm}\) is bounded between the macroscopic values of its two adjacent nodes, i.e., 
 	\begin{align}\label{eq:minmax-immerdate-state}
 		\min\{u_K,u_{K\pm\mathbf{e}_{\alpha}}\}\leq B_{\alpha,K}^\pm\leq \max\{u_K,u_{K\pm\mathbf{e}_{\alpha}}\}.
 	\end{align}
 	\noindent\textbf{Proof.}
 	Without loss of generality, we take \(B_{K,\alpha}^+\) as a representative case (the analysis for \(B_{K,\alpha}^-\) is similar). If $u_{K+\mathbf{e}_{\alpha}} = u_K$, then $B_{K,\alpha}^+ = u_K$.
 	Otherwise, let
 	\[
 	\beta = 2d \int_0^1 \frac{\mathrm{d} f_{\alpha+d}^{\text{eq}}}{\mathrm{d} u} \Bigl( u_K + t(u_{K+\mathbf{e}_\alpha} - u_K) \Bigr) \mathrm{d} t,\qquad \alpha=1,2,\ldots,d.
 	\]
 	By \eqref{eq:facecondition}, we have $0\le\beta\le1$. Applying the fundamental theorem of calculus then yields 
 	\[B_{K,\alpha}^+ = (1-\beta) u_K + \beta u_{K+\mathbf{e}_{\alpha}},\]
 	which establishes the desired bound. This completes the proof.\hfill$\square$
 	
 	For the upwind-type equilibrium distribution \(f_i^{\text{eq}}\) in \eqref{chapfirst-upwind-feq}, condition \eqref{eq:facecondition} reduces to
 	\begin{equation*}
 		c \ge 2d \max_{\alpha} \sup_{u\in I} \lvert F_\alpha'(u) \rvert.
 	\end{equation*}
 	This condition is generally stronger than the monotonic and nondecreasing conditions imposed on \(f_i^{\text{eq}}\) in \eqref{chapfirst-feq-upwind-monot}. In this work, we therefore employ \eqref{eq:facecondition} merely as a sufficient and readily verifiable condition.
 	
 	To ensure that the hybrid LB method \eqref{chap3-collision}-\eqref{chap3-stream} is BP, we now determine the hybrid interface coefficients \(\eta_{\alpha}^\pm\) (\(\alpha=1,2,\ldots,d\)). First, we require the interface states \(B_{\alpha,K}^{\pm}\) to lie in an given interval \([\ell_K, r_K]\), i.e.,
 	\begin{equation}\label{eq:basefeasible}
 		\ell_K \le B_{\alpha,K}^{\pm} \le r_K, \qquad \alpha = 1,2, \dots, d.
 	\end{equation}
 	
 	In the following, we specify the endpoints \(\ell_K\) and \(r_K\) for the GBP, LBP, and RLBP constraints, respectively. These constraints are then used to determine \(\eta_{\alpha}^{\pm}\) (\(\alpha = 1, 2,\dots, d\)), yielding the corresponding GBP, LBP, and RLBP limiters.

 	\begin{itemize}
 		\item For the GBP limiter, we take \(\ell_K^{\text{G}} = m_0\) and \(r_K^{\text{G}} = M_0\). If \(u_K \in [m_0, M_0]\) and the conditions of Lemma 1 hold on \(I = [m_0, M_0]\), then \eqref{eq:basefeasible} is satisfied. What is more, for interfaces $e_{\alpha}^\pm=K\pm\mathbf{e}_\alpha/2$ and $J=K\pm\mathbf{e}_\alpha$, we also require
 			\begin{align}\label{eq:faceconstraints}
 				\ell_K^{\text{G}}\le B_{K,\alpha}^\pm\mp 2d\eta_{\alpha}^\pm A_{e_{\alpha}^\pm}\le r_K^{\text{G}},\qquad
 				\ell_J^{\text{G}} \le B_{J,\alpha}^\mp \pm 2d\eta_\alpha^\pm A_{e_\alpha^\pm} \le r_J^{\text{G}}.
 			\end{align}
 		Selecting the maximum coefficient in $[0,1]$ that satisfies the constraints on both sides gives
 		\begin{equation}\label{eq:eta}
 			\eta_\alpha^+=\begin{cases}
 				1,&A_{e_\alpha^+}=0,\\
 				\displaystyle\min\left\{1,\frac{B_{K,\alpha}^+-\ell_K^{\text{G}}}{2dA_{e_\alpha^+}},
 				\frac{r_J^{\text{G}}- B_{J,\alpha}^-}{2dA_{e_\alpha^+}}\right\},&A_{e_\alpha^+}<0,\\ 
 				\displaystyle\min\left\{1,\frac{r_K^{\text{G}}-B_{K,\alpha}^+}{-2dA_{e_\alpha^+}},
 				\frac{B_{J,\alpha}^--\ell_J^{\text{G}}}{-2dA_{e_\alpha^+}}\right\},&A_{e_\alpha^+}>0,
 			\end{cases}\qquad
 			\eta_\alpha^-=\begin{cases}
 				1,&A_{e_\alpha^-}=0,\\
 				\displaystyle\min\left\{1,\frac{B_{J,\alpha}^--r_J^{\text{G}}}{2dA_{e_\alpha^-}},
 				\frac{\ell_K^{\text{G}}- B_{K,\alpha}^-}{2dA_{e_\alpha^-}}\right\},&A_{e_\alpha^-}<0,\\ 
 				\displaystyle\min\left\{1,\frac{\ell_J^{\text{G}}-B_{J,\alpha}^+}{-2dA_{e_\alpha^-}},
 				\frac{B_{K,\alpha}^--r_K^{\text{G}}}{-2dA_{e_\alpha^-}}\right\},&A_{e_\alpha^-}>0.
 			\end{cases}
 		\end{equation}
 		
 		We now establish the bound prerservation of the proposed hybrid LB method, with the interface coefficients \(\eta_{\alpha}^{\pm}\) (\(\alpha = 1, 2,\dots, d\)) chosen according to \eqref{eq:eta}.\\ 	
 		\textbf{Theorem 3 (Bound preservation of the hybrid LB method).} Assume that \eqref{eq:facecondition} holds on \([m_0, M_0]\). If \(u_K^0 \in [m_0, M_0]\), then, using the endpoints \(\ell_K^{\text{G}} = m_0\) and \(r_K^{\text{G}} = M_0\) in the GBP  limiter  and the interface coefficients given by \eqref{eq:eta}, the numerical solution generated by  \eqref{chap3-collision}-\eqref{chap3-stream} satisfies \(u_K^n \in [m_0, M_0]\) for all \(n \ge 0\).\\	
 		\textbf{Proof.} From \eqref{eq:uk-update} and \eqref{eq:1st-update}, we obtain
 		\begin{equation}\label{eq:limitedaverage}
 			u_K^{n+1} = \frac{1}{2d} \sum_{\alpha=1}^{2d} \Big[ \big( B_{K,\alpha}^+ - 2d \eta_{\alpha}^+ A_{e_\alpha^+} \big) + \big( B_{K,\alpha}^- + 2d \eta_{\alpha}^- A_{e_\alpha^-} \big) \Big].
 		\end{equation}
 		Assume that the macroscopic quantities \(u_K^n\) are bounded. Lemma 1 guarantees that all low-order interface states \(B_{K,\alpha}^{\pm}\) lie within the global interval \([m_0, M_0]\). Combining this with the constraint \eqref{eq:faceconstraints}, we see that \eqref{eq:limitedaverage} is the equally weighted convex combination of values lying in \([m_0, M_0]\), and hence remains in \([m_0, M_0]\). Together with the induction hypothesis, this completes the proof.\hfill$\square$
 		\item 	For the LBP limiter, the local endpoints at node \(K\) are defined as
 		\begin{equation}\label{eq:lmp}
 			\ell_K^{\text{L}} = \min\{ u_K, B_{1,K}^+, B_{1,K}^-, \dots, B_{d,K}^+, B_{d,K}^- \}, \qquad
 			r_K^{\text{L}} = \max\{ u_K, B_{1,K}^+, B_{1,K}^-, \dots, B_{d,K}^+, B_{d,K}^- \}.
 		\end{equation}
 		
 		\textbf{Theorem 4 (Local bound preservation of the hybrid LB method).} Using \eqref{eq:eta} and \eqref{eq:lmp}, we have
 		\[
 		\ell_K^{\text{L}} \le u_K^{n+1} \le r_K^{\text{L}},
 		\]
 		at each time step. In addition,  If the conditions of Lemma 1 hold, then
 		\begin{equation}\label{eq:neighbor}
 			\min_{J\in\mathcal N_K} u_J^n \le u_K^{n+1} \le \max_{J\in\mathcal N_K} u_J^n, \text{\quad with \quad}
 			\mathcal N_K = \{ K, K\pm\mathbf{e}_\alpha \}, \qquad \alpha = 1,2, \dots, d.
 		\end{equation}
 		
 		\textbf{Proof.} By construction, expression \eqref{eq:lmp} includes all low-order interface states \(B_{\alpha,K}^{\pm}\) (\(\alpha = 1,2, \dots, d\)) of node \(K\), and therefore satisfies \eqref{eq:basefeasible}. Lemma 1 together with \eqref{eq:limitedaverage} yields \(\ell_K^{\text{L}} \le u_K^{n+1} \le r_K^{\text{L}}\). Moreover, since \eqref{eq:minmax-immerdate-state} holds and \(u_K\) itself lies in the interval \([\min_{J\in\mathcal N_K} u_J, \max_{J\in\mathcal N_K} u_J]\), we obtain the inclusion
 		\[
 		[\ell_K^{\text{L}}, r_K^{\text{L}}] \subseteq [\min_{J\in\mathcal N_K} u_J, \max_{J\in\mathcal N_K} u_J],
 		\]
 		which yields \eqref{eq:neighbor}. Applying \eqref{eq:neighbor} successively over the periodic grid, we recover the bound preservation of the initial global bound.\hfill$\square$

 		\item For the RLBP limiter, we take a relaxed parameter \(C_{\text{relax}} \ge 0\) and define
 		\begin{equation}\label{eq:rlmp}
 			\ell_K^{\text{R}} = \ell_K^{\text{L}} - C_{\text{relax}} \lvert \ell_K^{\text{L}} \rvert, \qquad
 			r_K^{\text{R}} = r_K^{\text{L}} + C_{\text{relax}} \lvert r_K^{\text{L}} \rvert.
 		\end{equation}
 		
 		\textbf{Theorem 5 (Relaxed local bound preservation of the hybrid LB method).} Taking \eqref{eq:rlmp} as the endpoints in \eqref{eq:eta}, we have, at each time step,
 		\begin{equation}\label{bound-rlmp}
 			\ell_K^{\text{R}} \le u_K^{n+1} \le r_K^{\text{R}}.
 		\end{equation}
 		
 		\textbf{Proof.} From Theorem 3, we have \(\ell_K^{\text{R}} \le \ell_K^{\text{L}} \le B_{\alpha,K}^{\pm} \le r_K^{\text{L}} \le r_K^{\text{R}}\). Hence the low-order interface states \(B_{\alpha,K}^{\pm}\) (\(\alpha = 1, 2,\dots, d\)) lie within the relaxed interval \([\ell_K^{\text{R}}, r_K^{\text{R}}]\). Lemma 1 together with \eqref{eq:limitedaverage} then yields the conclusion directly. We note that the bound \eqref{bound-rlmp} does not guarantee a strict global bound in general.\hfill$\square$

 	\end{itemize}

 	\subsection{Algorithm implementation procedure}
 	We summarize the main implementation steps of the proposed BP hybrid LB method as follows.
 	
 	Given $u_K^0\in[m_0,M_0]$, we initialize $f_{i,K}^0=f_i^{\text{eq}}(u_K^0)$. At each time step:
 	\begin{itemize}
 		\item Compute the macroscopic quantities \(u_K^n = \sum_{i=0}^{2d} f_{i,K}^n\) at node \(K\) and time \(t_n\).
 		\item Compute the first-order updates \(L_{i,K} = f_i^{\text{eq}}(u_K^n)\) from \eqref{chap3-first-scheme}, and the second-order updates \(H_{i,K}\) from \eqref{chap3-second-scheme}.
 		\item Compute \(\mathbf{D}_K = \mathbf{H}_K - \mathbf{L}_K\) for each node \(K\), the corrections \(A_{e_\alpha^{\pm}}\) for each interface, and the interface states \(B_{K,\alpha}^{\pm}\).
 		\item Select the upper and lower bounds using the GBP , LBP, or RLBP limiter.
 		\item Compute \eqref{eq:eta} once per geometric interface and share the result between the two adjacent nodes. Set \(\eta_\alpha^{\pm} = 1\) when \(A_{e_\alpha^{\pm}} = 0\) to avoid division by zero.
 		\item Apply the collision step via \eqref{chap3-collision}, then reconstruct \(f^{\star}_0\) from \eqref{chap3-collision-0}.
 		\item Apply the streaming step via \eqref{chap3-stream} and obtain \(u_K^{n+1} = \sum_{i=0}^{2d} f_{i,K}^{n+1}\).
 	\end{itemize}

	\section{Numerical results and discussion}\label{chap4}
	In this section, we present a series of numerical tests to verify the convergence order, interface capturing performance, and discontinuity resolution capability of the proposed BP hybrid LB method. To evaluate numerical performance, the \(L^1\), \(L^2\), and \(L^\infty\) error norms are defined as
	\begin{align*}
		\text{Err}_1(u)=\Delta x^d\sum_{K\in\Omega}|u_K-u_K^{\text{ref}}|,\qquad \text{Err}_2(u)=\left(\Delta x^d\sum_{\mathbf{x}\in\Omega}|u_K-u_K^{\text{ref}}|^2\right)^{1/2},\qquad \text{Err}_{\infty}(u)=\max_{K\in\Omega}|u_K-u_K^{\text{ref}}|,		
	\end{align*}
	where \(u_K^{\text{ref}}\) denotes the reference solution at the node \(K\).
	
	\subsection{Convergence order and discontinuous solution capturing verification}\label{chapnum-ex1}
	In this example, we consider the following one-dimensional (1D) and two-dimensional (2D) Burgers' equations:
	\begin{align*} 
		&\partial_tu+u\partial_xu=0 \quad \text{with}\quad u_0(x)=\frac{1}{2}\sin(2\pi x),\qquad x\in[0,1],\\
		&\partial_tu+u\partial_xu+u\partial_yu=0 \quad \text{with}\quad u_0(x,y)=\frac{1}{4\pi}\sin(2\pi x)\sin(2\pi y),\qquad (x,y)\in[0,1]^2.
	\end{align*}

	Since the upper and lower bounds of these two problems cross zero,  in the simulations, the equilibrium distribution functions and source distribution functions are taken as \eqref{chapfirst-upwind-feq-delta} and \eqref{chapsecond-source-delta}, with \(\delta = 0.1\) for the 1D case and \(\delta = 1/(10\pi)\) for the 2D case, respectively.
	
\begin{figure}[h]
		\centering 
		\includegraphics[width=5.4in]{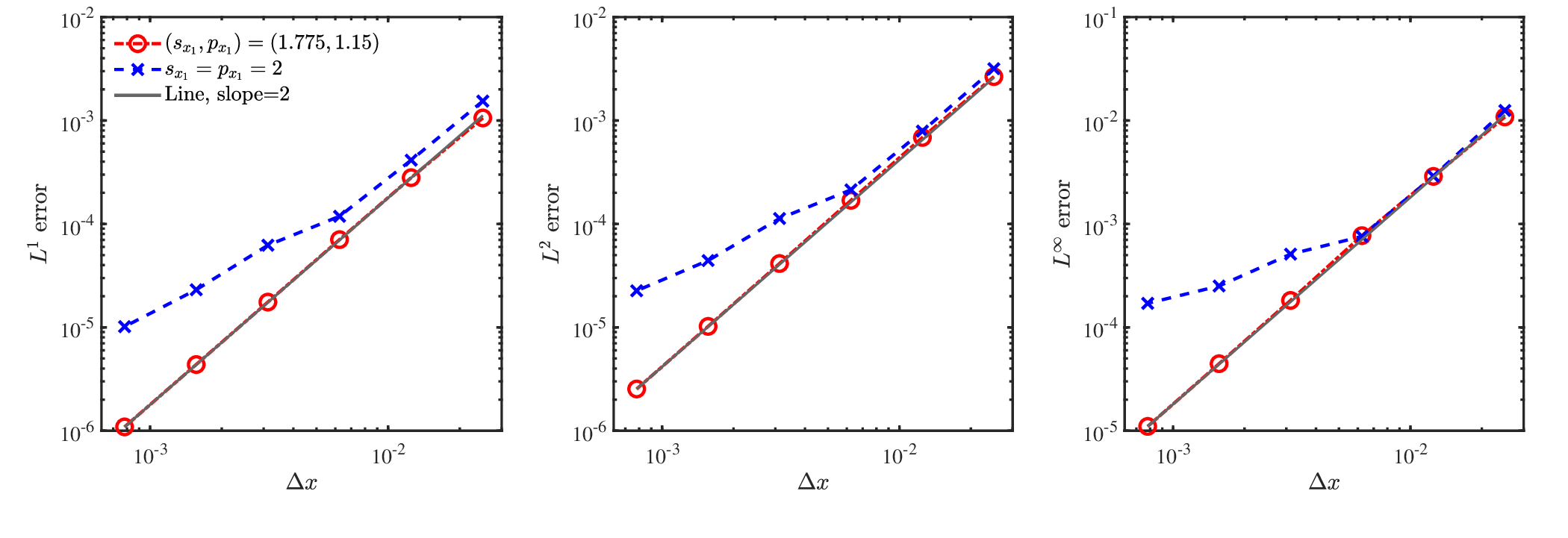}   
		\caption{\(L^1\) (left), \(L^2\) (middle), and \(L^\infty\) (right) errors of the hybrid LB method.}
		\label{ex1-1d-my-seq2-convergence}
	\end{figure}

\begin{figure}[h]
		\centering 
		\includegraphics[width=4.0in]{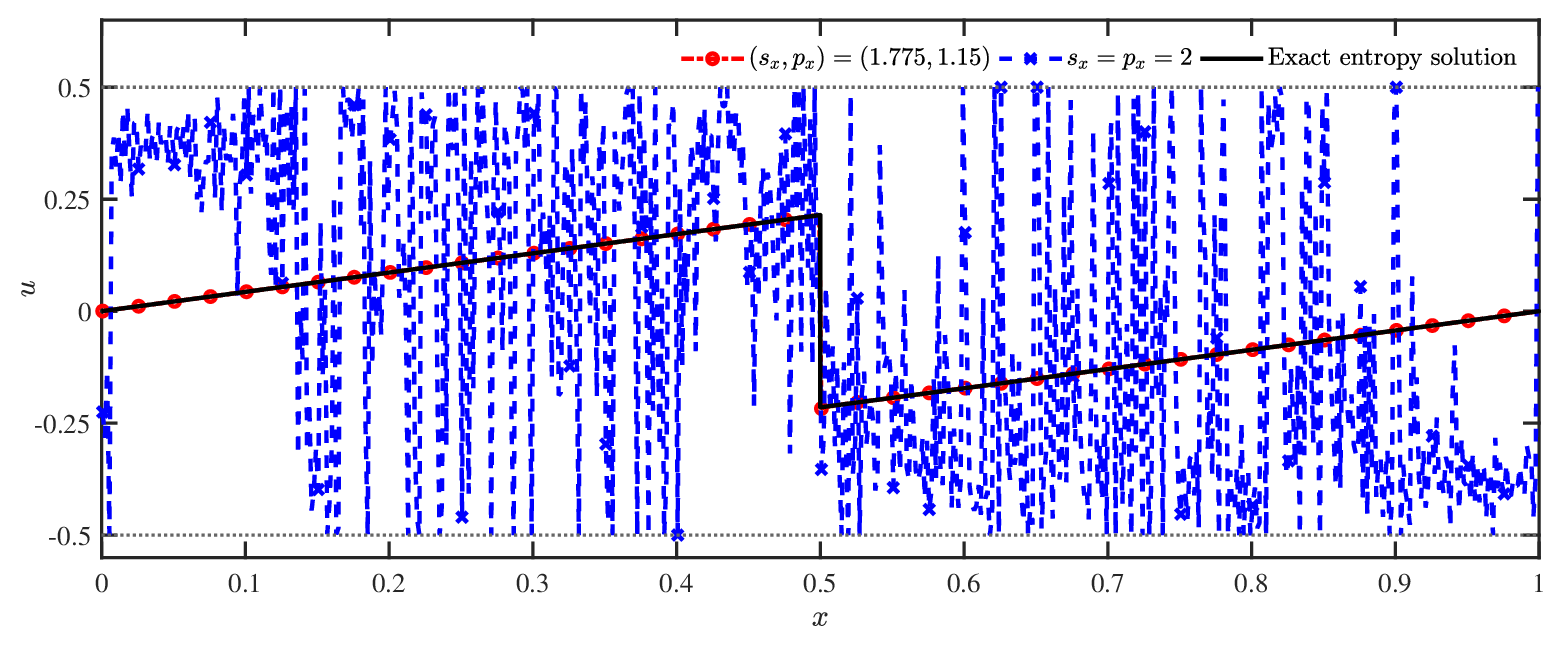}  
		\caption{Numerical solutions of the BP hybrid LB method.}
		\label{ex1-1d-my-seq2-profiles}
	\end{figure}

	For the 1D case, the initial solution lies in \([-1/2, 1/2]\), and the first gradient singularity occurs at \(t_\ast = 1/\pi\). We take  \(c = 2\) and \(t = 0.2 < t_\ast\) to test the convergence order on smooth solutions. Figure \ref{ex1-1d-my-seq2-convergence} compares the convergence of the hybrid LB method for two relaxation parameter sets: \((s_{x_1}, p_{x_1}) = (1.775, 1.15)\) and \((2, 2)\), where the GBP limiter is considered. The first set yields a convergence rate close to second-order under the \(L^1\), \(L^2\), and \(L^\infty\) norms, while the second fails to maintain consistent second-order convergence.
	We further compare the numerical solutions at \(t = 2.0\) (after shock formation) in Figure \ref{ex1-1d-my-seq2-profiles}, where \(\Delta x = 1/640\). The present method captures the shock structure well and agrees closely with the exact entropy solution, which is constructed via the Hopf-Lax potential of the Burgers' equation \cite{Lax1957Conservation,Evans2010PDE}. In contrast, the method with \(s_{x_1} = p_{x_1} = 2\) still keeps the solution in \([-1/2, 1/2]\), but produces clear non-physical oscillations. This confirms that the proposed method with adjustable relaxation parameters in the MRT model is more stable than setting all relaxation parameters to 2,  and that a proper choice of relaxation parameters remains essential for both accuracy and shock resolution.

\begin{figure}[h]
		\centering 
		\includegraphics[width=4.5in]{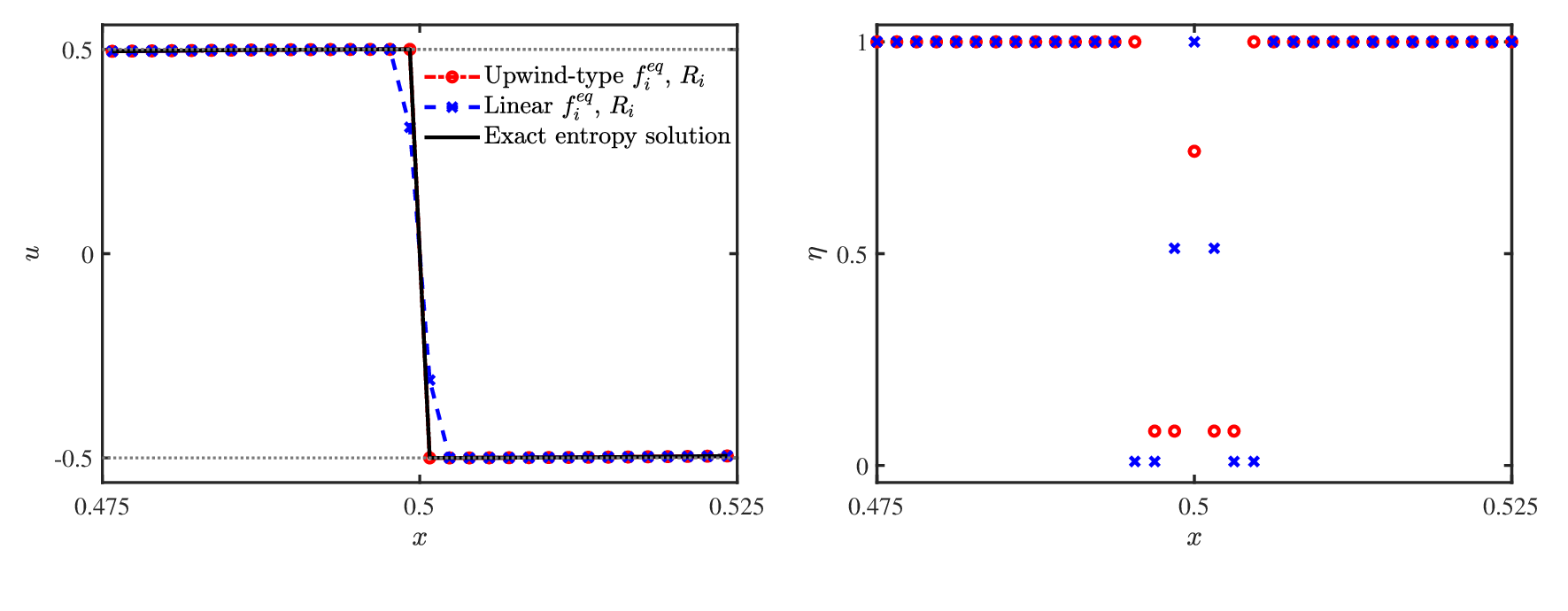}    
		\caption{Numerical solutions (left) and distribution of interface coefficients \(\eta\) (right) obtained by the hybrid LB method.}
		\label{ex1-1d-my-feq-profiles}
	\end{figure}

	In addition, we conduct a comparative simulation between the upwind-type [see  \eqref{chapfirst-upwind-feq-delta}-\eqref{chapsecond-source-delta}] and linear froms [see \eqref{chapfirst-linear-feq}-\eqref{chapsecond-source-linear}] within the hybrid LB method framework. Figure \ref{ex1-1d-my-feq-profiles} presents the comparison results for both numerical solutions and distributions of interface coefficient \(\eta\) at the time \(t = 0.5>t_*\). One can observe that, in the region near the shock at \(x = 0.5\), the GBP hybrid LB method with upwind-type \(f_i^{\text{eq}}\) and \(R_i\) outperforms the one with linear forms. What is more, around \(x=0.5\), the upwind-type forms yields larger value of  \(\eta\) than the linear forms, and also has fewer interface coefficients below 1. This indicates that after shock formation at \(t = 0.5\), the proposed hybrid LB method with upwind-type \(f_i^{\text{eq}}\) and \(R_i\) requires fewer corrections from the low-order LB method.


	For the 2D case, we set  \(c = 4\) and  \((s_{x_1}, s_{x_2}, p_{x_1}, p_{x_2}) = (1.06, 1.06, 1, 1)\). The convergence rates of the present hybrid LB method in the \(L^1\), \(L^2\), and \(L^\infty\) norms are tested here for three limiters: GBP, LBP, and RLBP (\(C_{\text{relax}}=0.001\)).	
	Table \ref{chapnum-ex1-table-errors-present} summarizes the errors and convergence orders at the time \(t=1.0\).	Specifically, GBP and RLBP yield nearly second-order convergence in the \(L^1\), \(L^2\), and \(L^\infty\) norms, which indicates that both limiers preserve the second-order accuracy for smooth solutions.	Under the stricter LBP limiter, \(L^1\) and \(L^2\) convergence remains near second-order, while the \(L^\infty\) order drops slightly to around 1.6. This is consistent with the limitation that strict local extremum constraints impose on high-order corrections near smooth extrema. In contrast, the moderately relaxed RLBP limiter restores nearly second-order \(L^\infty\) convergence.
	
\begin{figure}[h]
		\centering
		\subfloat[{$y=159/640$}]
		{\begin{minipage}{0.3\linewidth}
				\centering
				\includegraphics[width=\linewidth]{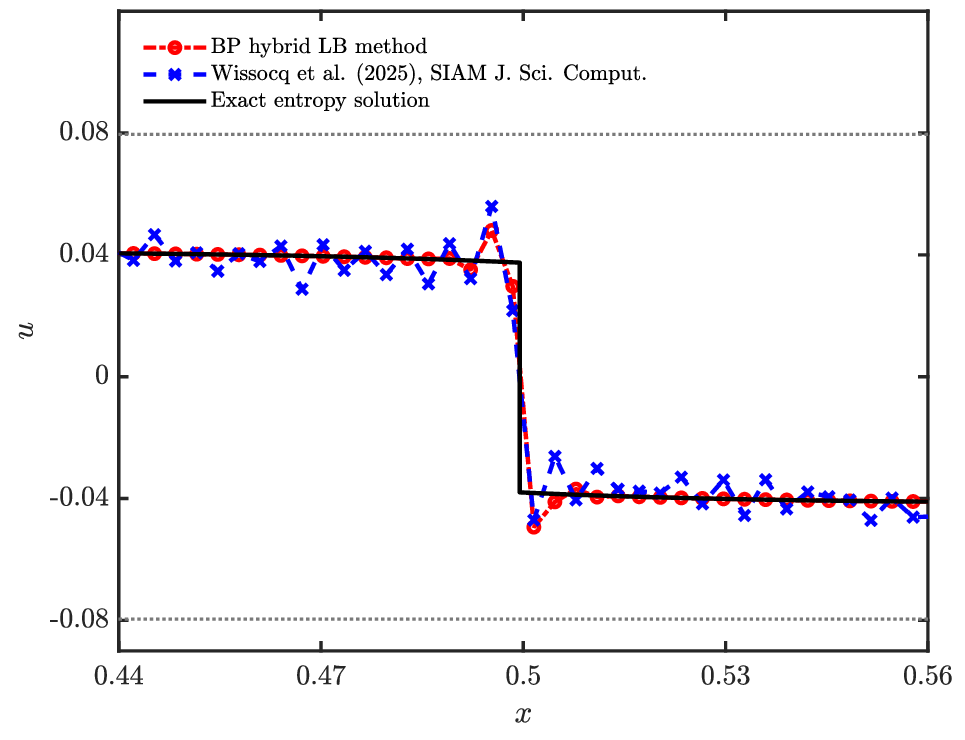}  
		\end{minipage}}   
		\subfloat[{$x=159/640$ } ]
		{\begin{minipage}{0.3\linewidth}
				\centering
				\includegraphics[width=\linewidth]{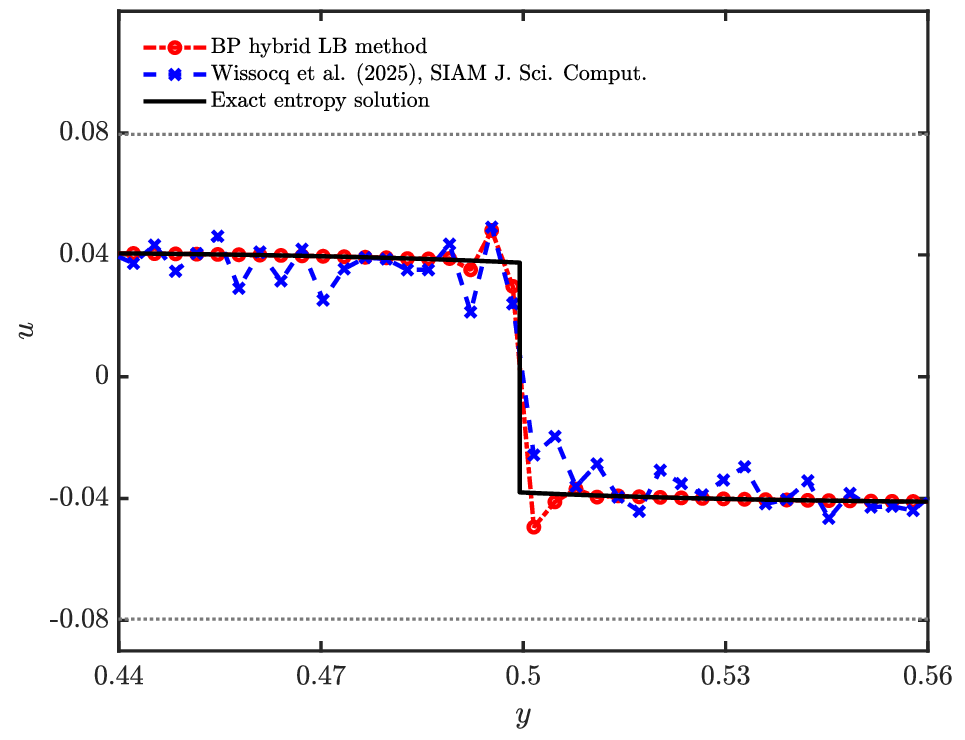} 
		\end{minipage}}     
		\caption{Profiles of the present hybrid LB method and the method of Wissocq et al. \cite{WissocqLiuAbgrall2025} under the GBP limiter.}
		\label{chapnum-ex1-mp}
	\end{figure}

\begin{figure}[h]
		\centering
		\subfloat[{$y=159/640$}]
		{\begin{minipage}{0.3\linewidth}
				\centering
				\includegraphics[width=\linewidth]{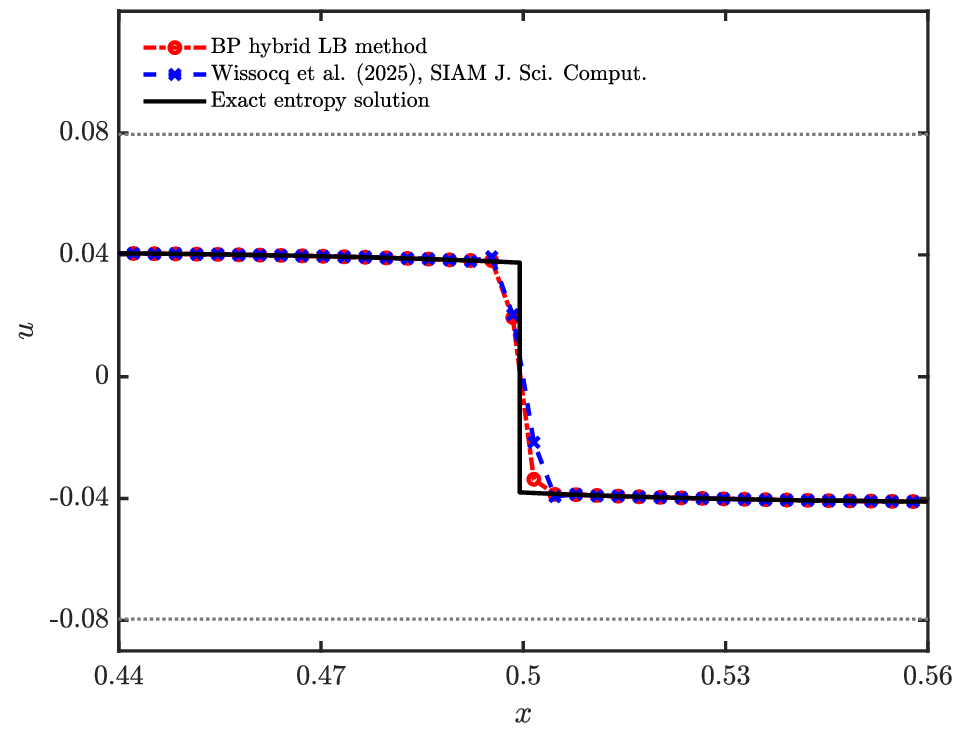}   
		\end{minipage}}   
		\subfloat[{$x=159/640$ } ]
		{\begin{minipage}{0.3\linewidth}
				\centering
				\includegraphics[width=\linewidth]{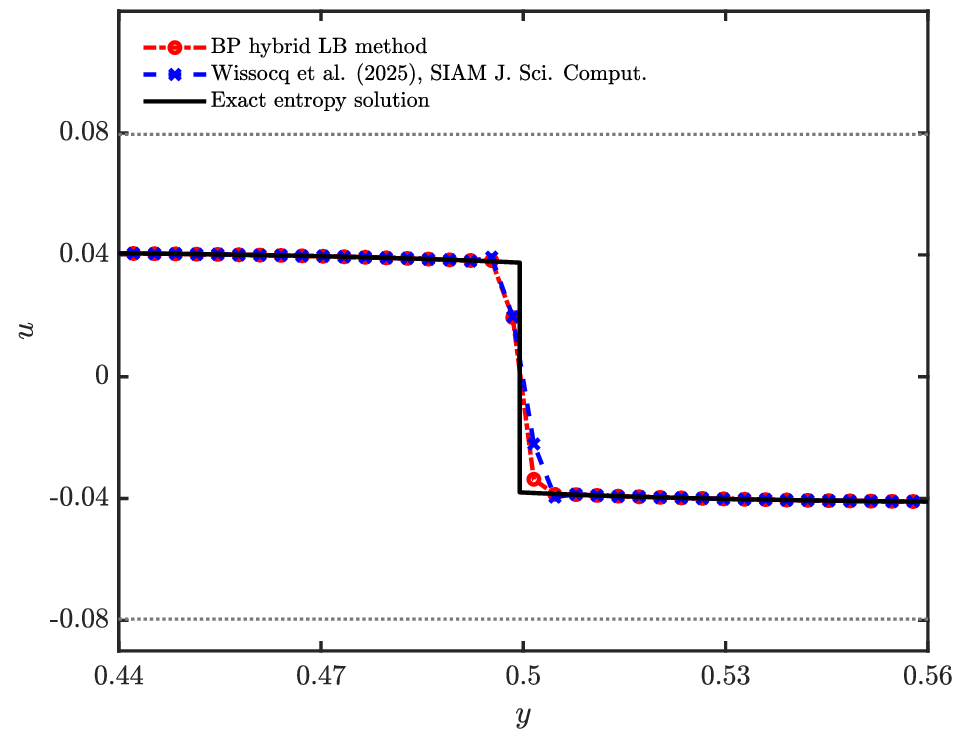}	
		\end{minipage}}     
		\caption{Profiles of the present hybrid LB method and the method of Wissocq et al. \cite{WissocqLiuAbgrall2025} under the LBP limiter.}
		\label{chapnum-ex1-lmp}
	\end{figure}

\begin{figure}[h]
		\centering
		\subfloat[{$y=159/640$}]
		{\begin{minipage}{0.3\linewidth}
				\centering
				\includegraphics[width=\linewidth]{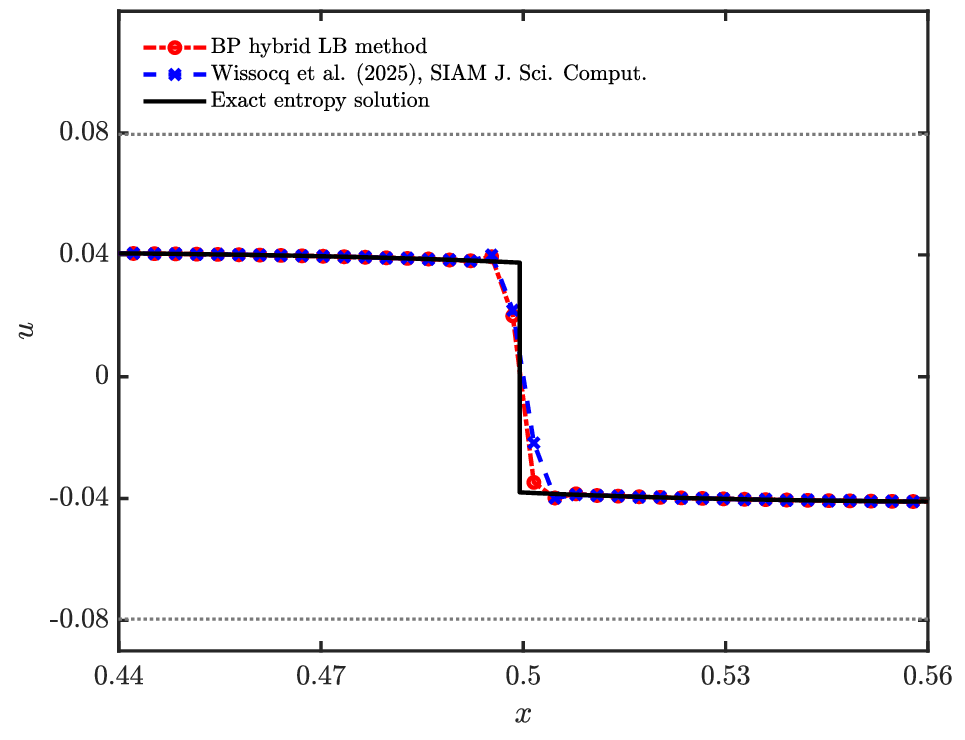}   
		\end{minipage}}   
		\subfloat[{$x=159/640$ } ]
		{\begin{minipage}{0.3\linewidth}
				\centering
				\includegraphics[width=\linewidth]{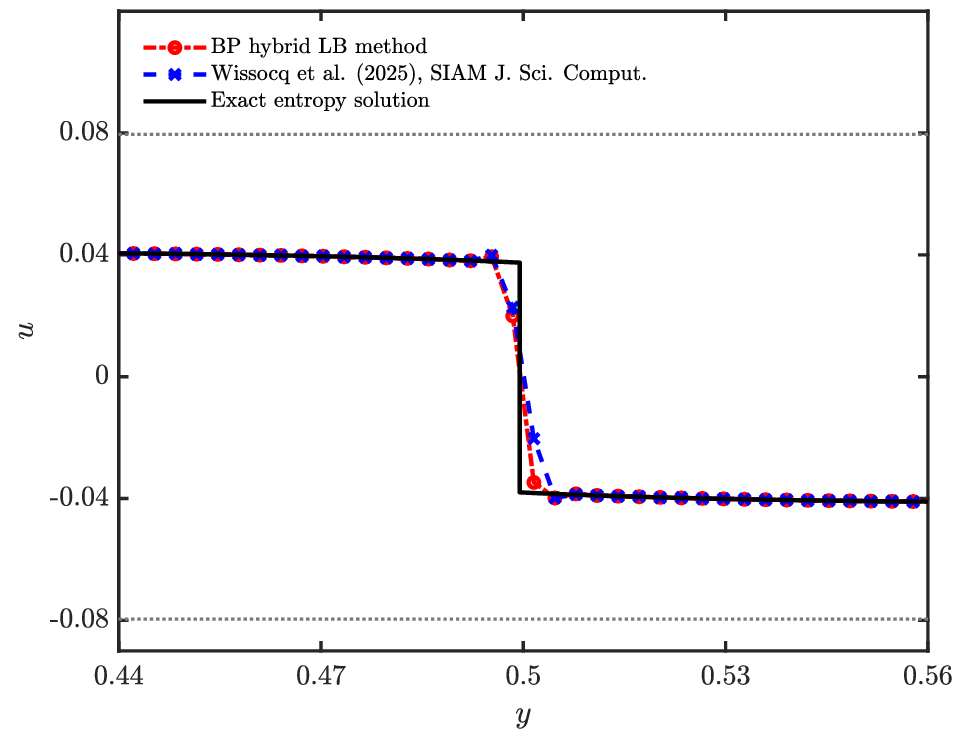}	
		\end{minipage}}    
		\caption{Profiles of the present hybrid LB method and the method of Wissocq et al. \cite{WissocqLiuAbgrall2025} under the RLBP limiter.}
		\label{chapnum-ex1-rlmp}
	\end{figure}

	Finally, under the lattice spacing \(\Delta x = 1/320\), we plot the local profiles of the present hybrid LB method and the scalar scheme of Wissocq et al. \cite{WissocqLiuAbgrall2025} at \(t = 4.0\)  in Figures \ref{chapnum-ex1-mp}, \ref{chapnum-ex1-lmp}, and \ref{chapnum-ex1-rlmp}, corresponding to the GBP, LBP, and RLBP limiters, respectively. For the method in Ref. \cite{WissocqLiuAbgrall2025}, the equilibrium distribution functions are designed as
	\begin{align*}
		\begin{cases}
			f_0^{\text{eq}}=\gamma u\\
			f^{\text{eq}}_{1}=(1-\gamma)u/4+ F_1(u)/(2c),\quad f^{\text{eq}}_{1}=(1-\gamma)u/4- F_1(u)/(2c),\\
			f^{\text{eq}}_{3}=(1-\gamma)u/4+ F_2(u)/(2c),\quad f^{\text{eq}}_{4}=(1-\gamma)u/4- F_2(u)/(2c), 
		\end{cases} 
	\end{align*}
	with  \(\gamma \in [1/2, 1)\) fixed at \(1/2\) for all examples below. And the BP condition under this setting requires \(c \ge 4\).
	
	As shown in Figure \ref{chapnum-ex1-mp}, both methods can capture the main features of the solution. In particular, the present method yields profiles in better agreement with the exact entropy solution near steep transition regions. However,  the method of Wissocq et al. \cite{WissocqLiuAbgrall2025} exhibits noticeable oscillations and local deviations in several regions, which arises primarily from its linear equilibrium distribution functions \(f_i^{\text{eq}}\) and a second-order SRT-LB method with \(\mathbf{S}=2\mathbf{I}_{2\times 2}\). For the numerical results under LBP and RLBP limiters shown in Figures \ref{chapnum-ex1-lmp} and \ref{chapnum-ex1-rlmp}, one can see that the two limiters provide stronger oscillation suppression than the GBP limiter, and the present method also outperforms the method of Wissocq et al. \cite{WissocqLiuAbgrall2025}. Overall, our proposed BP hybrid LB method can maintain boundedness of macroscopic quantity while providing good local waveform resolution.
		\begin{table}[H]    
		\begin{center}   
			\caption{\(L^1\), \(L^2\), and \(L^\infty\) errors of the present hybrid LB method under the GBP, LBP, and RLBP limiters.}  	 
			\label{chapnum-ex1-table-errors-present}
			\begin{tabular}{lccccccccccccccccc}\toprule
				Limiter&\(\Delta x\)&\(\text{Err}_1(u)\)&order&\(\text{Err}_2(u)\)&order&\(\text{Err}_{\infty}(u)\)&order\\
				\midrule
				\multirow{5}{*}{GBP}&1/20&4.6919\(\times 10^{-4}\)&-&8.0807\(\times 10^{-4}\)&-&2.8325\(\times 10^{-3}\)&-\\
				&1/40&1.2574\(\times 10^{-4}\)&1.8997&2.1993\(\times 10^{-4}\)&1.8774&8.9872\(\times 10^{-4}\)&1.6561\\
				&1/80&3.1922\(\times 10^{-5}\)&1.9779&5.5089\(\times 10^{-5}\)&1.9972&2.1975\(\times 10^{-4}\)&2.0320\\
				&1/160&8.0388\(\times 10^{-6}\)&1.9895&1.3708\(\times 10^{-5}\)&2.0067&4.9011\(\times 10^{-5}\)&2.1647\\
				&1/320&2.0177\(\times 10^{-6}\)&1.9943&3.4181\(\times 10^{-6}\)&2.0038&1.2158\(\times 10^{-5}\)&2.0112\\\hline
				\multirow{5}{*}{LBP}&1/20&5.8847\(\times 10^{-4}\)&-&9.8370\(\times 10^{-4}\)&-&3.5210\(\times 10^{-3}\)&-\\
				&1/40&1.4135\(\times 10^{-4}\)&2.0577&2.4394\(\times 10^{-4}\)&2.0117&1.0989\(\times 10^{-3}\)&1.6799\\
				&1/80&3.4217\(\times 10^{-5}\)&2.0464&5.9925\(\times 10^{-5}\)&2.0253&3.8687\(\times 10^{-4}\)&1.5061\\
				&1/160&8.5725\(\times 10^{-6}\)&1.9969&1.4806\(\times 10^{-5}\)&2.0170&1.2366\(\times 10^{-4}\)&1.6454\\
				&1/320&2.1809\(\times 10^{-6}\)&1.9748&3.7702\(\times 10^{-6}\)&1.9734&4.3429\(\times 10^{-5}\)&1.5097\\\hline
				\multirow{5}{*}{RLBP}&1/20&4.8225\(\times 10^{-4}\)&-&8.2218\(\times 10^{-4}\)&-&2.8397\(\times 10^{-3}\)&-\\
				&1/40&1.2574\(\times 10^{-4}\)&1.9393&2.1993\(\times 10^{-4}\)&1.9024&8.9872\(\times 10^{-4}\)&1.6598\\
				&1/80&3.1922\(\times 10^{-5}\)&1.9779&5.5089\(\times 10^{-5}\)&1.9972&2.1975\(\times 10^{-4}\)&2.0320\\
				&1/160&8.0388\(\times 10^{-6}\)&1.9895&1.3708\(\times 10^{-5}\)&2.0067&4.9011\(\times 10^{-5}\)&2.1647\\
				&1/320&2.0177\(\times 10^{-6}\)&1.9748&3.4181\(\times 10^{-6}\)&2.0038&1.2158\(\times 10^{-5}\)&2.0112\\
				\bottomrule
			\end{tabular}
		\end{center}  
	\end{table}

	\subsection{Interface capturing, bound preservation, and mass conservation verification}\label{chapnum-ex2}
	In the second example, we evaluate the performance of our proposed BP hybrid LB method for large-deformation interface advection. Two classic benchmark cases are considered: the 2D single vortex and a spherical interface in a three-dimensional (3D) nonlinear shear flow, both governed by the following pure advection equation:
	\begin{align}\label{chapnum-ex2-2d-eq}
		\partial_t\phi+\mathbf{v}\cdot\nabla\phi=0.
	\end{align}
	\subsubsection{2D single vortex}\label{chapnum-ex2-2d}
	For this 2D test, the single-vortex advection problem is considered \cite{Liang2023InterfaceCapturing}. The computational domain is \(\Omega = [0, L]^2\) with \(L = 200\), and periodic boundary conditions are imposed. The initial condition is a circular disk centered at \((x_c, y_c) = (100, 150)\) with radius \(R = 40\), whose continuous profile reads
	\begin{align*}
		\phi_0(x, y) = \mathbf{1}_{(x-100)^2 + (y-150)^2 < 40^2},
	\end{align*}
	the prescribed velocity field is given by
	\begin{align*}
		\begin{cases}
			v_x= U_0 \sin^2\left(\frac{\pi x}{L}\right) \sin\left(\frac{2\pi y}{L}\right) \cos\left(\frac{\pi t}{T}\right), \\
			v_y= -U_0 \sin^2\left(\frac{\pi y}{L}\right) \sin\left(\frac{2\pi x}{L}\right) \cos\left(\frac{\pi t}{T}\right). 
		\end{cases}
	\end{align*}
	We set \(U_0 = 0.04\), which yields an advective time scale \(T_f = L/U_0 = 5000\) and a total simulation time \(T = 6T_f = 30000\). The velocity vanishes at \(t = T/2\), after which the cosine factor becomes negative and the flow reverses direction automatically. This indicates that the effective advection time is given by
	\begin{align}\label{chapnum-ex2-reftime}
		\tau(t) = \int_0^t \cos(\pi s/T)\,ds = \frac{T}{\pi}\sin\left(\frac{\pi t}{T}\right). 
	\end{align}
Accordingly, the interface reaches its maximum deformation at \(T/2\), and returns to its initial position at \(t = T\).
	
	Simulations are performed on a uniform \(200 \times 200\) grid using the \(D2Q5\) lattice structure, with the lattice spacing \(\Delta x = 1\) and the lattice velocity \(c = 1\). The sufficient BP condition in \eqref{eq:facecondition} is well satisfied under this setup. For the hybrid LB method, the GBP limiter is adopted, and the relaxation parameters in the second-order MRT-LB method are given by \((s_{x_1},s_{x_2},p_{x_1},p_{x_2})=(1.6, 1.6, 0.8, 0.8)\).  
	
	To compute the numerical errors, the characteristic equation \eqref{chapnum-ex2-2d-eq} is integrated backward using a fourth-order Runge-Kutta method \cite{Hairer1993ODE,Butcher2016ODE}, and each cell center is traced back to the initial time. The reference solution is set to 1 if the traced point lies inside the initial disk and 0 otherwise. At the final time \(t = T\), the flow returns to its initial shape, so the initial profile can be used directly as the reference solution.

\begin{figure}[h]
		\centering 
		\includegraphics[width=5.4in]{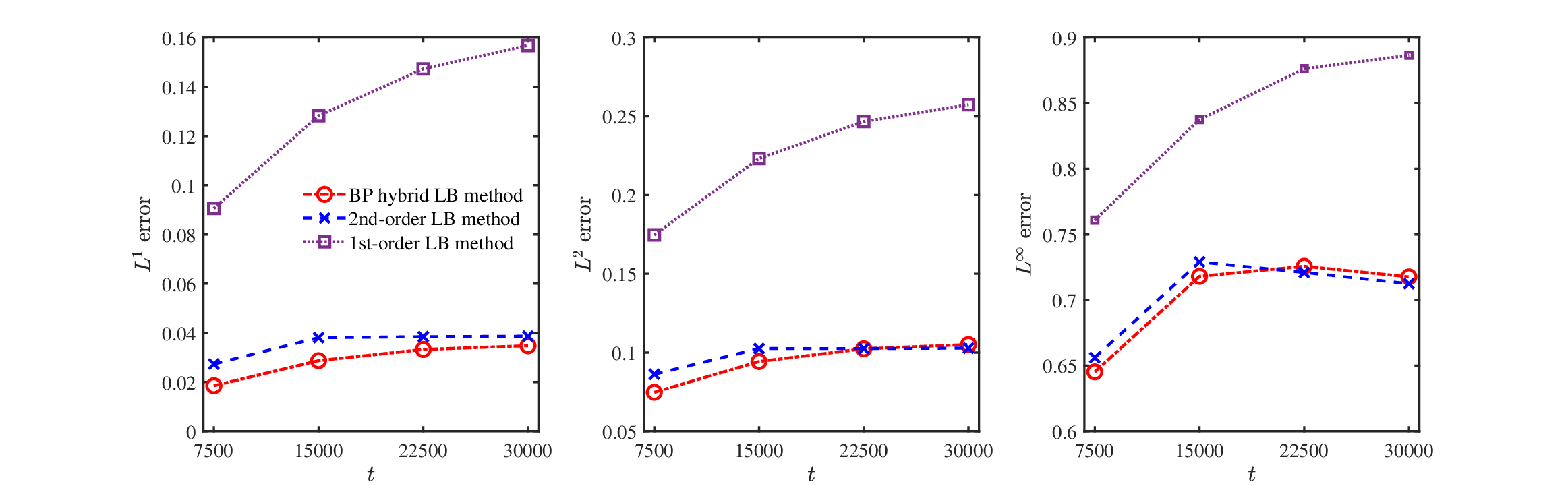}   
		\caption{\(L^1\) (left), \(L^2\) (middle), and \(L^\infty\) (right) errors of the BP hybrid, the first-order, and the  second-order LB methods.}
		\label{chapnum-ex2-2d-error}
	\end{figure}

	Figure \ref{chapnum-ex2-2d-error} presents the \(L^1\), \(L^2\), and \(L^\infty\) errors for the hybrid , the first-order, and the second-order LB methods at four values of the time: \(T/4\), \(T/2\), \(3T/4\), and \(T\). The results demonstrate that both the hybrid method and the second-order method achieve substantially higher numerical accuracy than the first-order method. Notably, the error magnitudes of the BP hybrid LB method and the second-order LB method are highly comparable, which shows that the extra error introduced by the GBP limiter is acceptable.
	
\begin{figure}[h]
		\centering
		\subfloat[{$T/4$}]
		{\begin{minipage}{0.75\linewidth}
				\centering
				\includegraphics[width=\linewidth]{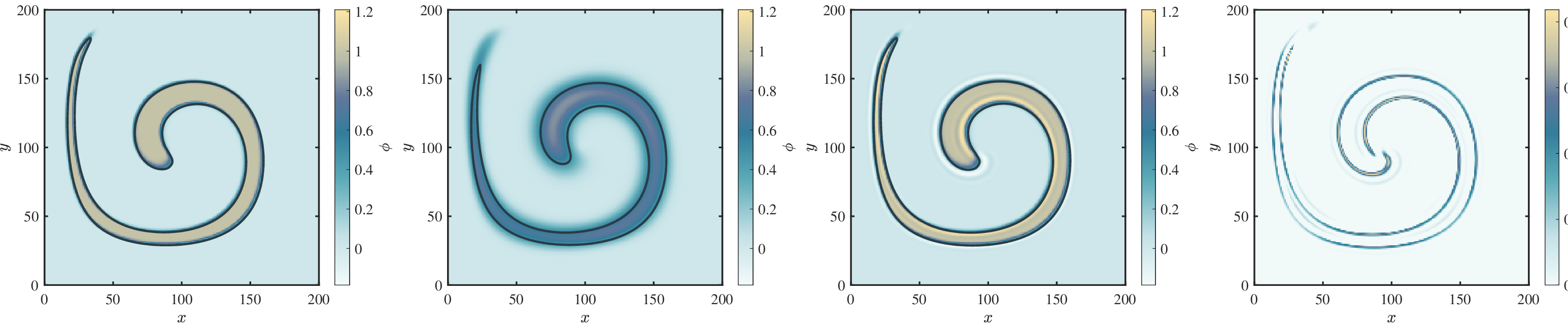}   
		\end{minipage}}   
	
		\subfloat[{$T/2$ } ]
		{\begin{minipage}{0.75\linewidth}
				\centering
				\includegraphics[width=\linewidth]{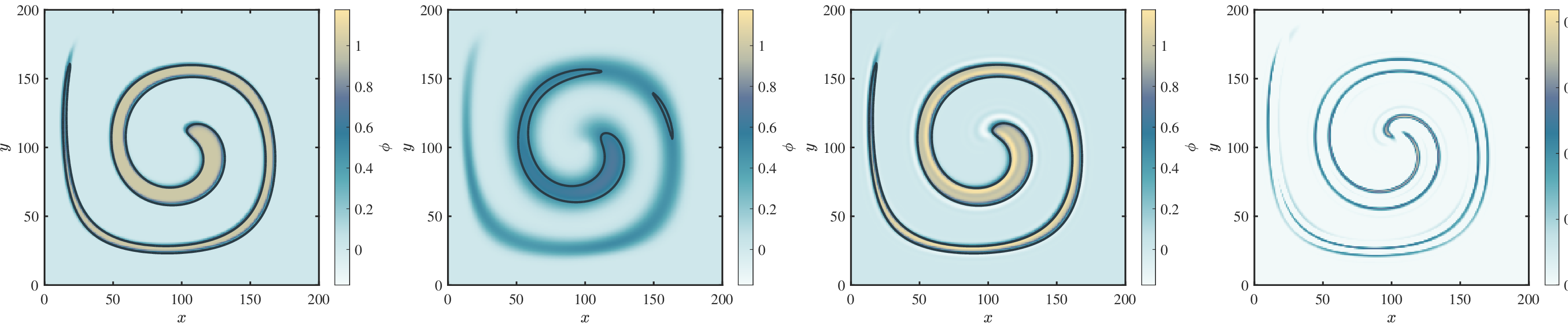}	
		\end{minipage}}    
		
		\subfloat[{$3T/4$ } ]
		{\begin{minipage}{0.75\linewidth}
				\centering
				\includegraphics[width=\linewidth]{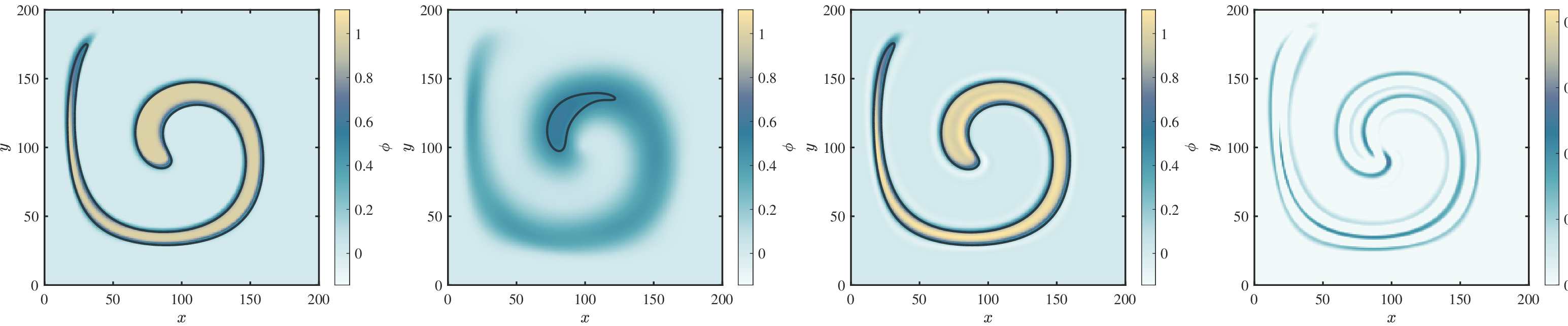}   
		\end{minipage}}   
		
		\subfloat[{$T$ } ]
		{\begin{minipage}{0.75\linewidth}
				\centering
				\includegraphics[width=\linewidth]{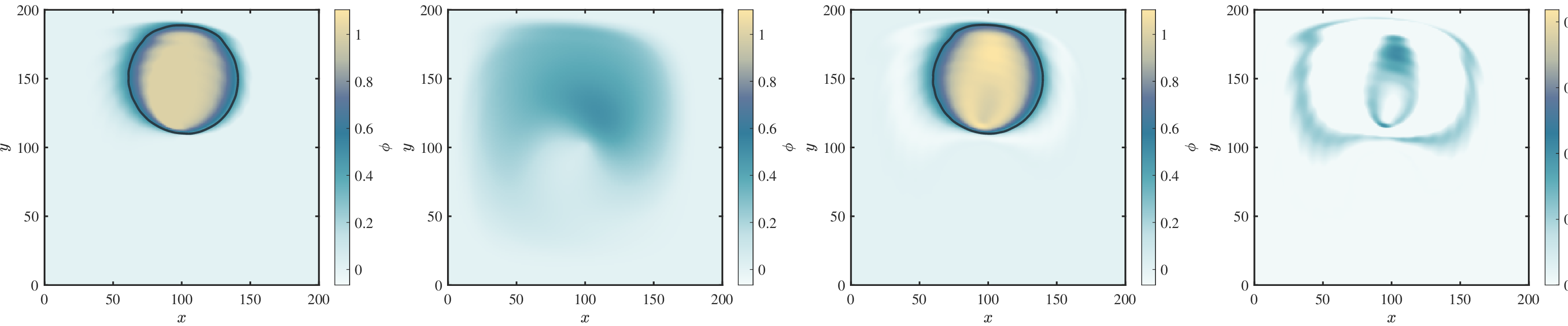} 
		\end{minipage}}   
		
		\caption{2D scalar fields. From left to right: the BP hybrid LB method,  the first-order LB method, the second-order LB method, and the distribution of \(\max(0, -\phi, \phi-1)\) for the second-order method.}
		\label{chapnum-ex2-fourtime-field}
	\end{figure}

\begin{figure}[h]
		\centering\par\vspace{-5mm}
		\subfloat[{$\min \phi$}]
		{\begin{minipage}{0.3\linewidth}
				\centering
				\includegraphics[width=\linewidth]{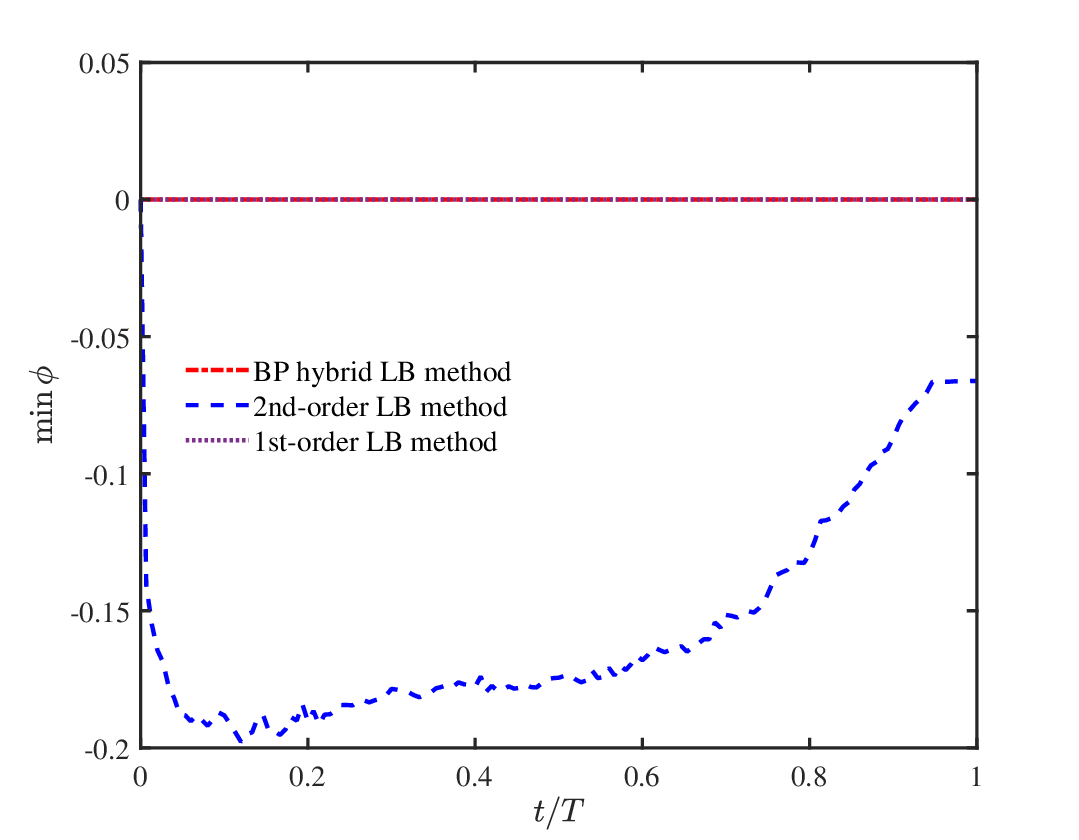}   
		\end{minipage}}   
		\subfloat[{$\max \phi$ } ]
		{\begin{minipage}{0.3\linewidth}
				\centering
				\includegraphics[width=\linewidth]{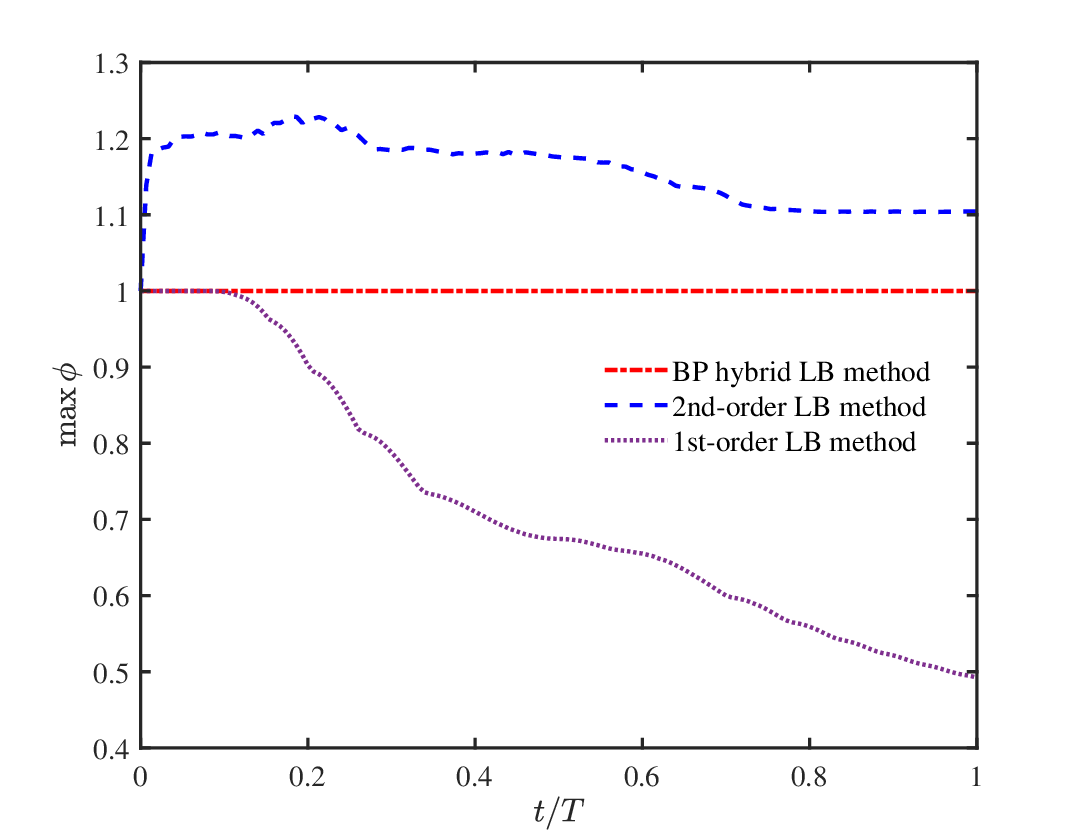}	
		\end{minipage}}    
		
		\caption{Minimum and maximum values of $\phi$ for the BP hybrid,  the first-order, and the second-order LB methods.}
		\label{chapnum-ex2-minmax-phi}
	\end{figure}

	Figure \ref{chapnum-ex2-fourtime-field} presents the scalar fields from the BP hybrid, the first-order, and the second-order LB methods at several times. For the second-order method, the quantity \(\max(0, -\phi, \phi-1)\) is also shown to measure bound violations. Figure \ref{chapnum-ex2-minmax-phi} plots the minimum and maximum values of \(\phi\) for all three methods over the full simulation. As seen in Figure \ref{chapnum-ex2-fourtime-field}, the initial disk is stretched and twisted into a long spiral during the first half period. The flow then reverses in the second half period, and the interface slowly unwinds back toward its initial shape. Both the BP hybrid and second-order methods keep the main spiral structure with relatively sharp interfaces. Moreover, Figure \ref{chapnum-ex2-minmax-phi} confirms that the BP hybrid method keeps the solution strictly between 0 and 1 throughout the simulation. In contrast, the second-order method shows pronounced overshoots and undershoots near the interface, as evidenced in both Figures \ref{chapnum-ex2-fourtime-field} and \ref{chapnum-ex2-minmax-phi}. The first-order method also keeps the solution strictly between 0 and 1, but it significantly smears the transition layer and weakens the interface along the thin filament. Notably, its peak \(\phi\) drops to about 0.5 at \(t = T\) in Figure \ref{chapnum-ex2-minmax-phi}. As a result, no contour at \(\phi = 0.5\) is visible for the first-order method at the final time in Figure \ref{chapnum-ex2-fourtime-field}(d).
	
\begin{figure}[h]
		\centering 
		\includegraphics[width=2.5in]{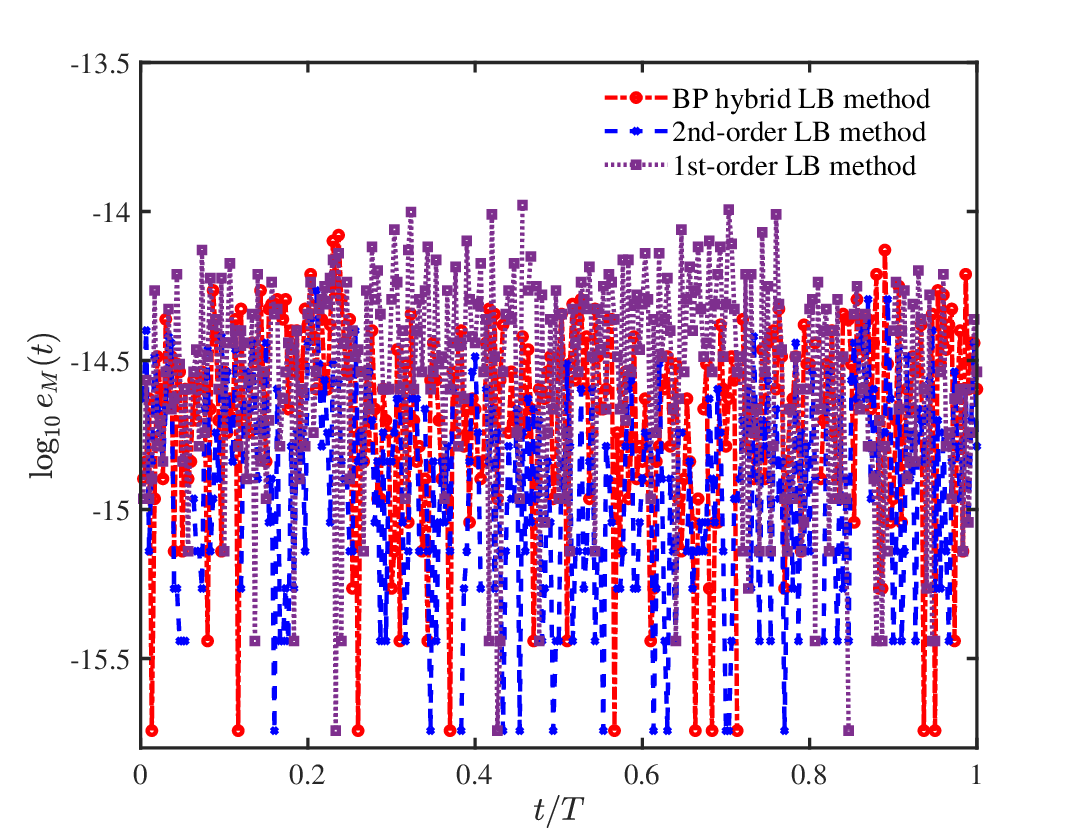}   
		\caption{Relative mass $e_M$ of the BP hybrid,  the first-order, and the second-order LB methods.}
		\label{chapnum-ex2-2d-mass}
	\end{figure}

	Figure \ref{chapnum-ex2-2d-mass} presents the relative mass errors for all three LB methods over the full simulation period. The discrete total mass and the corresponding relative mass error are defined as
	\begin{align*}
		M_{\Delta x}(t) = \Delta x^d \sum_{K \in \Omega} u_K(t), \qquad e_M = \frac{|M_{\Delta x}(t) - M_{\Delta x}(0)|}{M_{\Delta x}(0)}.
	\end{align*}
	As shown in the figure, the BP hybrid, the first-order, and the second-order LB methods all maintain mass conservation, which is consistent with the theoretical analysis.

	\subsubsection{3D sphere deformation in a shear flow}
	Now we consider a 3D spherical interface advection test \cite{Liang2023InterfaceCapturing}. The computational domain is \(\Omega = [0,L]^3\) with \(L = 100\). The initial condition corresponds to a sphere of radius 20 centered at \((30, 30, 50)\), given by
	\begin{align*}
		\phi_0(x, y, z) = \mathbf{1}_{(x-30)^2+(y-30)^2+(z-50)^2 < 20^2}.
	\end{align*}
	
	Let \(\chi_\alpha = \pi(x_\alpha/L - 1/2)\). The velocity field \(\mathbf{v}\) is given by
	\begin{align*}
		v_x &= U_0\pi \cos\chi_x \left(\sin\chi_z - \sin\chi_y\right) \cos(\pi t/T),  \\
		v_y &= U_0\pi \cos\chi_y \left(\sin\chi_x - \sin\chi_z\right) \cos(\pi t/T),  \\
		v_z &= U_0\pi \cos\chi_z \left(\sin\chi_y - \sin\chi_x\right) \cos(\pi t/T). 
	\end{align*}
		
	We set \(U_0 = 0.02\), which gives an advective time scale \(T_f = L/U_0 = 5000\) and a total simulation time \(T = 2T_f = 10000\). The time-dependent factor is the same as in the 2D case, as defined in \eqref{chapnum-ex2-reftime}, so the flow reverses at the half period and the exact solution returns to its initial spherical shape at the final time.	
	Simulations are performed on a uniform \(100^3\) grid using the \(D3Q7\) lattice structure, with \(\Delta x = 1\) and \(c = 1\). The sufficient BP condition \eqref{eq:facecondition} is then well satisfied.  For the hybrid LB method, the GBP limiter is  considered, and the relaxation parameters in the second-order LB method are set to \((s_{x_1}, s_{x_2}, s_{x_3}, p_{x_1}, p_{x_2}, p_{x_3}) = (1.6, 1.6, 1.6, 1, 1, 1)\).

\begin{figure}[h]
		\centering
		\includegraphics[width=5.4in]{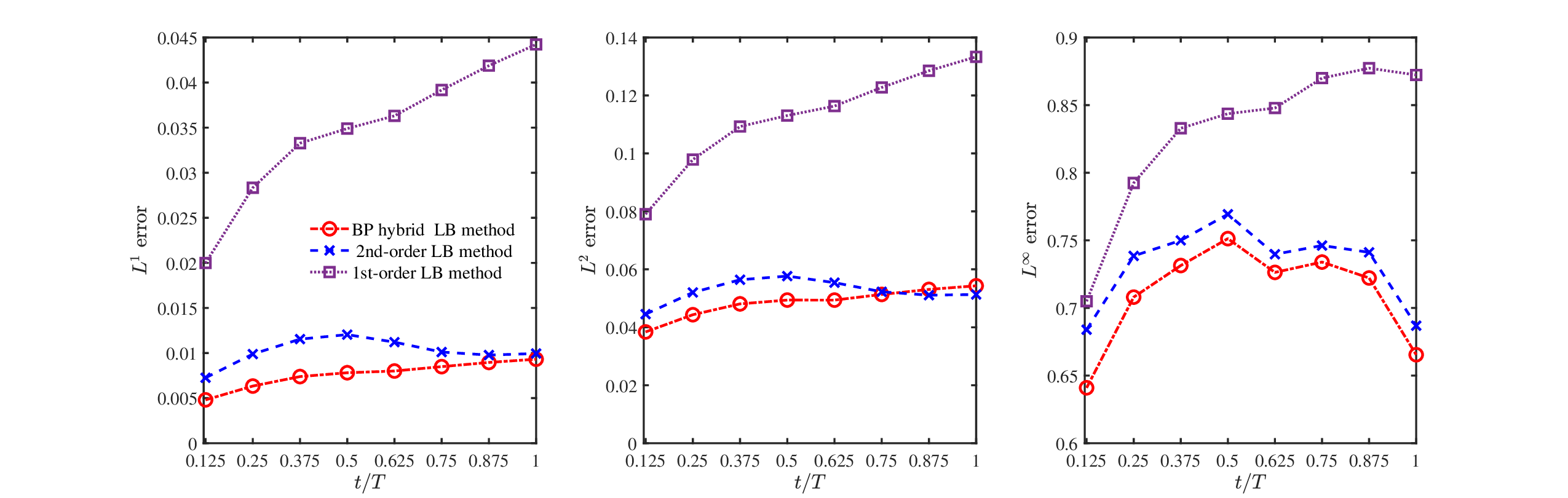}   
		\caption{\(L^1\) (left), \(L^2\) (middle), and \(L^\infty\) (right) errors of the BP hybrid, the first-order, and the  second-order LB methods.}
		\label{chapnum-ex2-3d-error}
	\end{figure}

	Similar to the 2D test case in section \ref{chapnum-ex2-2d}, Figure \ref{chapnum-ex2-3d-error} shows the \(L^1\), \(L^2\), and \(L^\infty\) errors for the BP hybrid, the first-order, and the second-order LB methods at the eight different values of the time: \(T/8\), \(T/4\), \(3T/8\), \(T/2\), \(5T/8\), \(3T/4\), and \(T\). At all time points, both the BP hybrid and second-order methods have much higher numerical accuracy than the first-order  method. In particular, the BP hybrid method generally has smaller errors than the second-order method for most of the simulation.

\begin{figure}[h]
		\centering
		\includegraphics[width=5.4in]{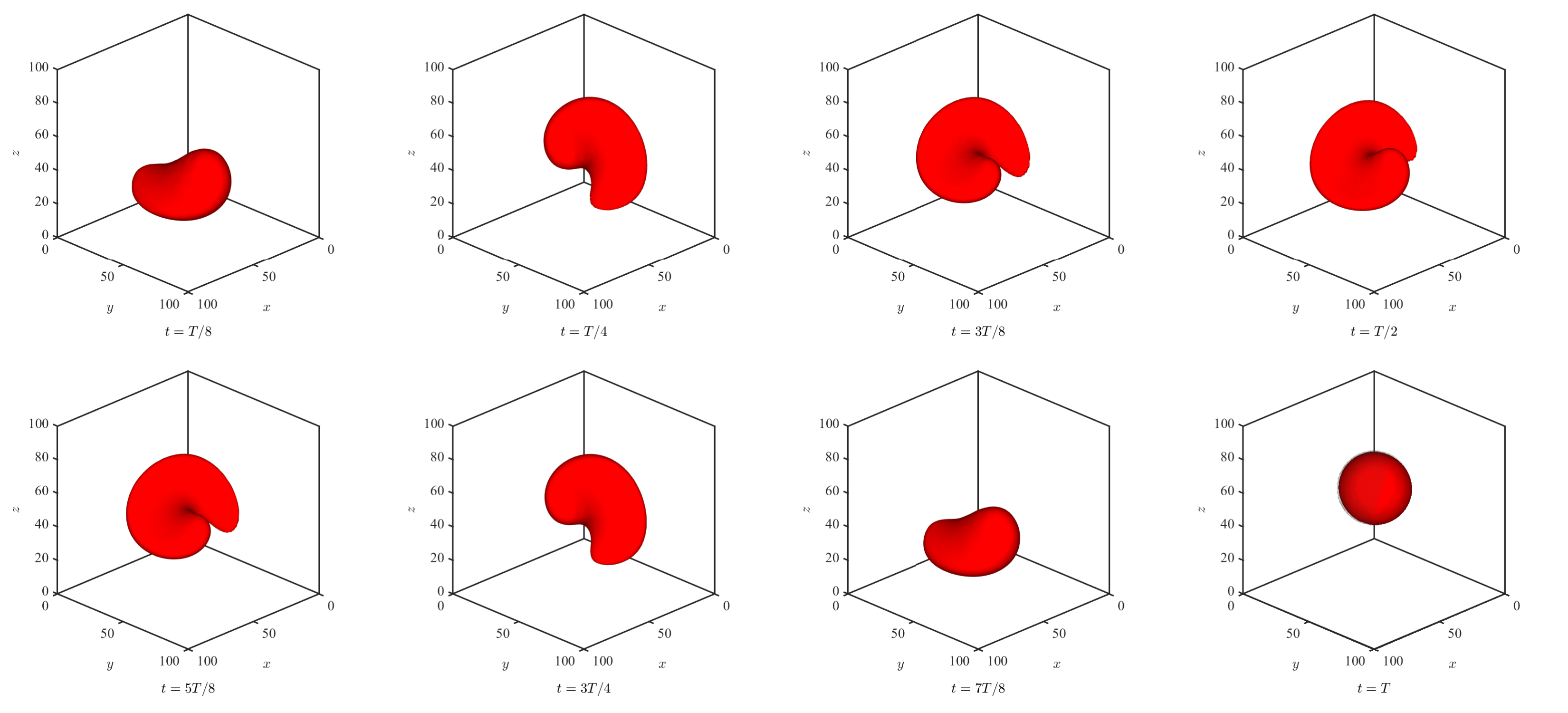}  
		\caption{3D isosurfaces at \(\phi = 0.5\) of the BP hybrid LB method.}
		\label{chapnum-ex2-eighttime-field-hybrid}
	\end{figure}

\begin{figure}[h]
		\centering\par\vspace{-5mm}
			\includegraphics[width=5.4in]{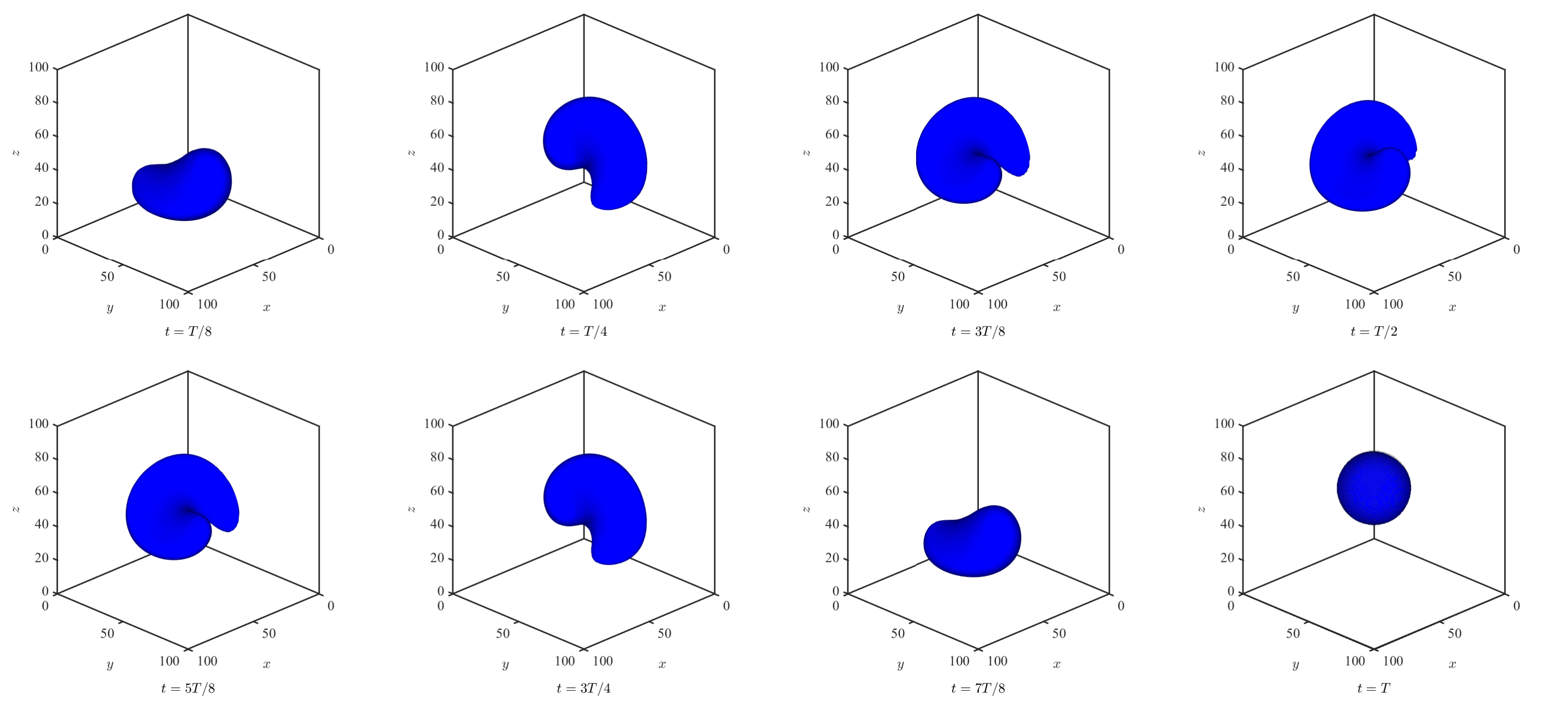}  
		\caption{3D isosurfaces at \(\phi = 0.5\) of the second-order LB method.}
		\label{chapnum-ex2-eighttime-field-high}
	\end{figure}

\begin{figure}[h]
		\centering
		\includegraphics[width=5.4in]{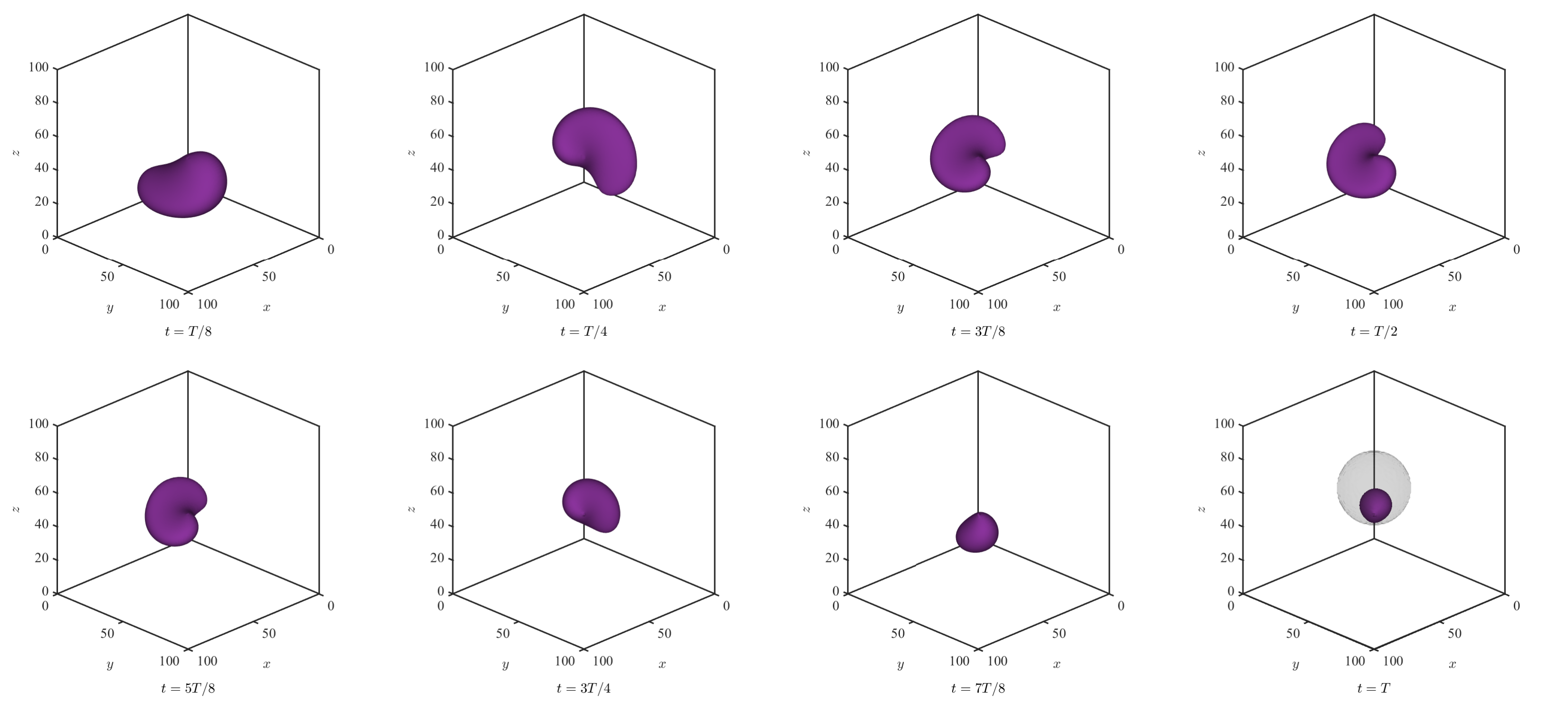}  
		\caption{3D isosurfaces at \(\phi = 0.5\) of the first-order LB method.}
		\label{chapnum-ex2-eighttime-field-low}
	\end{figure}

\begin{figure}[h]
	 	\centering\par\vspace{-5mm}
	 	\subfloat[{$\min \phi$}]
	 	{\begin{minipage}{0.3\linewidth}
	 			\centering
	 			\includegraphics[width=\linewidth]{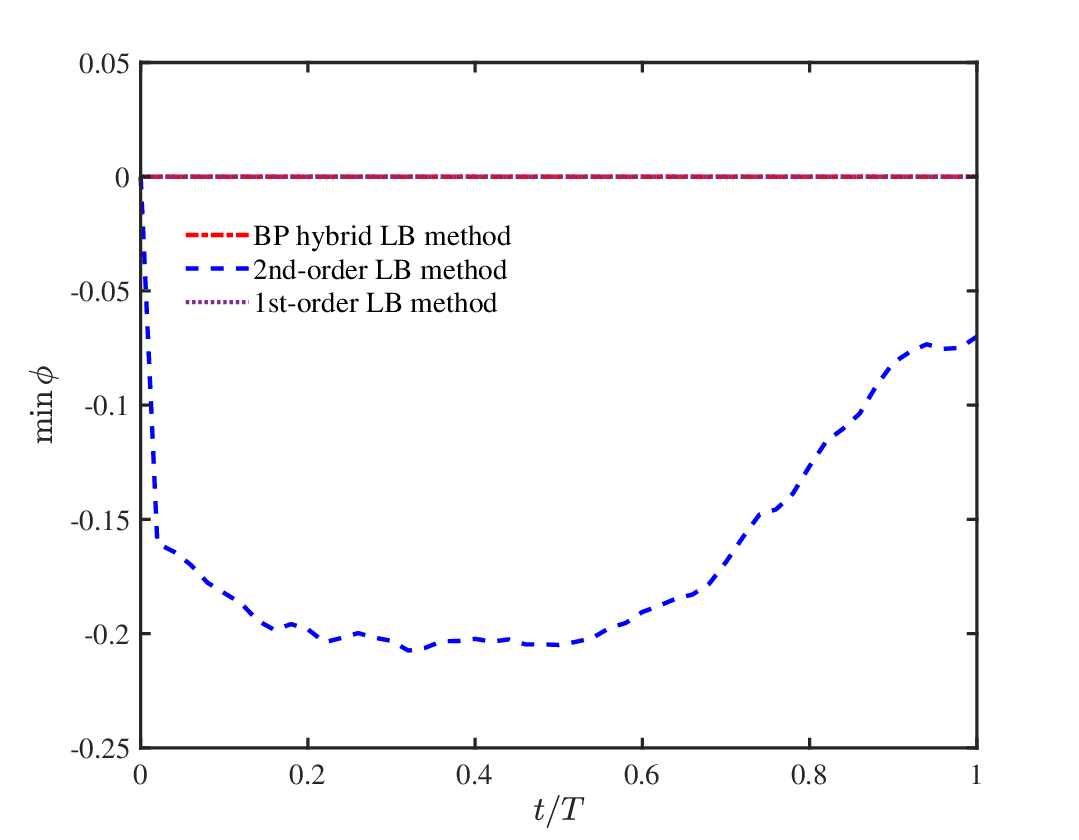}   
	 	\end{minipage}}   
	 	\subfloat[{$\max \phi$ } ]
	 	{\begin{minipage}{0.3\linewidth}
	 			\centering
	 			\includegraphics[width=\linewidth]{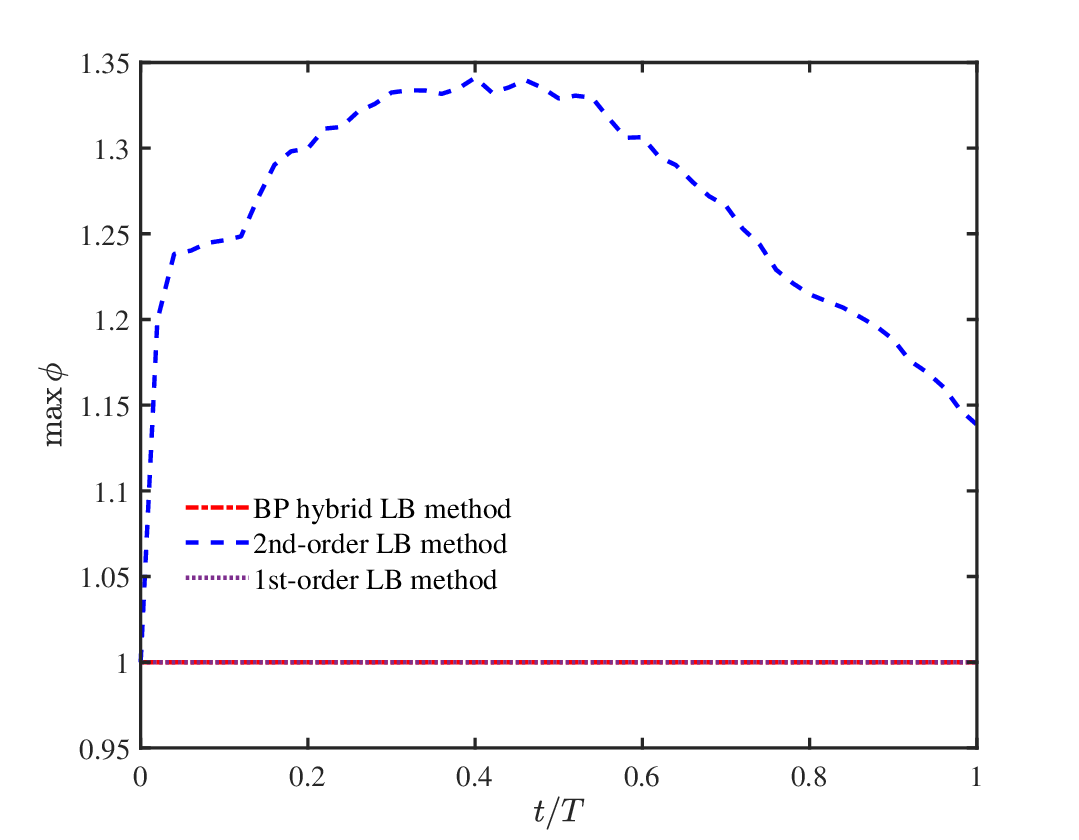}	
	 	\end{minipage}}    
	 	
	 	\caption{Minimum and maimum  values of $\phi$ for the BP hybrid,  the first-order, and the second-order LB methods.}
	 	\label{chapnum-ex2-3d-minmax-phi}
	 \end{figure}

	Figures \ref{chapnum-ex2-eighttime-field-hybrid}, \ref{chapnum-ex2-eighttime-field-high}, and \ref{chapnum-ex2-eighttime-field-low} show the \(\phi = 0.5\) isosurfaces from the three LB methods over the full simulation. The spherical interface is gradually stretched, bent, and twisted by the shear flow, reaching maximum deformation at \(T/2\). When the flow reverses in the second half period, the interface slowly returns to its initial shape. Figures \ref{chapnum-ex2-eighttime-field-hybrid} and \ref{chapnum-ex2-eighttime-field-high} show that the BP hybrid and second-order methods produce isosurfaces with very similar overall shapes. By contrast, the first-order method in Figure \ref{chapnum-ex2-eighttime-field-low} gives noticeably shrunk isosurfaces that lose fine details, which is due to the strong numerical diffusion.	
	We also plot the minimum and maximum values of \(\phi\) over time in Figure \ref{chapnum-ex2-3d-minmax-phi}. The BP hybrid method keeps the solution strictly between 0 and 1 throughout the simulation. However, the second-order method produces values outside \([0, 1]\) throughout the simulation. The first-order method stays within \([0, 1]\), but it introduces strong numerical dissipation, and the sphere shrinks noticeably relative to its initial shape, as also seen in Figure \ref{chapnum-ex2-eighttime-field-low}.

	 \begin{figure}[h]
	 	\centering 
	 	\includegraphics[width=2.5in]{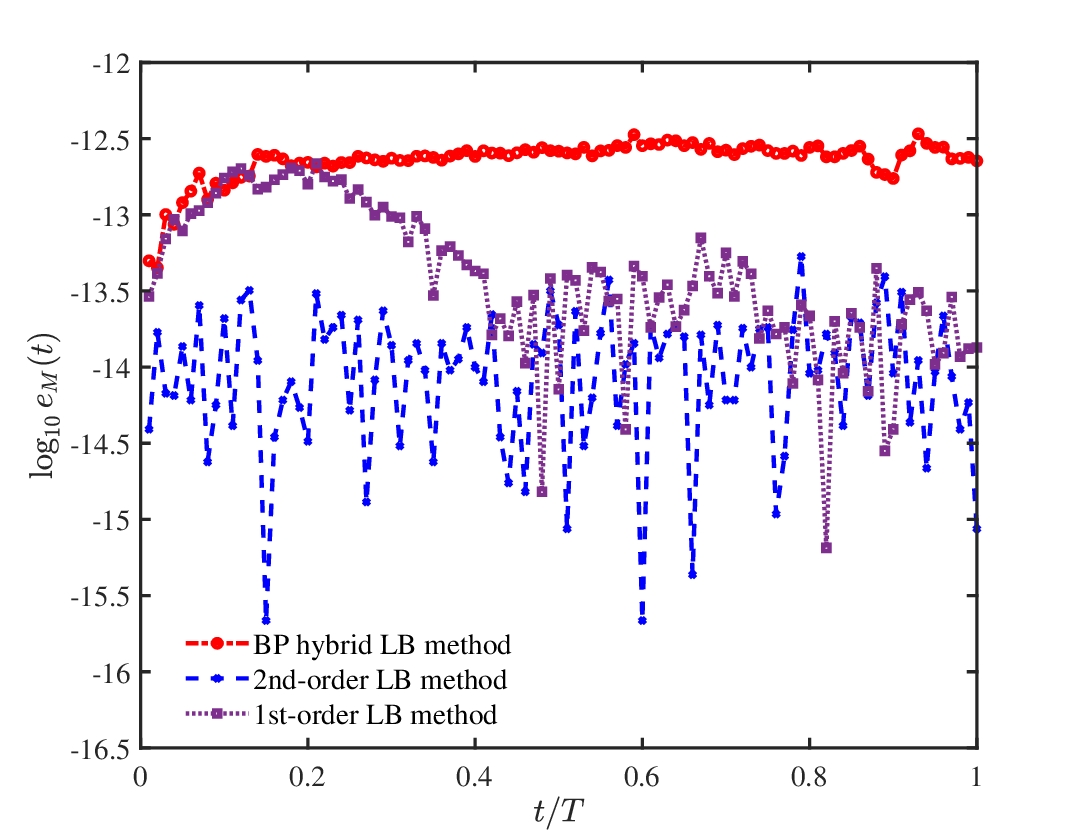}   
	 	\caption{Relative mass $e_M$ of the BP hybrid,  the first-order, and the second-order LB methods.}
	 	\label{chapnum-ex2-3d-mass}
	 \end{figure}
	 
	Figure \ref{chapnum-ex2-3d-mass} shows the relative mass errors for all three LB methods over the whole simulation period. All three methods maintain good mass conservation, which is consistent with the theoretical properties of the discrete scheme.

	In conclusion, numerical results from both the 2D and 3D test cases consistently show that our proposed BP hybrid LB method effectively suppresses nonphysical out-of-bounds behavior of macroscopic scalar quantity during long-duration reversible large-deformation advection, while maintaining mass conservation via the shared-surface flux mechanism. In particular, compared with the first-order LB method, the BP hybrid LB method significantly reduces interface deformation and numerical dissipation. Compared with the second-order LB method, the BP hybrid LB method preserves a similar main interface structure while strictly enforcing solution boundedness.
	\subsection{2D rotating wave: waveform resolution and bound preservation}\label{chapnum-ex3}
	To evaluate the numerical performance of the BP hybrid LB method under the combined effects of nonconvex fluxes and strong discontinuities, this example considers the Kurganov-Petrova-Popov (KPP) rotating wave problem \cite{Kurganov2007KPP}. 
	Specifically, we solve the following equation on the periodic computational domain \(\Omega = [-2, 2]^2\):
	\begin{align*}
		\partial_tu+\partial_x\sin u+\partial_y(\cos u -1)=0,
	\end{align*}
	with the initial condition given by
	\begin{align*}
		u_0(x, y) = 
		\begin{cases}
			7\pi/2, & x^2 + y^2 < 1, \\
			\pi/4, & \text{otherwise}.
		\end{cases}
	\end{align*}
	
In the simulations, we take \(\Delta x = 1/80\) and \(\Delta t = \Delta x/c\). The reference solution is computed from a high-resolution finite-volume run at \(2560 \times 2560\), using a Godunov flux for the nonconvex flux, minmod reconstruction, and a second-order SSP Runge-Kutta time step with \(\Delta t \le 0.18\,\Delta x_{\text{ref}}\) \cite{Gottlieb2001SSP,LeVeque2002FiniteVolume}.
	
\begin{figure}[h]
		\centering 
		\subfloat[\(y=-1/160\)]
		{\begin{minipage}{0.7\linewidth}
				\centering
				\includegraphics[width=\linewidth]{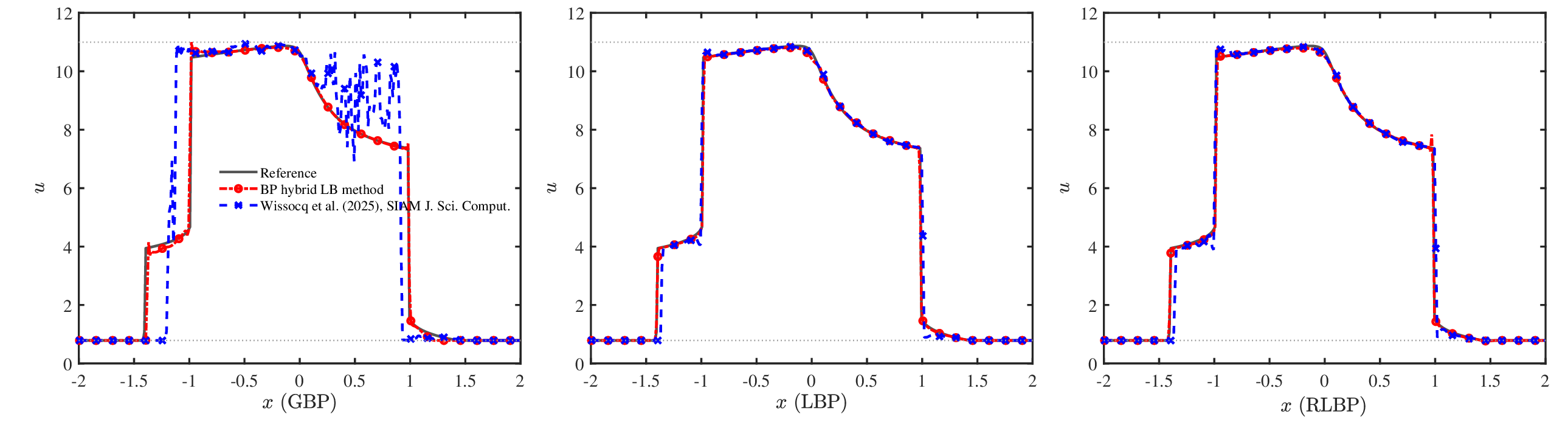}	
		\end{minipage}}    
		
		\subfloat[\(x=-1/160\)]
		{\begin{minipage}{0.7\linewidth}
				\centering
				\includegraphics[width=\linewidth]{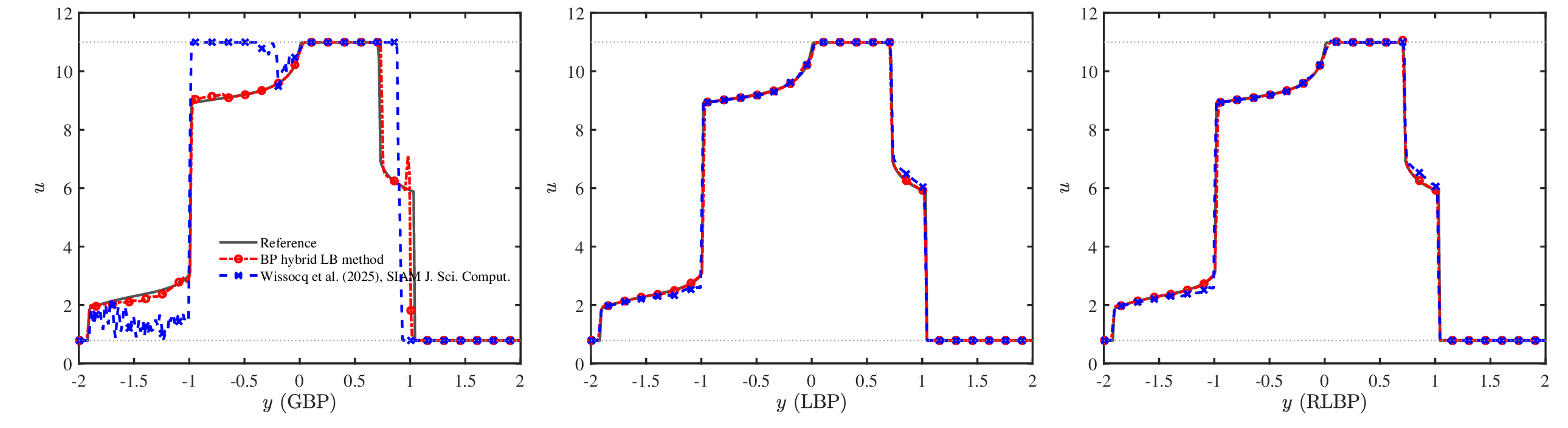}   
		\end{minipage}}    
		
		\caption{Profiles in the \(x\) (top) and \(y\) (bottom) directions at \(t = 1\) under \(c=4.0\). }
		\label{chapnum-ex3-profiles}
	\end{figure}

\begin{figure}[h]
		\centering\par\vspace{5mm}
		\subfloat[GBP  limiter]
		{\begin{minipage}{0.55\linewidth}
				\centering
				\includegraphics[width=\linewidth]{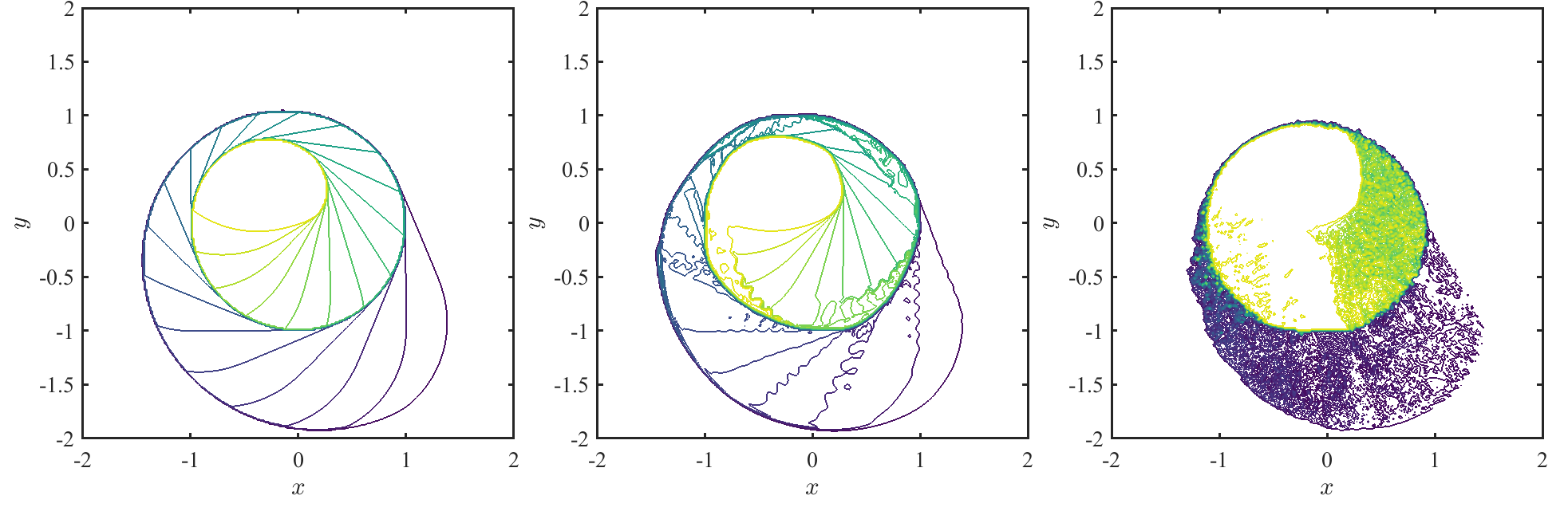}   
		\end{minipage}}   
		
		\subfloat[LBP limiter]
		{\begin{minipage}{0.55\linewidth}
				\centering
				\includegraphics[width=\linewidth]{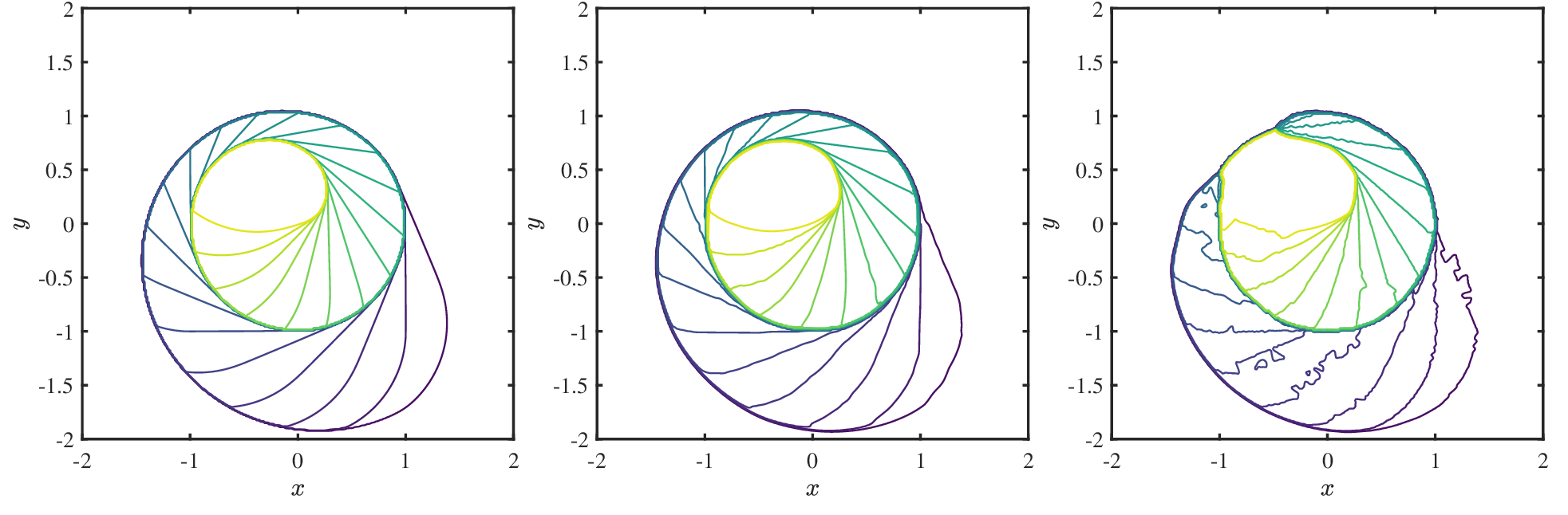}	
		\end{minipage}}    
		
		\subfloat[RLBP limiter]
		{\begin{minipage}{0.55\linewidth}
				\centering
				\includegraphics[width=\linewidth]{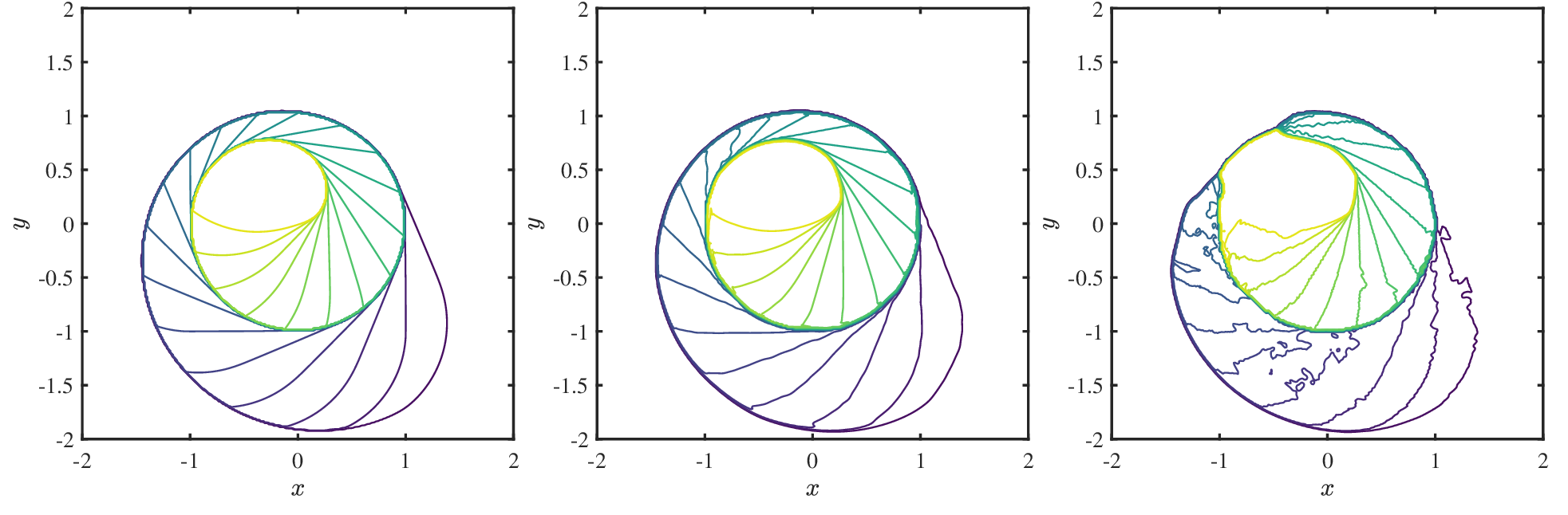}   
		\end{minipage}}    
		
		\caption{2D contour lines at \(t = 1\) under \(c = 4.0\). From left to right: the reference solution, the numerical solution of the hybrid LB method, and that of the method in Ref. \cite{WissocqLiuAbgrall2025}.}
		\label{chapnum-ex3-contours}
	\end{figure}

	This test case mainly focuses on comparing our proposed BP hybrid LB method with the scalar scheme developed by Wissocq et al. \cite{WissocqLiuAbgrall2025} based on the convex limiter framework. First, with \(\Delta x = 1/80\) and \(c = 4\), we plot the profiles and 2D contour lines of the numerical solutions obtained by the hybrid LB method under the GBP, LBP, and RLBP (\(C_{\text{relax}} = 0.001\)) limiters at \(t = 1.0\) in Figures \ref{chapnum-ex3-profiles} and \ref{chapnum-ex3-contours}. The corresponding relaxation parameters are given by
	\begin{align*}
		\begin{cases}
			(s_{x_1}^{\text{GBP}},s_{x_2}^{\text{GBP}},p_{x_1}^{\text{GBP}},p_{x_2}^{\text{GBP}})=(1.06,1,1,1), \\
			(s_{x_1}^{\text{LBP}},s_{x_2}^{\text{LBP}},p_{x_1}^{\text{LBP}},p_{x_2}^{\text{LBP}})=	(1.4,0.96875,1,1),\\
			(s_{x_1}^{\text{RLBP}},s_{x_2}^{\text{RLBP}},p_{x_1}^{\text{RLBP}},p_{x_2}^{\text{RLBP}})=	(1.126875,0.853125,1,1).
		\end{cases}
		\end{align*} 
		
	As seen in Figures \ref{chapnum-ex3-profiles} and \ref{chapnum-ex3-contours}(a), the hybrid LB method with the BP limiter can effectively control the global extrema but does not fully eliminate nonphysical oscillations within the domain. This behavior is observed in both the present method and that of Wissocq et al. \cite{WissocqLiuAbgrall2025}, but the internal structure obtained by our hybrid LB method is in better agreement with the reference solution. 
	When the LBP and RLBP limiters are applied, the contours obtained by our hybrid LB method are more continuous. Both the internal fan-shaped structure and the outer curved wavefronts show improved agreement with the reference solution. By contrast, the method of Wissocq et al. \cite{WissocqLiuAbgrall2025} still exhibits noticeable oscillations along the wavefronts. This discrepancy arises from their use of linear equilibrium distribution functions and a second-order LB method with all relaxation parameters set to 2. This indicates that, the upwind-type equilibrium distribution functions and the second-order MRT-LB method with the auxiliary source distribution functions adopted in this work can improve waveform resolution to a certain extent.
	
\begin{figure}[h]
		\centering\par\vspace{-5mm}	
		\subfloat[GBP limiter]
		{\begin{minipage}{0.7\linewidth}
				\centering
				\includegraphics[width=\linewidth]{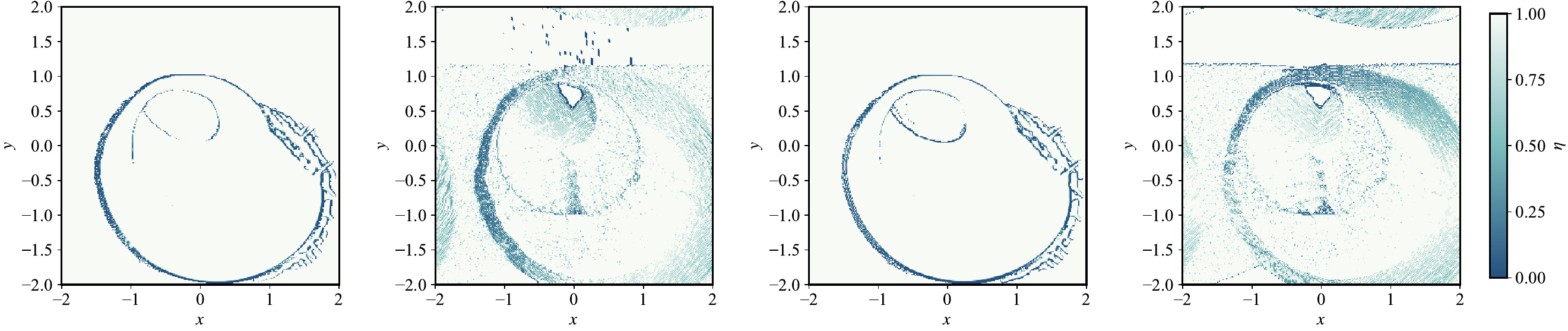}   
		\end{minipage}}   
		
		\subfloat[LBP limiter]
		{\begin{minipage}{0.7\linewidth}
				\centering
				\includegraphics[width=\linewidth]{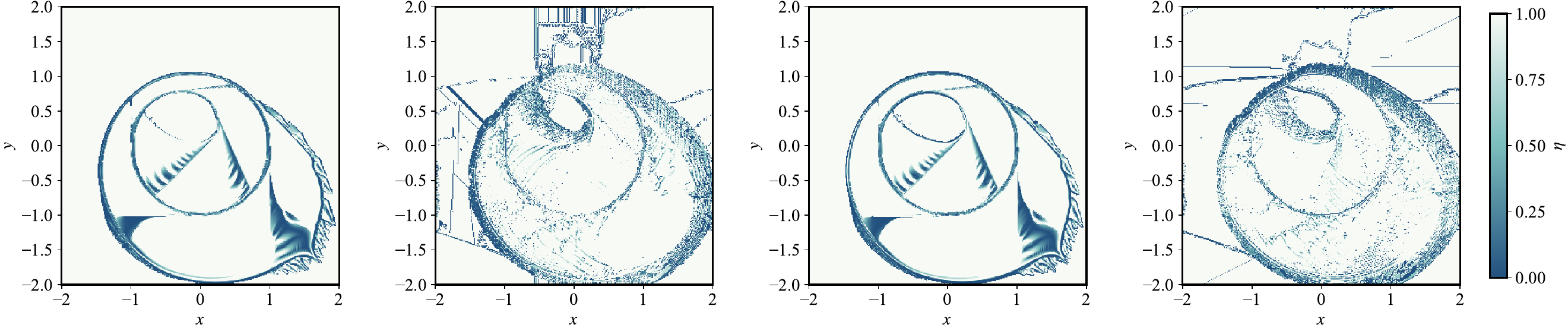}   
		\end{minipage}}    
		
		\subfloat[RLBP limiter]
		{\begin{minipage}{0.7\linewidth}
				\centering
				\includegraphics[width=\linewidth]{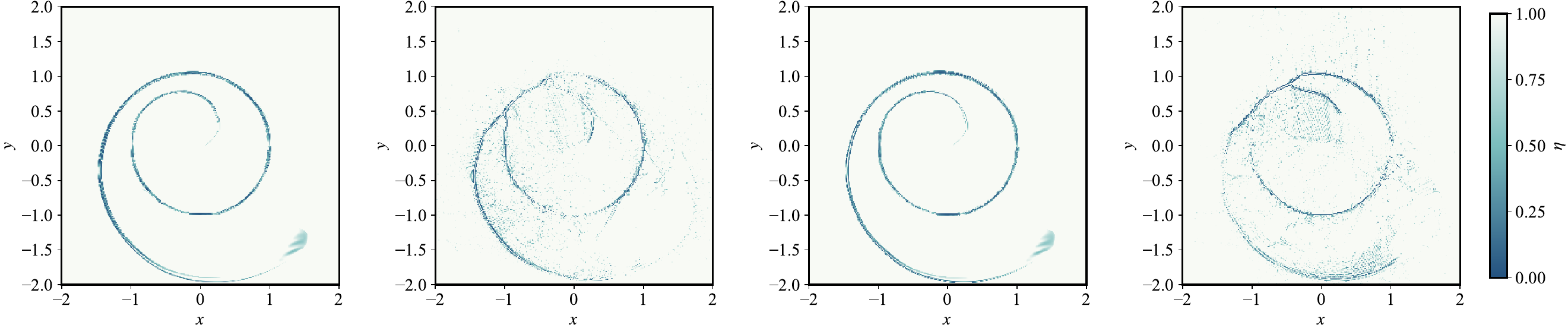}    
		\end{minipage}}    
		\par\vspace{-1mm}	
		\caption{Distributions of the interface coefficient \(\eta\) in the \(x\) and \(y\) directions at \(t = 1\). From left to right: \(x\)-direction results for the present hybrid LB method and Ref. \cite{WissocqLiuAbgrall2025}, followed by the corresponding \(y\)-direction results.}
		\label{chapnum-ex3-flux}
	\end{figure}

	In addition, for the present hybrid LB method and the method of Wissocq et al. \cite{WissocqLiuAbgrall2025} under the three convex limiters, we also present the distributions of interface coefficient \(\eta\) at the time \(t = 1\) in Figure \ref{chapnum-ex3-flux}. And the \(L^1\), \(L^2\), and \(L^\infty\) errors are lisited in Table \ref{chapnum-ex3-table-errors}, where different values of the time ( \(t = 0.1, 0.2, 0.5, 1.0,2.0\)) are considered. 	It is observed from Figure \ref{chapnum-ex3-flux} that the smaller values of \(\eta\) in the present method are mainly distributed along curved wavefronts and their adjacent regions with strong solution gradients. The coefficient drops to zero on some interfaces, indicating that high-order corrections are fully suppressed at these locations. This behavior is more pronounced under the LBP limiter, suggesting that the limiter mechanism can adapt to the local structure of the solution. When the RLBP limiter relaxes the local admissible intervals, the regions with strong restriction are significantly reduced. In contrast, the interface coefficients in the method of Wissocq et al. \cite{WissocqLiuAbgrall2025} appear as more discrete spots and scattered bands, and their alignment with the main wavefronts is less clear than in the proposed method. In particular, for \(c = 4\), where both methods satisfy the BP conditions, Table \ref{chapnum-ex3-table-errors} shows that the present method yields smaller errors than the method of Wissocq et al. \cite{WissocqLiuAbgrall2025} under all three limiters and at all five times.

	For the results in Figures \ref {chapnum-ex3-profiles}, \ref{chapnum-ex3-contours}, and Table \ref{chapnum-ex3-table-errors}, we would also like to point out that, under the parameters selected for this example, there is virtually no difference between the RLMP limiter and the LMP limiter. This means that the LMP limiter is sufficient for this example in terms of waveform resolution.
	
	Finally, to examine the BP performance of the proposed hybrid method at a lower lattice velocity, we adopt a flux-correction limiter based on the nodal BP constraints (see \ref{app1} for details). This strategy only requires the hybrid method to satisfy the sufficient BP condition of the first-order LB method, i.e., \(c \ge \sup_u \bigl(|\cos u| + |\sin u|\bigr)\) for the present example, whereas the BP condition for the first-order LB method in Refs. \cite{Chai2016,Chai2020} is \(c \ge 2\).
	
\begin{figure}[h]
		\centering\par\vspace{-5mm}
		\subfloat[\(y=-1/160 \)]	
			{\begin{minipage}{0.7\linewidth}
				\centering
				\includegraphics[width=\linewidth]{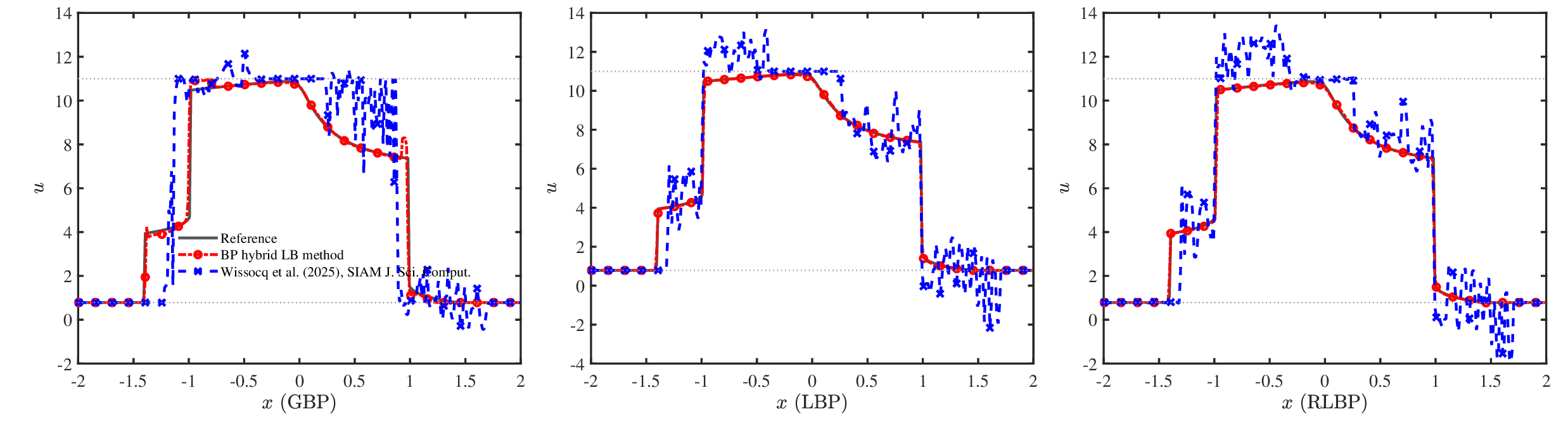}	
		\end{minipage}}    
		
		\subfloat[\(x=-1/160 \)]
		{\begin{minipage}{0.7\linewidth}
				\centering
				\includegraphics[width=\linewidth]{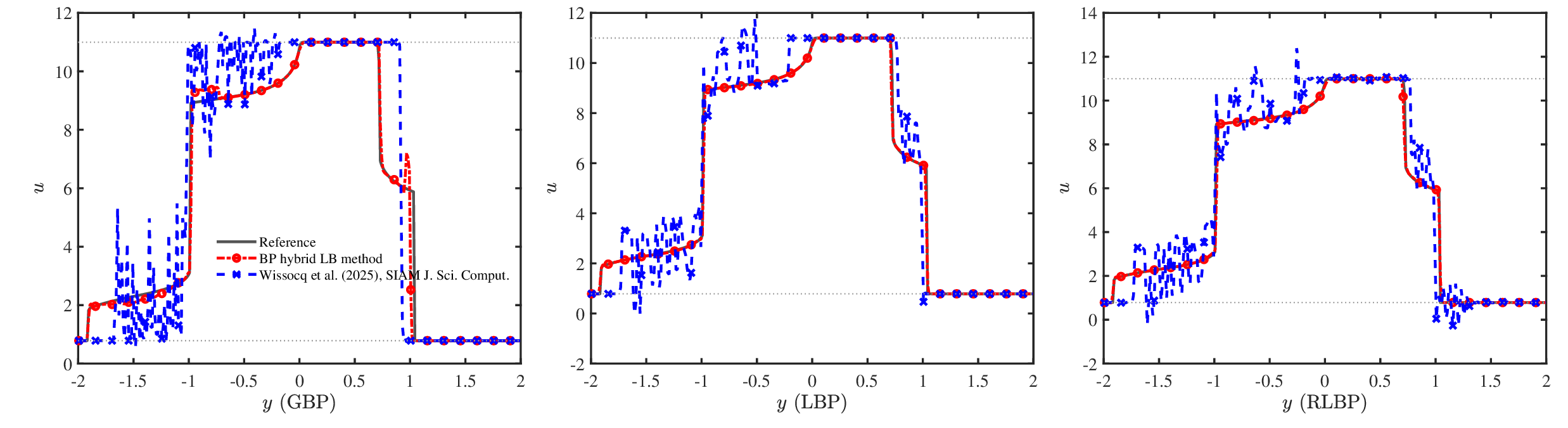}   
		\end{minipage}}   
	\caption{Profiles in the \(x\) (top) and \(y\) (bottom) directions at \(t = 1\) under \(c=1.5\). }
		\label{chapnum-ex3-c1p5-profiles}
	\end{figure}

\begin{figure}[h]
		\centering\par\vspace{-5mm}
		\subfloat[GBP limiter]
		{\begin{minipage}{0.55\linewidth}
				\centering
				\includegraphics[width=\linewidth]{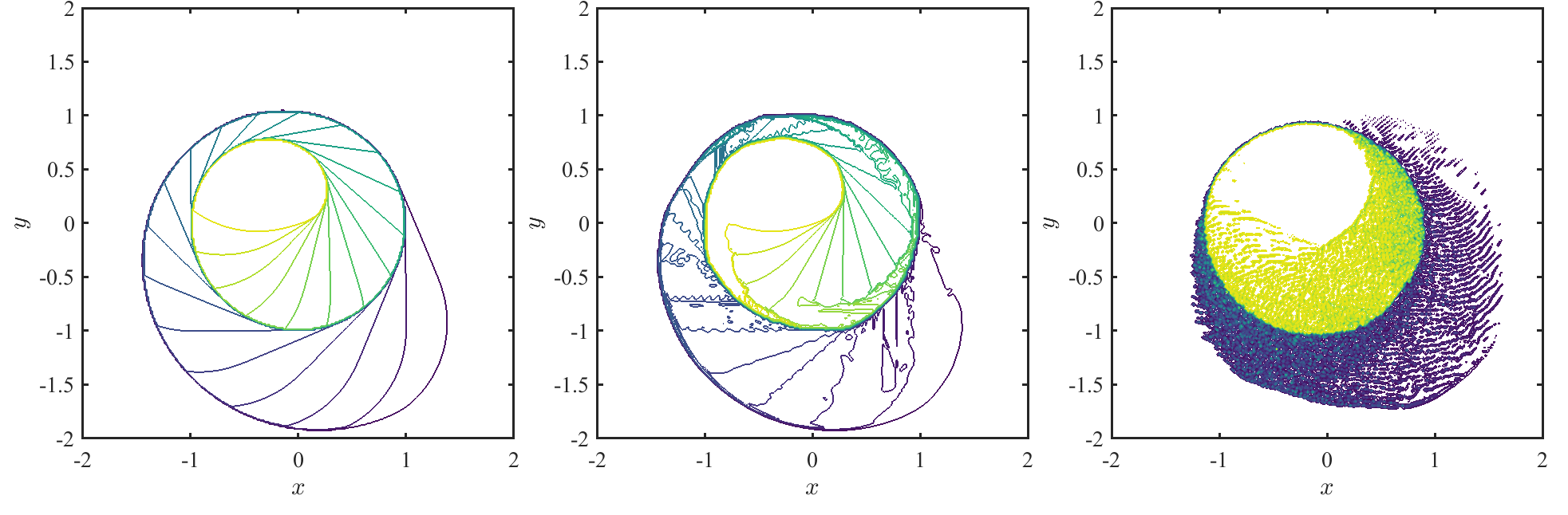}   
		\end{minipage}} 
		
		\subfloat[LBP limiter]
		{\begin{minipage}{0.55\linewidth}
				\centering
				\includegraphics[width=\linewidth]{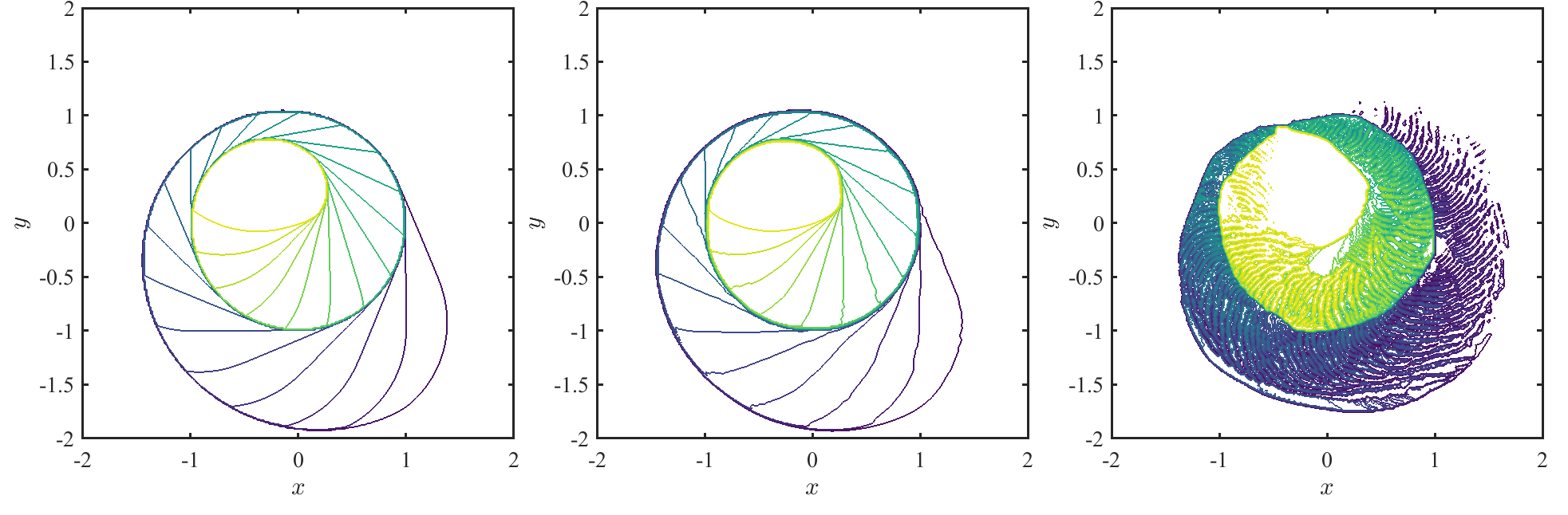}	
		\end{minipage}}  
		
		\subfloat[RLBP limiter]
		{\begin{minipage}{0.55\linewidth}
				\centering
				\includegraphics[width=\linewidth]{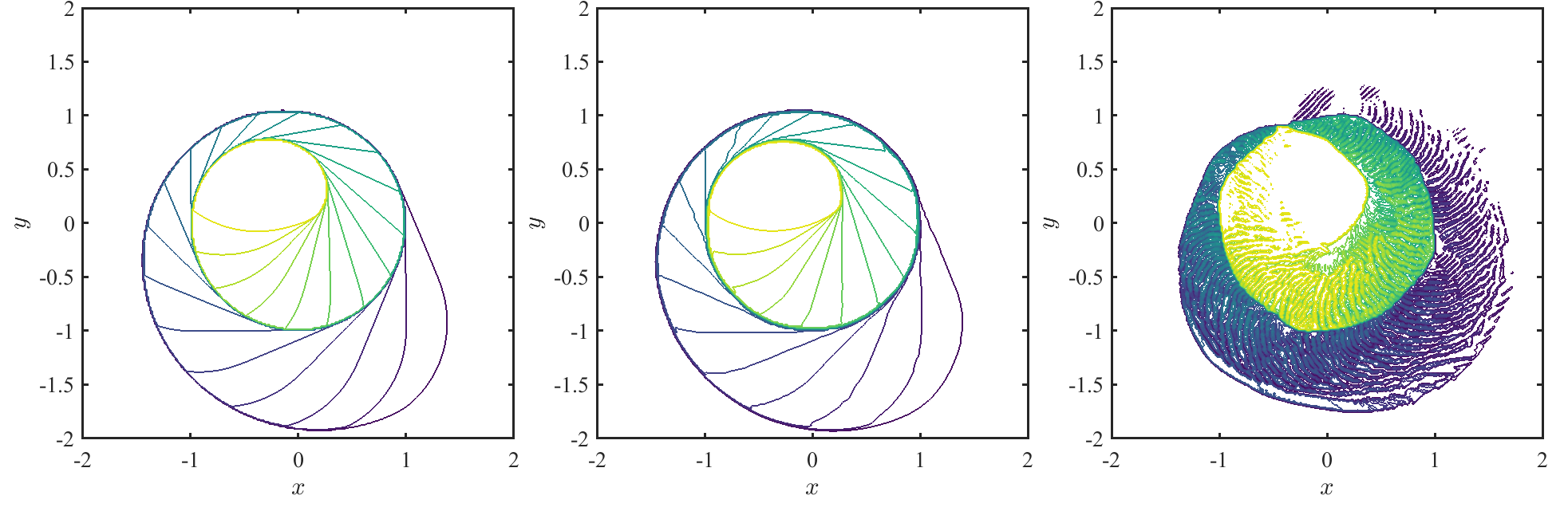}   
		\end{minipage}}    
		
		\caption{2D contour lines at \(t = 1\) under \(c = 1.5\). From left to right: the reference solution, the numerical solution of the hybrid LB method, and that of the method in Ref. \cite{WissocqLiuAbgrall2025}.}
		\label{chapnum-ex3-c1p5-contours}
	\end{figure}

	Specifically, we set \(c = 1.5 > \sqrt{2}\) for comparison in the numerical simulations, with all other experimental settings unchanged. Figures \ref{chapnum-ex3-c1p5-profiles} and \ref{chapnum-ex3-c1p5-contours} present the profiles and contours at \(t = 1.0\) under the three convex limiters. The results show that the present hybrid LB method still maintains the prescribed upper and lower bounds and captures the main structure of the rotating wave. In contrast, the method in Ref.  \cite{WissocqLiuAbgrall2025} fails to maintain numerical stability under this parameter setting, as its stability condition requires \(c \geq 2\). This indicates that, for this test case, the proposed method preserves the boundedness of macroscopic quantity while allowing for larger time steps.

\begin{table}[H]   
	\small
	\begin{center}   
		\caption{\(L_1\), \(L_2\), and \(L_\infty\) errors of the present hybrid LB method and the method of Wissocq et al. \cite{WissocqLiuAbgrall2025}.}  	 
		\label{chapnum-ex3-table-errors}
		\begin{tabular}{clcccccccc}\toprule
			\multirow{2}{*}{\(t\)}&\multirow{2}{*}{Limiter}&\multicolumn{2}{c}{\(\text{Err}_1(u)\)}&&\multicolumn{2}{c}{\(\text{Err}_2(u)\)}&&\multicolumn{2}{c}{\(\text{Err}_{\infty}(u)\)}\\
			\cline{3-4}\cline{6-7}\cline{9-10}
			&&Present&Wissocq et al. \cite{WissocqLiuAbgrall2025}&&Present&Wissocq et al. \cite{WissocqLiuAbgrall2025}&&Present&Wissocq et al. \cite{WissocqLiuAbgrall2025}\\
			\midrule
			\multirow{3}{*}{0.1}&GBP&0.01721339&0.04870854&&0.17264241&0.33521732&&5.17821615&7.38069626\\
			&LBP&0.01823612&0.02477874&&0.15837634&0.23682737&&4.31692445&5.92826249 \\
			&RLBP&0.01920009&0.02548645&&0.16277402&0.23938444&&4.38483043&5.93120560\\
			\hline
			\multirow{3}{*}{0.2}&GBP&0.02984293&0.10474751&&0.25331445&0.53916328&&	5.63827001&7.44895678\\
			&LBP&0.01934161&0.03920917&&0.15933238&0.34049615&&4.06491265&6.11051770\\	
			&RLBP&0.02047530&0.04078358&&0.16588342&0.34406827&&4.30349204&6.27917502\\\hline
			\multirow{3}{*}{0.5}&GBP&0.06978324&0.29720806&&0.41794320&0.99197038&&6.37206013&	7.80360300	\\
			&LBP&0.02137018&0.06871254&&0.16755532&0.49083587&&4.02156084&6.21236511\\
			&RLBP&0.02352954&0.07438920&&0.17351770&0.50160461&&4.90414487&6.70815507\\\hline
			\multirow{3}{*}{1.0}&GBP&0.11275271&0.61510555&&0.56746395&1.46435413&&6.76845258&8.17641260\\
			&LBP&0.02251768&0.10734206&&0.16682238&0.62656382&&4.13101930&6.07316572\\
			&RLBP&0.02433739&0.11549362&&0.18137969&0.65065248&&5.44559279&6.71146761\\\hline
			\multirow{3}{*}{2.0}&GBP&0.14740929&1.20037240&&0.67437720&2.10771582&&7.20921704&8.34466350\\
			&LBP&0.02326158&0.14447500&&0.17071886&0.64720491&&5.02509093&5.57706816\\
			&RLBP&0.02468825&0.15977157&&0.18419654&0.68419857&&5.04792179&5.60421719\\
			\bottomrule
		\end{tabular}
	\end{center}  
\end{table} 
	
	\section{Conclusion}\label{chap5}
In this paper, we have constructed a BP hybrid LB method for the \(d\)-dimensional scalar hyperbolic conservation law \eqref{chapintro-eq}. It builds on a BP first-order LB method with the upwind-type equilibrium distribution functions \eqref{chapfirst-upwind-feq}, and incorporates a second-order MRT-LB method with an auxiliary source term \eqref{chapsecond-mrt-lbm} through convex limiters. The resulting hybrid LB method  \eqref{chap3-collision}-\eqref{chap3-stream}
maintains the collision-streaming structure, satisfies the discrete mass conservation, and preserves the prescribed macroscopic bounds under the corresponding parameter conditions. Within this framework, we have presented the GBP, LBP, and RLBP limiters together with their macroscopic BP proofs.

Numerical experiments test the convergence order,  interface capturing performance, and discontinuity resolution capability of the proposed method. Compared with purely first-order and second-order LB methods, the hybrid LB method controls out-of-bounds behavior while reducing numerical dissipation, and captures discontinuities and interface structures well. A further comparison of the different limiters indicates that, although the GBP constraints can guarantee global bound preservation, non-physical oscillations may still appear within the interval, while the LBP and RLBP constraints can effectively suppress such oscillations.

In future work, we will further study the discrete entropy constraints of the hybrid LB method and the conditions for convergence to an entropy solution, and explore the balance among entropy stability, bound preservation, and numerical dissipation. In addition, combining positivity-preserving low-order updates with convex limiters, we will investigate the extension to hyperbolic systems such as the compressible Euler equations, as well as the compatibility conditions among conservation, physical admissibility, and higher-order accuracy.

	\section*{Acknowledgments}
	The computation is completed in the HPC Platform of Huazhong University of Science and Technology. This work was financially
	supported by the National Natural Science Foundation of China (Grant Nos. 12672330 and 12501599), the Postdoctoral Fellowship
	Program of CPSF (Grant No. GZB20250714), China Postdoctoral Science Foundation (Grant No. 2025M773077), the Open Research
	Fund of Stake Key Laboratory of Mesoscience and Engineering (Grant No. MESO-25-D04), the Interdisciplinary Research Program of
	Hust (Grant No. 2024JCYJ001) and the Fundamental Research Funds for the Central Universities of Hust (Grant No. 2025JYCXJJ054).
 \appendix
 \section{A flux-correction limiter based on nodal bound-preserving constraints}\label{app1}
 
 Unlike the interface intermediate state limiting approach, here we do not impose the bounds in \eqref{eq:faceconstraints} on each  \(B_{\alpha,K}^{\pm}\) individually. Instead, we limit the total positive and negative flux corrections to the low-order update in each node \cite{Zalesak1979}. The interface coefficient is then determined jointly by the admissible corrections from the two adjacent nodes.
 
 Assume that the monotonicity condition  \eqref{chapfirst-feq-mono-reduce} holds, so that the equilibrium distribution functions $f_i^{\mathrm{eq}}$ are nondecreasing.
 Let  \( m_K=\min_{J\in\mathcal N_K}u_J\) and \(M_K=\max_{J\in\mathcal N_K}u_J\) with \(\mathcal{N}_K = \{K,\, K\pm\mathbf{e}_1, \dots, K\pm\mathbf{e}_d\}\), we can obtain \begin{equation}\label{eq:lowlocal}
 	m_K
 	=
 	\sum_{i=0}^{2d}f_i^{\mathrm{eq}}(m_K)
 	\le u_K^{L}
 	\le
 	\sum_{i=0}^{2d}f_i^{\mathrm{eq}}(M_K)
 	=
 	M_K.
 \end{equation}
 
 Define
 \[
 d_{K,\alpha}
 =
 \left(L_{\alpha,K}
 -L_{\alpha+d,K+\mathbf e_\alpha}\right)
 -
 \left(L_{\alpha,K-\mathbf e_\alpha}
 -L_{\alpha+d,K}\right),
 \qquad
 v_{K,\alpha}=u_K-d\,d_{K,\alpha},
 \qquad \alpha=1,2,\ldots,d.
 \]
 The conservative form of the low-order update yields
 \[
 \frac1d\sum_{\alpha=1}^{d}v_{K,\alpha}
 =
 u_K-\sum_{\alpha=1}^{d}d_{K,\alpha}
 =
 u_K^L.
 \]
 Now, we condiser the admissible bounds for the GBP, LBP, and RLBP
 limiters as follows:
 \begin{subequations}\label{eq:appendixbounds}
 	\begin{align}
 		&\mathrm{GBP:}\quad
 		\widetilde{\ell}_K^{\text{G}}=m_0,
 		\qquad
 		\widetilde r_K^{\text{G}}=M_0,
 		\\
 		&\mathrm{LBP: }\quad
 		\widetilde{\ell}_K^{\mathrm L}
 		=
 		\max\left\{
 		m_K,\,
 		\min\{u_K,v_{K,1},\ldots,v_{K,d}\}
 		\right\},
 		\qquad
 		\widetilde r_K^{\mathrm L}
 		=
 		\min\left\{
 		M_K,\,
 		\max\{u_K,v_{K,1},\ldots,v_{K,d}\}
 		\right\},
 		\\
 		&\mathrm{RLBP: }\quad
 		\widetilde{\ell}_K^{\mathrm R}
 		=
 		\widetilde{\ell}_K^{\mathrm L}
 		-
 		C_{\mathrm{relax}}
 		\left|\widetilde{\ell}_K^{\mathrm L}\right|,
 		\qquad \widetilde r_K^{\mathrm R}
 		=
 		\widetilde r_K^{\mathrm L}
 		+
 		C_{\mathrm{relax}}
 		\left|\widetilde r_K^{\mathrm L}\right|,
 		\qquad C_{\mathrm{relax}}\ge0.
 	\end{align}
 \end{subequations}
 Since \(u_K^L\) is the arithmetic mean of \(v_{K,1},v_{K,2}, \dots, v_{K,d}\), it lies between their minimum and maximum. Together with \eqref{eq:lowlocal}, this gives \(m_K \le \widetilde{\ell}_K^{\mathrm{L}} \le u_K^L \le \widetilde{r}_K^{\mathrm{L}} \le M_K\). In particular, the RLBP interval contains the LBP interval. For the GBP limiter, the inclusion follows from \eqref{eq:lowlocal} provided the current solution \(u_K^L\) lies in \([m_0, M_0]\). 

 For the hybrid LB method in the conservation form of a numerical flux \eqref{eq:uk-update}, which can be rewittrn as
 \begin{align*}
 	u_K^{n+1}
 	=
 	u_K^L+S_K^+-S_K^-,
 \end{align*}
 where 
 \begin{align*}
 	S_K^+
 	=
 	\sum_{\alpha=1}^{d}
 	\left(
 	\eta_\alpha^-
 	[A_{e_\alpha^-}]_+
 	+
 	\eta_\alpha^+
 	[-A_{e_\alpha^+}]_+
 	\right),
 	\qquad	S_K^-
 	=
 	\sum_{\alpha=1}^{d}
 	\left(
 	\eta_\alpha^-
 	[-A_{e_\alpha^-}]_+
 	+
 	\eta_\alpha^+
 	[A_{e_\alpha^+}]_+
 	\right).
 \end{align*}
  On this basis, we define
\[	P_K^+
 	=
 	\sum_{\alpha=1}^{d}
 	\left(
 	[A_{e_\alpha^-}]_+
 	+
 	[-A_{e_\alpha^+}]_+
 	\right),
 	\qquad
 	P_K^-
 =
 	\sum_{\alpha=1}^{d}
 	\left(
 	[-A_{e_\alpha^-}]_+
 	+
 	[A_{e_\alpha^+}]_+
 	\right),
 	\qquad 	Q_K^+
 	=\overline r_K-u_K^L,
 	\qquad
 	Q_K^-=u_K^L-\overline{\ell}_K.
 \]
Here the bound-pair \((\overline{\ell}_K, \overline{r}_K)\) is \((\widetilde{\ell}_K^{\text{G}}, \widetilde{r}_K^{\text{G}})\) for the GBP limiter, \((\widetilde{\ell}_K^{\mathrm{L}}, \widetilde{r}_K^{\mathrm{L}})\) for the LBP limiter, and \((\widetilde{\ell}_K^{\mathrm{R}}, \widetilde{r}_K^{\mathrm{R}})\) for the RLBP limiter. We then set
 \begin{equation}\label{eq:cell-ratios}
 	R_K^\pm
 	=
 	\begin{cases}
 		\displaystyle
 		\min\left\{1,\frac{Q_K^\pm}{P_K^\pm}\right\},
 		& P_K^\pm>0,
 		\\[2mm]
 		1, & P_K^\pm=0.
 	\end{cases}
 \end{equation}
 For the interfaces
 $e_\alpha^\pm=K\pm\mathbf e_\alpha/2$ (\(\alpha=1,2,\ldots,d\)) with $J=K\pm\mathbf e_\alpha$, we choose
 \begin{equation}\label{eq:cell-eta}
 	\eta_\alpha^+
 	=
 	\begin{cases}
 		1, & A_{e_\alpha^+}=0,\\
 		\min\{R_K^+,R_J^-\}, & A_{e_\alpha^+}<0,
 		\\
 		\min\{R_K^-,R_J^+\}, & A_{e_\alpha^+}>0,
 	\end{cases}
\qquad 
 \eta_\alpha^-
 =
 \begin{cases}
 	1, & A_{e_\alpha^-}=0\\
 	\min\{R_{J}^+,R_K^-\},
 	& A_{e_\alpha^-}<0,\\
 	\min\{R_{J}^-,R_K^+\},
 	& A_{e_\alpha^-}>0.
 \end{cases}
 \end{equation}
 It follows directly from \eqref{eq:cell-eta} that all interface coefficients lie in \([0,1]\). Moreover, \eqref{eq:cell-eta} gives
 \[
 0\le S_K^\pm\le R_K^\pm P_K^\pm\le Q_K^\pm.
 \]
  Consequently,
 \[
 \overline{\ell}_K
 =
 u_K^L-Q_K^-
 \le
 u_K^L-S_K^-
 \le
 u_K^{n+1}
 =
 u_K^L+S_K^+-S_K^-
 \le
 u_K^L+S_K^+
 \le
 u_K^L+Q_K^+
 =
 \overline r_K.
 \]
Thus, the GBP limiter preserves the initial global bounds, the LBP limiter satisfies the local maximum principle, and the RLBP limiter preserves the prescribed relaxed local bounds. And in general, the RLBP limiter does not imply preservation of the initial global bounds. We further note that the flux corrections based on nodal constraints inherit the BP conditions of the lower-order LB method. However, treating the positive and negative corrections separately may suppress the contribution of the higher-order LB method.

 \bibliographystyle{elsarticle-harv} 
 \bibliography{reference}
 
\end{document}